\documentclass[11pt]{amsart}
 \usepackage{amsbsy,amssymb,amscd,amsfonts,latexsym,amstext,delarray,
 amsmath,graphicx}
\input{amssymb.sty}
\usepackage[dvips]{epsfig}

\input xypic

\newtheorem{thm}{Theorem}[section]
\newtheorem{prop}[thm]{Proposition}
\newtheorem{cor}[thm]{Corollary}
\newtheorem{lem}[thm]{Lemma}
\newtheorem{defn}[thm]{Definition}
\newtheorem{rem}[thm]{Remark}
\newtheorem{ex}[thm]{Example}

\numberwithin{equation}{section}

\def\bC{{\mathbb C}}

\def\bG{{\mathbb G}}
\def\bH{{\mathbb H}}

\def\bP{{\mathbb P}}

\def\bS{{\mathbb S}}
\def\bT{{\mathbb T}}

\def\A{{\mathbb A}}
\def\C{{\mathbb C}}
\def\R{{\mathbb R}}
\def\F{{\mathbb F}}
\renewcommand{\H}{{\mathbb H}}
\def\K{{\mathbb K}}
\def\N{{\mathbb N}}
\renewcommand{\P}{{\mathbb P}}
\def\Q{{\mathbb Q}}
\def\Z{{\mathbb Z}}
\def\Zb{{\mathbb Z}}

\def\cal{\mathcal}

\def\Aut{{\mbox{Aut}}}
\def\dim{{\mbox{dim}}}

\def\Hom {{\mbox{Hom}}}

\def\End{{\mbox{End}}}
\def\Trace{{\mbox{Trace}}}

\def\wh{\widehat}

\def\cala{{\cal A}}

\def\call{{\cal L}}

\def\frach{{\mathfrak h}}

\def\bbbone{\mbox{\rm 1\hspace {-.6em} l}}

\def\fU{{\mathfrak U}}

\def\cutint{{\int \!\!\!\!\!\! -}}
\def\qqq{{\,,\quad \forall}}
\def\cA{{\mathcal A}}
\def\cB{{\mathcal B}}
\def\cC{{\mathcal C}}

\def\cE{{\mathcal E}}

\def\cG{{\mathcal G}}
\def\cH{{\mathcal H}}

\def\cK{{\mathcal K}}
\def\cL{{\mathcal L}}
\def\cM{{\mathcal M}}

\def\cO{{\mathcal O}}

\def\cR{{\mathcal R}}
\def\cS{{\mathcal S}}
\def\cT{{\mathcal T}}
\def\cU{{\mathcal U}}

\def\cX{{\mathcal X}}
\def\cY{{\mathcal Y}}

\newcommand{\ie}{{\it i.e.\/}\ }
\newcommand{\eg}{{\it e.g.\/}\ }
\newcommand{\cf}{{\it cf.\/}\ }

\def\text{\hbox}

\def\Aut{{\rm Aut}}

\def\End{{\rm End}}

\def\Gal{{\rm Gal}}
\def\GL{{\rm GL}}

\def\Hom{{\rm Hom}}

\def\Ind{{\rm Ind}}

\def\PGL{{\rm PGL}}
\def\PSL{{\rm PSL}}

\def\Res{{\rm Res}}

\def\SL{{\rm SL}}

\def\Sp{{\rm Spec}}

\def\Trace{{\rm Tr}}
\def\Tr{{\rm Tr}}
\def\tr{{\rm tr}}

\title[The noncommutative ``Pardis'']{A walk in the noncommutative garden}
\author[Connes]{Alain Connes}
\author[Marcolli]{Matilde Marcolli}
\address{A.~Connes: Coll\`ege de France \\
3, rue d'Ulm \\ Paris, F-75005 France
\\ I.H.E.S. and Vanderbilt
University} \email{alain\@@connes.org}
\address{M.~Marcolli: Max--Planck Institut f\"ur Mathematik  \\
Vivatsgasse 7 \\
Bonn, D-53111 Germany} \email{marcolli\@@mpim-bonn.mpg.de}

\begin{document}
\maketitle

\tableofcontents

\newpage

\section{Introduction}\label{intro}

\begin{verse}
{\em If you cleave the hearth of one drop of water \\
a hundred pure oceans emerge from it.} \\
\smallskip
(Mahmud Shabistari, {\em Gulshan-i-raz})
\end{verse}

\bigskip

We have decided to contribute to the volume of the IPM lectures on
noncommutative geometry a text that collects a list of examples of
noncommutative spaces. As the quote of the Sufi poet here above suggests,
it is often better to approach a new subject by analyzing specific
examples rather than presenting the general theory. We hope that the
diversity of examples the readers will encounter in this text will
suffice to convince them of the fact that noncommutative geometry is
a very rich field in rapid evolution, full of interesting and yet
unexplored landscapes. Many of the examples collected here have not
yet been fully explored from the point of view of the general
guidelines we propose in Section \ref{road} and the main point of
this text is to provide a great number of open questions. The reader
should interpret this survey as a suggestion of possible interesting
problems to investigate, both in the settings described here, as
well as in other examples that are available but did not fit in this
list, and in the many more that still await to be discovered.
Besides the existing books on NCG such as \cite{Co90},
\cite{Manin91NC}, \cite{Co94}, \cite{khal},
\cite{landi} \cite{FGV} \cite{Mar}, two new
books are being written: one by the two authors of this paper \cite{CoMarBook},
and one by Connes and Moscovici \cite{CoMoBook}.

\section{Handling noncommutative spaces in the wild:
basic tools}\label{road}

We are going to see in many example how one obtains the algebra
of coordinates $\cA$ of a noncommutative space $X$.
Here we think of $\cA$ as being the algebra
of ``smooth functions'', which will usually be
a dense subalgebra of a $C^*$-algebra $\bar\cA$.

\smallskip

Here are some basic steps that one can perform in order to
acquire a good understanding of a given noncommutative space $X$
with algebra of coordinates $\cA$.

\begin{itemize}

\item[1)] Resolve the diagonal of $\cA$ and compute the cyclic cohomology.

\item[2)] Find a geometric model of $X$ up to homotopy.

\item[3)] Construct the spectral geometry $(\cA,\cH,D)$.

\item[4)] Compute the time evolution and analyze the thermodynamics.

\end{itemize}

\medskip

1) The first step means finding a resolution of the $\cA$-bimodule
$\cA$ by projective $\cA$-bimodules making it possible to compute
the Hochschild homology of $\cA$ effectively. In general, such
resolutions will be of Kozsul type and the typical example is the
resolution of the diagonal for the algebra $C^\infty(X)$ of smooth
functions on a compact manifold as in the $C^\infty$ version
\cite{ConnesCH} of the Hochschild, Kostant, Rosenberg theorem (\cf
\cite{HKR}). It makes it possible to know what is the analogue of
differential forms and of de Rham currents on the space $X$ and to
take the next step of computing the cyclic homology and cyclic
cohomology of $\cA$, which are the natural replacements for the de
Rham theory. For foliation algebras
this was done long ago (\cf \cite{Co1}, \cite{BryNi}, \cite{Crai}).
It ties in with the natural double complex of transverse currents.

It is not always easy to perform this step of finding
a resolution and computing Hochschild and cyclic (co)homology.
For instance, in the case of algebras given by generators and 
relations this uses the whole theory of Kozsul duality,
which has been successfully extended to $N$-homogeneous algebras (\cf 
\cite{D-V}, \cite{D-V1}, \cite{D-V2}, \cite{D-V3}, \cite{D-V4},
\cite{BDV}). 

One specific example in which it would be very interesting to
resolve the diagonal is the modular Hecke algebras (Section
\ref{hecke}). In essence, finding a resolution of the diagonal in
the algebra of modular forms of arbitrary level, equivariant with
respect to the action of the group $\GL_2(\A_f)$ of finite adeles,
would yield formulas for the compatibility of Hecke operators with
the algebra structure. This is a basic and hard problem of the
theory of modular forms.

\medskip

Cyclic cohomology (and homology) is a well developed theory which
was first designed to handle the leaf spaces of foliations as well
as group rings of discrete groups (\cf \cite{khal}). The theory
admits a purely algebraic version which is at center stage in
``algebraic" noncommutative geometry, but it is crucial in the
analytic set-up to construct cyclic cocycles with good compatibility
properties with the topology of the algebra. For instance, when the
domain of definition of the cocycle is a dense subalgebra stable
under holomorphic functional calculus, it automatically gives an
invariant of the $K$-theory of the underlying $C^*$-algebra (\cf
\cite{Co1}).

\smallskip

2) The essence of the second step is that many noncommutative spaces
defined as a ``bad quotients'' (\cf Section \ref{quotients})
can be desingularized, provided one is
ready to work up to homotopy. Thus for instance if the space $X$ is
defined as the quotient
$$
X=\,Y/\,\sim
$$
of an ``ordinary" space $Y$ by an equivalence relation $\sim$ one
can often find a description of the same space $X$ as a quotient
$$
X=\,Z/\,\sim
$$
where the equivalence classes are now {\em contractible} spaces. The
homotopy type of $Z$ is then uniquely determined and serves as a
substitute for that of $X$ (see \cite{BauCo}).

For instance, if the equivalence relation on $Y$ comes from the free
action of a torsion free discrete group $\Gamma$, the space $Z$ is
simply a product over $\Gamma$ of the form
$$
Z=\,Y\times_\Gamma\,E\Gamma ,
$$
where $E\Gamma$ is a contractible space on which $\Gamma$ acts
freely and properly.

The main point of this second step is that it gives a starting point
for computing the $K$-theory of the space $X$ \ie of the
$C^*$-algebra $A=\bar\cA$ playing the role of the algebra of continuous
functions on $X$. Indeed, for each element of the $K$-homology
of the classifying space $Z$, there is a general construction
of an index problem for ``families
parameterized by $X$" that yields an assembly map (\cf \cite{BauCo})
\begin{equation}\label{mu-map}
\mu\;:\; K_*(Z)\to K_*(A)
\end{equation}
This Baum--Connes map is an isomorphism in a lot of cases (with
suitable care of torsion, \cf \cite{BCH}) including all
connected locally compact groups, all amenable groupoids
and all hyperbolic discrete groups.
It thus gives a computable guess for the $K$-theory of $X$.

The next
step is not only to really compute $K(A)$ but also to get a good
model for the ``vector bundles" on $X$ \ie the finite projective
modules over $A$. This step should then be combined with the above
first step to compute the Chern character using connections,
curvature, and eventually computing moduli spaces of Yang-Mills
connections as was done for instance for the NC-torus in
\cite{CoRi}.

\smallskip

3) The third step makes it possible to pass from the soft part of
differential
geometry to the harder ``Riemannian" metric aspect. The sought for
spectral geometry $(\cA,\cH,D)$ has three essential features:

\begin{itemize}

\item The $K$-homology class of $(\cA,\cH,D)$.

\item The smooth structure.

\item The metric.

\end{itemize}

One should always look for a spectral triple whose $K$-homology
class is as non-trivial as possible. Ideally it should extend to a
class for the double algebra $\cA\otimes \cA^o$ and then be a generator
for Poincar\'e duality. In general this is too much to ask for, since
many interesting spaces do not fulfill Poincar\'e duality. The main
tool for determining the stable homotopy class of the spectral
triple is Kasparov's bivariant $KK$-theory. Thus it is quite
important to already have taken step 2 and to look for classes whose
pairing with $K$-theory is as non-trivial as can be. For the smooth
structure, there is often a natural guess for a subalgebra
$A^\infty\subset A$ of the $C^*$-algebra $A=\bar\cA$ that will play the role
of the algebra of smooth functions. It should in general contain the
original algebra $\cA$ but should have the further property of being
stable under the holomorphic functional calculus. This ensures that
the inclusion $A^\infty\subset A$ is an isomorphism in $K$-theory
and makes it possible to complete the classification of smooth vector bundles.

The role of the unbounded operator $D$ for the smooth structure is
that it defines the geodesic flow by the formula
$$
F_t(a)=\,e^{it|D|}\,a\,e^{-it|D|}\qqq a\in A^\infty
$$
and one expects that smoothness is governed by the smoothness of the
operator valued map $\R\ni t \mapsto F_t(a)$. The main result of the
general theory is the local index formula of \cite{CoMo}, which
provides the analogue of the Pontrjagin classes of smooth manifolds
in the noncommutative framework.

The problem of determining $D$ from the knowledge of the
$K$-homology class is very similar to the choice of a connection on
a bundle. There are general results that assert the existence of an
unbounded selfadjoint $D$ with bounded commutators with $\cA$ from
estimates on the commutators with the phase $F$. The strongest is obtained
(\cf \cite{Co94}) just assuming that the $[F,a]$ are in an ideal
called ${\rm Li}\,\cH$ and it ensures the existence of a theta-summable
spectral triple which is what one needs to get started.

It is not always possible to find a finitely-summable spectral
triple, first because of growth conditions on the algebra
\cite{Co-fr}, but also since the finitely-summable condition is very
analogous to  type II in the theory of factors. In very general
cases, like the noncommutative space coming from foliations, one can
however go from type III to type II by passing to the total space of
the space of transverse metrics and then use the theory of
hypoelliptic operators \cite{CoMoHopf}.

\medskip

Another way to attack the problem of determining $D$ is to consider
the larger algebra generated by $\cA$ and $D$, write a-priori
relations between $\cA$ and $D$ and then look for irreducible
representations that fall in the correct stable homotopy class.
Ideally one should minimize the spectral action functional
\cite{C-C2} in this homotopy class thus coming close to gravity. In
practise one should use anything available and the example of the
NC-space given by the quantum group $SU_q(2)$ shows that things can
be quite subtle \cite{SDLSV2}.

Once the spectral triple $(\cA,\cH,D)$ has been determined, the basic
steps are the following, one should compute

\begin{itemize}

\item The dimension spectrum $\Sigma\subset \C$.

\item The local index formula.

\item The inner fluctuations, scalar curvature, and spectral action.

\end{itemize}

\bigskip

4) Often a noncommutative space comes with a measure class, which in
turn determines a time evolution $\sigma_t$, namely a 1-parameter
family of automorphisms of the $C^*$-algebra $A=\bar\cA$. In the type II
situation one can apply the discussion of step 3 above and, in the
finite dimensional case, use the operator $D$ to represent
functionals in the measure class in the form
$$
\varphi(a)=\, \int\!\!\!\!\!\!-\ \,a \,|D|^{-p} \qqq \,a\in \cA
$$
where $\int\!\!\!\!\!\!-\  $ is the noncommutative integral \ie the
Dixmier trace and $p$ is the ``dimension".

In the general case one should expect to be in the type III
situation in which  the time evolution $\sigma_t$ is highly
non-trivial. We shall see some examples, for instance in Section
\ref{Qlatt}. Given the data $(\bar\cA,\sigma_t)$ it is natural to regard
it as a quantum statistical mechanical system, with $\bar\cA$ as algebra
of observables and $\sigma_t$ as time evolution. One can then look
for equilibrium states for the system, for a given value $\beta$ of
the thermodynamic parameter (inverse temperature).

If the algebra $\bar\cA$ is concretely realized as an algebra of bounded
operators on a Hilbert space $\cH$, then one can consider the
Hamiltonian $H$, namely the (unbounded) operator on $\cH$ that is the
infinitesimal generator of the time evolution. If the operator
$\exp(-\beta H)$ is trace class, then one has equilibrium states for
the system $(\bar\cA,\sigma_t)$ written in the usual Gibbs form
$$
\varphi_\beta(a)=\frac{\Tr(a\exp(-\beta H))}{\Tr(\exp(-\beta H))},
$$
where $Z(\beta)=\Tr(\exp(-\beta H))$ is the partition function of the
system. The notion of equilibrium state continues to make sense
when $\exp(-\beta H)$ is not necessarily trace class, and is given by
the more subtle notion of KMS (Kubo--Martin--Schwinger) states.

These are states on $\bar\cA$, namely continuous functionals
$\varphi: \bar\cA\to \C$ with $\varphi(1)=1$ and $\varphi(a^*a)\geq
0$, satisfying the KMS$_\beta$ condition that, for all $a,b\in
\bar\cA$ there exists a function $F_{a,b}(z)$ which is holomorphic
on the strip $0< \Re(z) <i\beta$ continuous and bounded on the
closed strip and such that, for all $t\in\R$,
\begin{equation}\label{KMScond}
  F_{a,b}(t)=\varphi(a\sigma_t(b)) \ \ \text { and } \ \
F_{a,b}(t+i\beta)=\varphi(\sigma_t(b)a).
\end{equation}
KMS states at zero temperature can be defined as weak limits as
$\beta\to \infty$ of KMS$_\beta$ states.

\begin{figure}
\begin{center}
\includegraphics[scale=0.9]{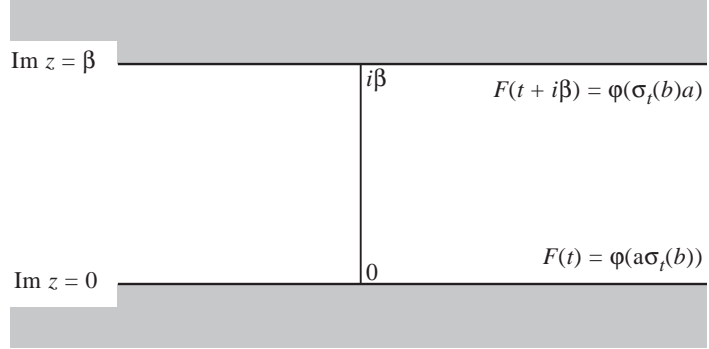}
\end{center}
\caption{The KMS condition \label{FigKMS}}
\end{figure}

One can construct using KMS states very refined invariants of
noncommutative spaces. For a fixed $\beta$, the KMS$_\beta$ states
form a simplex, hence one can consider only the extremal KMS$_\beta$
states $\cE_\beta$, from which one recovers all the others by convex
combinations. An extremal KMS$_\beta$ state is always factorial and
the type of the factor is an invariant of the state. The simplest
situation is the type I. One can show under minimal hypothesis
(\cite{CCM}) that extremal KMS$_\beta$ states continue to survive
when one lowers the temperature \ie one increases $\beta$. Thus, in
essence, when cooling down the system this tends to become more and
more ``classical" and in the $0$-temperature limit $\cE_\beta$ gives
a good replacement of the notion of classical points for a
noncommutative space. We shall see in Section \ref{shimura} how, in
examples related to arithmetic, the ``classical points'' described
by the zero temperature KMS states of certain quantum statistical
mechanical systems recover classical arithmetic varieties. The
extremal KMS states at zero temperature, evaluated on suitable
arithmetic elements in the noncommutative algebra, can be shown in
significant cases to have an interesting Galois action, related to
interesting questions in number theory (\cf \cite{BC}, \cite{CMR},
\cite{CoMa}).

In joint work with Consani, we showed in \cite{CCM}
how to define an analog in characteristic
zero of the action of the Frobenius on the etale cohomology by a
process involving the above thermodynamics. One key feature is that
the analogue of the Frobenius is the ``dual" of the above time
evolution $\sigma_t$. The process involves cyclic homology and its
three basic steps are (\cite{CCM})

\begin{itemize}

\item Cooling.

\item Distillation.

\item Dual action of  $\R_+^*$ on the cyclic homology of the
distilled space.

\end{itemize}

\medskip
When applied to the simplest system
(the Bost--Connes system of \cite{BC})
this yields a cohomological interpretation of the spectral
realization of the zeros of the Riemann zeta function
(\cite{CoRH}, \cite{CCM}).

\bigskip

\section{Phase spaces of microscopic systems}\label{heisenberg}

What can be regarded historically as the first example of a
non--commutative space is the Heisenberg formulation of the
observational Ritz-Rydberg law of spectrocopy. In fact, quantum
mechanics showed that indeed the parameter space, or phase space
of the mechanical system given by a single atom fails to be a
manifold. It is important to convince oneself of this fact and to
understand that this conclusion is indeed dictated by the {\it
experimental findings of spectroscopy}.

At the beginning of the twentieth century a wealth of experimental
data was being collected on the spectra of various chemical
elements. These spectra obey experimentally discovered laws, the
most notable being the Ritz-Rydberg combination principle. The
principle can be stated as follows; spectral lines are indexed by
pairs of labels. The statement of the principle then is that
certain pairs of spectral lines, when expressed in terms of
frequencies, do add up to give another line in the spectrum.
Moreover, this happens  precisely when the labels are of the form
$i,j$ and $j,k$.

\begin{center}
\begin{figure}
\includegraphics[scale=0.9]{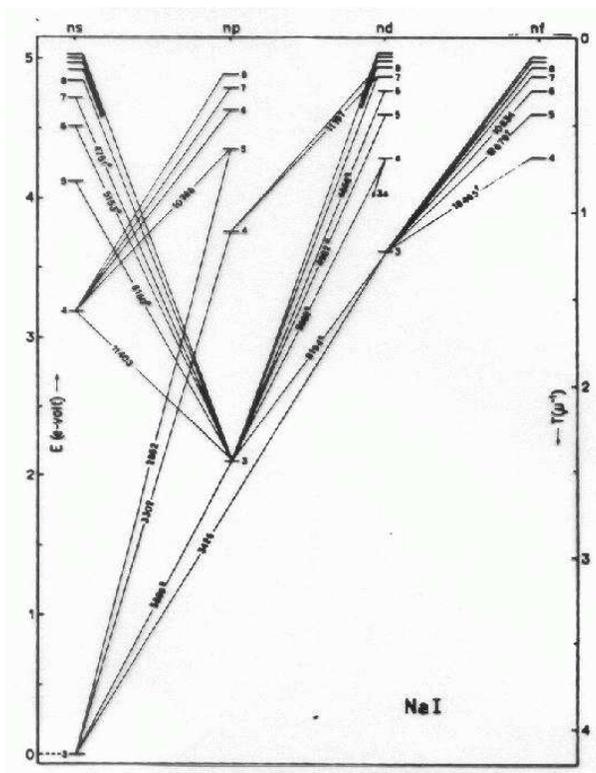}
\caption{Spectral lines and the Ritz-Rydberg law.
\label{FigRitz}}
\end{figure}
\end{center}

In the seminal paper \cite{Heisenberg} of 1925, Heisenberg
considers the classical prediction for the radiation emitted by a
moving electron in a field, where the observable dipole moment can
be computed, with the motion of the electron given in Fourier
expansion. The classical model would predict (in his notation)
frequencies distributed according to the law
\begin{equation}\label{emissionCM}
 \nu(n,\alpha) = \alpha \nu(n) = \alpha \frac{1}{h} \frac{dW}{dn}.
\end{equation}
When comparing the frequencies obtained in this classical model
with the data, Heisenberg noticed that the classical law
\eqref{emissionCM} did not match the phenomenon observed.

The spectral rays provide a `picture' of an atom: if atoms were
classical systems, then the picture formed by the spectral lines
would be (in our modern mathematical language) a {\it group},
which is what \eqref{emissionCM} predicts. That is, the classical
model predicts that the observed frequencies should simply add,
obeying a group law, or, in Heisenberg's notation, that
\begin{equation}\label{emission2}
 \nu(n,\alpha)+ \nu(n,\beta)= \nu(n,\alpha+\beta).
\end{equation}
Correspondingly, observables would form the convolution algebra of
a group.

What the spectral lines were instead providing was the picture of
a {\it groupoid}. Heisenberg realized that the classical law of
(\ref{emissionCM}) \eqref{emission2} would have to be replaced by
the quantum--mechanical
\begin{equation}\label{emissionQM}
 \nu(n,n-\alpha) = \frac{1}{h} \left( W(n)- W(n-\alpha) \right).
\end{equation}
This replaces the group law with that of a groupoid, replacing the
classical \eqref{emission2} by the quantum--mechanical
\begin{equation}\label{emissionQ}
 \nu(n,n-\alpha)+ \nu(n-\alpha, n-\alpha-\beta)=
\nu(n,n-\alpha-\beta).
\end{equation}
Similarly, the classical Fourier modes $\fU_{\alpha}(n)
e^{i\omega(n)\alpha t}$ were replaced by $ \fU
(n,n-\alpha)e^{i\omega(n,n-\alpha)t}$.

The analysis of the emission spectrum given by Heisenberg was in
very good agreement with the Ritz--Rydberg law, or combination
principle, for spectral lines in emission or absorption spectra.

In the same paper, Heisenberg also extends his redefinition of the
multiplication law for the Fourier coefficients to coordinates and
momenta, by introducing transition amplitudes that satisfy similar
product rules. This is, in Born's words, ``his most audacious
step'': in fact, it is precisely this step that brings
non--commutative geometry on the scene.

It was Born who realized that what Heisenberg described in his
paper corresponded to replacing classical coordinates with
coordinates which no longer commute, but which obey the laws of
matrix multiplication. In his own words reported in
\cite{vdWaerden},

\begin{quote}
After having sent Heisenberg's paper to the Zeitschrift f\"ur
Physik for publication, I began to ponder about his symbolic
multiplication, and was soon so involved in it that I thought the
whole day and could hardly sleep at night. For I felt there was
something fundamental behind it ... And one morning ... I suddenly
saw light: Heisenberg's symbolic multiplication was nothing but
the matrix calculus.
\end{quote}

\bigskip

Thus, spectral lines are parameterized by two indices
$L_{\alpha\beta}$ satisfying a cocycle relation
$L_{\alpha\beta}+L_{\beta\gamma}= L_{\alpha\gamma}$, and a
coboundary relation expresses each spectral line as a difference
$L_{\alpha\beta}=\nu_\alpha -\nu_\beta$. In other words, the
Ritz--Rydberg law gives the groupoid law \eqref{emissionQ}, or
equivalently,
$$ (i,j)\bullet (j,k)= (i,k) $$
and the convolution algebra of the group is replaced by
observables satisfying the matrix product
$$ (AB)_{ik} = \sum_j A_{ij} B_{jk}, $$
for which in general commutativity is lost:
$$ AB\neq BA. $$
The Hamiltonian $H$ is a matrix with the frequencies on the
diagonal, and observables obey the evolution equation
$$ \frac{d}{dt} A = i [H,A]. $$

Out of Heisenberg's paper and Born's interpretation of the same in
terms of matrix calculus, emerged the statement of Heisenberg's
uncertainty principle in the form of a commutation relation of
matrices
$$ [P,Q] = \frac{h}{2\pi i} I. $$
The matrix calculus and the uncertainty principle were formulated
in the subsequent paper of Born and Jordan \cite{Born}, also
published in 1925. This viewpoint on quantum mechanics was later
somewhat obscured by the advent of the Schr\"odinger equation.
The Schr\"odinger approach shifted the emphasis back to the more
traditional technique
of solving partial differential equations, while the more modern
viewpoint of Heisenberg implied a much more serious change
of paradigm, affecting our most basic understanding of the notion of
space. Heisenberg's approach can be regarded as the historic origin of
noncommutative geometry.

\bigskip

\section{Noncommutative quotients}\label{quotients}

A large source of examples of noncommutative spaces is given by
quotients of equivalence relations. One starts by an ordinary
commutative space $X$ (\eg a smooth manifold or more generally a
locally compact Hausdorff topological space). This can be described
via its algebra of functions $C(X)$, and abelian $C^*$-algebra.
Suppose then that we are interested in taking a quotient $Y=X/\sim$
of $X$ with respect to an equivalence relation. In general, one
should not expect the quotient to be nice. Even when $X$ is a smooth
manifold, the quotient $Y$ need not even be a Hausdorff space. In
general, one would like to still be able to characterize $Y$ through
its ring of functions.  One usually defines $C(Y)$ to be functions
on $X$ that are invariant under the equivalence relation ,
\begin{equation}\label{commCY}
C(Y)= \{ f\in C(X): f(a)=f(b), \forall a\sim b \} .
\end{equation}
Clearly, for a ``bad'' equivalence relation one typically
gets this way only constant functions $C(Y)=\C$.

There is a better way to
associate to the quotient space $Y$ a ring of functions which is
nontrivial for any equivalence relation. This requires dropping the
commutativity requirement. One can then consider functions of two
variables $f_{ab}$ defined on the graph of the equivalence relation,
with a product which is no longer the commutative pointwise product,
but the noncommutative convolution product dictated by the groupoid of
the equivalence relation.
In general the elements in the algebra of functions
\begin{equation}\label{ncCY}
``C(Y)\text{''}= \{ (f_{ab}) : a\sim b \}
\end{equation}
act as bounded operators on the Hilbert space $L^2$ of the
equivalence class. This also guarantees the convergence in the
operator norm of the convolution product
$$ \sum_{a\sim b\sim c} f_{ab} g_{bc}. $$

We give a few examples to illustrate the difference between the
traditional construction and the one of noncommutative geometry.

\begin{ex} {\em
Consider the space $Y=\{ x_0, x_1\}$ with the equivalence relation
$x_0\sim x_1$. With the first point of view the algebra of
functions on the quotient is $\C$, and in the second point of view
it is $\cB =M_2(\C)$, that is
\begin{equation}
\cB = \left\{ f = \left( \begin{array}{cc} f_{aa} & f_{ab} \\
f_{ba} & f_{bb} \end{array} \right) \right\} \, . \label{twopoint}
\end{equation}
These two algebras are not the same, though in this case they are
Morita equivalent. }\end{ex}

Notice that, when one computes the spectrum of the algebra
(\ref{twopoint}), it turns out that it is composed of only one
point, so the two points $a$ and $b$ have been identified.
This first trivial example represents the typical situation where
the quotient space is ``nice'': the two constructions give Morita
equivalent algebras. In this sense, Morita
equivalent algebras are regarded as ``the same'' (or better
isomorphic) spaces in non--commutative geometry.

\begin{ex} {\em
Consider the space $Y=[0,1] \times \{ 0,1\}$ with the equivalence
relation $(x,0)\sim (x,1)$ for $x\in (0,1)$. Then in the first
viewpoint the algebra of functions is again given just by the
constant functions $\C$, but in the second case we obtain
\begin{equation}  \{ f\in C([0,1])\otimes M_2(\C) : \, f(0) \text{ and
} f(1) \text{ diagonal } \}. \label{alg1} \end{equation} }
\label{ex2} \end{ex}

In this case, these algebras are not Morita equivalent. This can be
seen by computing their $K$--theory. This means that the approach of
non--commutative spaces produces something genuinely new, as soon as
the quotient space ceases to be ``nice'' in the classical sense.

In general, the first kind of construction of functions on the
quotient space is cohomological in nature: one seeks for functions
satisfying certain equations or constraints. Usually there are
very few solutions. The second approach, instead, typically
produces a very large class of functions.

\bigskip

\section{Spaces of leaves of foliations}\label{foliations}

There is a very rich collection of examples of noncommutative spaces
given by  the leaf spaces of foliations. The connection thus
obtained between noncommutative geometry and the geometric theory of
foliations is very far reaching for instance through the role of
Gelfand-Fuchs cohomology, of the Godbillon-Vey invariant and of the
passage from type III to type II using the transverse frame bundle.
It is this class of examples  that triggered the initial development
of cyclic cohomology (\cf \cite{khal} section 4), of the local index
formula in noncommutative geometry as well as the theory of
characteristic classes for Hopf algebra actions.

\medskip

The construction of the algebra associated to a foliation is a
special case of the construction of section \ref{quotients} but both
the presence of holonomy and the case when the graph of the
foliation is non-hausdorff require special care, so we shall recall
the basic steps below.

Let $V$ be a smooth manifold and $TV$ its tangent bundle, so that
for each $x \in V$, $T_x V$ is the tangent space of $V$ at $x$. A
smooth subbundle $F$ of $TV$ is called {\it integrable} iff one of
the following equivalent conditions is satisfied:

\smallskip

\begin{enumerate}

\item[a)] Every $x \in V$ is contained in a submanifold $W$ of $V$ such that
$$
T_y (W) = F_y \qquad \forall \, y \in W \, .
$$

\smallskip

\item[b)] Every $x \in V$ is in the domain $U \subset V$ of a
submersion $p : U \to {\mathbb R}^q$ ($q = {\rm codim} \, F$) with
$$
F_y = {\rm Ker} (p_*)_y \qquad \forall \, y \in U \, .
$$

\smallskip
\item[c)] $C^{\infty} (F) = \{ X \in C^{\infty} (TV) \, , \ X_x \in
F_x \quad \forall \, x \in V \}$ is a Lie algebra.

\smallskip

\item[d)] The ideal $J(F)$ of smooth exterior differential forms which
vanish on $F$ is stable by exterior differentiation.
\end{enumerate}

\medskip

Any $1$-dimensional subbundle $F$ of $TV$ is integrable, but for
$\dim F \geqq 2$ the condition is non trivial, for instance if $P
\overset{p}{\to} B$ is a principal $H$-bundle (with compact
structure group $H$) the bundle of horizontal vectors for a given
connection is integrable iff this connection is flat.

\smallskip

A foliation of $V$ is given by an integrable subbundle $F$ of $TV$.
The leaves of the foliation $(V,F)$ are the maximal connected
submanifolds $L$ of $V$ with $T_x (L) = F_x $, $\forall \, x \in L$,
and the partition of $V$ in leaves $$V = \cup
L_{\alpha}\,,\quad\alpha \in X$$ is characterized geometrically by
its ``local triviality'': every point $x \in V$ has a neighborhood
$U$ and a system of local coordinates
$(x^j)_{j = 1 , \ldots , \dim V}$ called
{\it foliation charts}, so
that the partition of $U$ in connected components of
leaves corresponds to the partition of $${\mathbb
R}^{\dim V} = {\mathbb R}^{\dim F} \times {\mathbb R}^{\text{codim}
\, F}$$ in the parallel affine subspaces ${\mathbb R}^{\dim F}
\times {\rm pt}$. These  are the leaves of the
restriction of $F$ and are called {\it plaques}.

The set $X=V/F$ of leaves of a foliation $(V,F)$
is in most cases a noncommutative space. In other words even though
as a set it has the cardinality of the continuum it is in general
not so at the effective level and it is in general impossible to
construct a countable set of measurable functions on $V$ that form a
complete set of invariants for the equivalence relation coming from
the partition of $V$ in leaves $V = \cup L_{\alpha}$. Even in the
simple cases in which the set $X=V/F$ of leaves is classical it
helps to introduce the associated algebraic tools in order to get a
feeling for their role in the general singular case.

\smallskip

To each foliation $(V,F)$ is canonically associated a $C^*$ algebra
$C^* (V,F)$ which encodes the topology of the space of leaves. The
construction is basically the same as the general one for quotient
spaces of Section \ref{quotients}, but there are interesting nuances
coming from the presence of holonomy in the foliation context. To
take this into account one first constructs a manifold $G$, $\dim
\,G = \dim \,V + \dim \,F$, called the graph (or holonomy groupoid)
of the foliation, which refines the equivalence relation coming from
the partition of $V$ in leaves $V = \cup L_{\alpha}$. This
construction is due to Thom, Pradines and Winkelnkemper,
see \cite{Winkel}.

An element $\gamma$ of $G$ is given by two points $x = s(\gamma)$,
$y = r(\gamma)$ of $V$ together with an equivalence class of smooth
paths: $\gamma (t)\in V$, $t \in [0,1]$; $\gamma (0) = x$, $\gamma
(1) = y$, tangent to the bundle $F$ ({\it i.e.} with $\dot\gamma (t)
\in F_{\gamma (t)}$, $\forall \, t \in {\mathbb R}$) up to the
following equivalence: $\gamma_1$ and $\gamma_2$ are equivalent iff
the {\it holonomy} of the path $\gamma_2 \circ \gamma_1^{-1}$ at the
point $x$ is the {\it identity}. The graph $G$ has an obvious
composition law. For $\gamma , \gamma' \in G$, the composition
$\gamma \circ \gamma'$ makes sense if $s(\gamma) = r(\gamma')$. If
the leaf $L$ which contains both $x$ and $y$ has no holonomy, then
the class in $G$ of the path $\gamma (t)$ only depends on the pair
$(y,x)$. The condition of trivial holonomy
is generic in the topological sense of dense $G_{\delta}$'s.
In general, if one fixes $x = s(\gamma)$, the map from $G_x
= \{ \gamma , s(\gamma) = x \}$ to the leaf $L$ through $x$, given
by $\gamma \in G_x \mapsto y = r(\gamma)$, is the holonomy covering
of $L$.

Both maps $r$ and $s$ from the manifold $G$ to $V$ are smooth
submersions and the map $(r,s)$ to $V \times V$ is an immersion
whose image in $V \times V$ is the (often singular) subset
\begin{equation}\label{subset}
\{ (y,x)\in V \times V: \, \text{ $y$ and $x$ are on the same leaf}\}.
\end{equation}
In first
approximation one can think of elements of $C^* (V,F)$ as continuous
matrices $k(x,y)$, where $(x,y)$ varies in the set \eqref{subset}.
We now describe this $C^*$ algebra in all details. We
assume, for notational convenience, that the manifold $G$ is
Hausdorff. Since this fails to be the case in very interesting
examples, we   also explain briefly how to remove this
hypothesis.

\smallskip

The basic elements of $C^* (V,F)$ are smooth half densities
$f \in C_c^{\infty} (G , \Omega^{1/2})$ with
compact support on $G$. The
bundle $\Omega_G^{1/2}$ of half densities on $G$ is obtained as
follows. One first defines a line bundle $\Omega_V^{1/2}$ on V. For
$x\in V$ one lets $\Omega_x^{1/2}$ be the one dimensional complex
vector space of maps from the exterior power $\wedge^k \, F_x$, $k =
\dim F$, to ${\mathbb C}$ such that
$$
\rho \, (\lambda \, v) = \vert \lambda \vert^{1/2} \, \rho \, (v)
\qquad \forall \, v \in \wedge^k \, F_x \, , \quad \forall \,
\lambda \in {\mathbb R} \, .
$$
Then, for $\gamma \in G$, one can identify $\Omega_{\gamma}^{1/2}$ with the one
dimensional complex vector space $\Omega_y^{1/2} \otimes
\Omega_x^{1/2}$, where $\gamma : x \to y$. In other words
$$
\Omega_G^{1/2}=\, r^*(\Omega_V^{1/2})\otimes s^*(\Omega_V^{1/2})\,.
$$
 Of course the bundle $\Omega_V^{1/2}$ is trivial on $V$, and we
could choose once and for all  a trivialisation $\nu$ turning
elements of $C_c^{\infty} (G , \Omega^{1/2})$ into functions.
Let us
however stress that the use of half densities makes all the
construction completely canonical.

For $f,g \in C_c^{\infty} (G , \Omega^{1/2})$, the convolution
product $f * g$ is defined by the equality
$$
(f * g) (\gamma) = \int_{\gamma_1 \circ \gamma_2 = \gamma}
f(\gamma_1) \, g(\gamma_2) \, .
$$
This makes sense because, for fixed $\gamma : x \to y$ and fixing $v_x
\in \wedge^k \, F_x$ and $v_y \in \wedge^k \, F_y$, the product
$f(\gamma_1) \, g(\gamma_1^{-1} \gamma)$ defines a $1$-density on
$G^y = \{ \gamma_1 \in G , \, r (\gamma_1) = y \}$, which is smooth
with compact support (it vanishes if $\gamma_1 \notin$ support $f$),
and hence can be integrated over $G^y$ to give a scalar, namely $(f * g)
(\gamma)$ evaluated on $v_x , v_y$.

The $*$ operation is defined by $f^* (\gamma) =
\overline{f(\gamma^{-1})}$, {\it i.e.} if $\gamma : x \to y$ and
$v_x \in \wedge^k \, F_x$, $v_y \in \wedge^k \, F_y$ then $f^*
(\gamma)$ evaluated on $v_x , v_y$ is equal to
$\overline{f(\gamma^{-1})}$ evaluated on $v_y , v_x$. We thus get a
$*$ algebra $C_c^{\infty} (G , \Omega^{1/2})$. For each leaf $L$ of
$(V,F)$ one has a natural representation of this $*$ algebra on the
$L^2$ space of the holonomy covering $\tilde L$ of $L$. Fixing a
base point $x \in L$, one identifies $\tilde L$ with $G_x = \{
\gamma , s(\gamma) = x \}$ and defines
$$
(\pi_x (f) \, \xi) \, (\gamma) = \int_{\gamma_1 \circ \gamma_2 =
\gamma} f(\gamma_1) \, \xi (\gamma_2) \qquad \forall \, \xi \in L^2
(G_x),
$$
where $\xi$ is a square integrable half density on $G_x$. Given
$\gamma : x \to y$ one has a natural isometry of $L^2 (G_x)$ on $L^2
(G_y)$ which transforms the representation $\pi_x$ in $\pi_y$.

By definition $C^* (V,F)$ is the $C^*$ algebra completion of
$C_c^{\infty} (G , \Omega^{1/2})$ with the norm $$\Vert f \Vert =
\sup_{x \in V} \, \Vert \pi_x (f) \Vert\,.$$

Note that $C^* (V,F)$ is always norm separable and admits a natural
smooth subalgebra, namely the algebra $C_c^\infty (V,F)=C_c^{\infty} (G ,
\Omega^{1/2})$ of smooth compactly supported half densities.

If the leaf $L$ has trivial holonomy then the representation
$\pi_x$, $x \in L$, is irreducible. In general, its commutant is
generated by the action of the (discrete) holonomy group $G_x^x$ in
$L^2 (G_x)$. If the foliation comes from a submersion $p : V \to B$,
then its graph $G$ is $\{ (x,y) \in V \times V , \, p(x) = p(y)\}$
which is a submanifold of $V \times V$, and $C^* (V,F)$ is identical
with the algebra of the continuous field of Hilbert spaces $L^2
(p^{-1} \{ x \})_{x \in B}$. Thus (unless $\dim \,F = 0$) it is
isomorphic to the tensor product of $C_0 (B)$ with the elementary
$C^*$ algebra of compact operators. If the foliation comes from an
action of a Lie group $H$ in such a way that the graph is identical
with $V \times H$ (this is not always true even for flows)  then
$C^* (V,F)$ is identical with the reduced crossed product of $C_0
(V)$ by $H$. Moreover the construction of $C^* (V,F)$ is local in
the following sense.

If $V' \subset V$ is an open set and $F'$ is the restriction of $F$
to $V'$, then the graph $G'$ of $(V' , F')$ is an open set in the
graph $G$ of $(V,F)$, and the inclusion $C_c^{\infty} (G' ,
\Omega^{1/2}) \subset C_c^{\infty} (G , \Omega^{1/2})$ extends to an
isometric $*$ homomorphism of $C^* (V' , F')$ in $C^* (V,F)$. The
proof is straightforward and also applies in the case of non
Hausdorff graph.

Let us now briefly explain how the construction of the $C^*$ algebra
$C^* (V,F)$ has to be done in the case when the graph of the
foliation is not Hausdorff. This case is rather rare, since it never
occurs if the foliation is real analytic. However, it does occur in
cases which are topologically interesting for foliations, such as
the Reeb foliation of the 3 sphere, which are constructed by
patching together foliations of manifolds with boundaries $(V_i ,
F_i)$ where the boundary $\partial V_i$ is a leaf of $F_i$. In fact
most of the constructions done in geometry to produce smooth
foliations of given codimension on a given manifold give a non
Hausdorff graph. The $C^*$-algebra $C^* (V,F)$ turns out in this
case to be obtained as a fibered product of the $C^* (V_i , F_i)$.

\begin{figure}
\begin{center}
\includegraphics[scale=0.6]{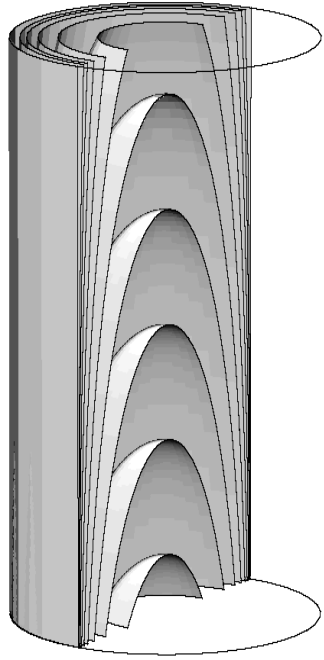}
\end{center}
\caption{The Reeb foliation \label{FigReeb}}
\end{figure}

In the general non-Hausdorff case the graph $G$ of $(V,F)$, being
non Hausdorff  may have only very few continuous fonctions with
compact support. However, being a manifold, we can give a local chart
$U \overset{\chi}{\to} {\mathbb R}^{\dim G}$. Take a smooth function
$\varphi \in C_c^{\infty} ({\mathbb R}^{\dim G})$, ${\rm Supp} \,
\varphi \subset \chi (U)$ and consider the function on $G$ equal to
$\varphi \circ \chi$ on $U$ and to $0$ outside $U$. If $G$ were
Hausdorff, this would generate all of $C_c^{\infty} (G)$ by taking
linear combinations, and in general we take this linear span as the
{\it definition} of $C_c^{\infty} (G)$. Note that we do not get {\it
continuous} functions, since there may well be a sequence $\gamma_n
\in U$ with two limits, one in ${\rm Supp} \, \varphi \circ \chi$
one in the complement of $U$. The above definition of $C_c^{\infty}
(G)$ obviously extends to get $C_c^{\infty} (G,\Omega^{1/2})$ the
space of smooth $\frac{1}{2}$ densities on $G$. One then shows that
the convolution $\varphi_1 * \varphi_2$ of $\varphi_1 , \varphi_2
\in C_c^{\infty} (G , \Omega^{1/2})$ is in $C_c^{\infty} (G ,
\Omega^{1/2})$.

Then we proceed exactly as in the Hausdorff case, and construct the
representation $\pi_x$ of the $*$ algebra $C_c^{\infty} (G ,
\Omega^{1/2})$ in the Hilbert space $L^2 (G_x)$. We note that though
$G$ is not Hausdorff, each $G_x$ is {\it Hausdorff},  being the
holonomy covering of the leaf through $x$.

For each $\varphi \in C_c^{\infty} (G , \Omega^{1/2})$ and $x \in
V$, $\pi_x (\varphi)$ is an ordinary smoothing operator, bounded in
$L^2 (G_x)$.

Exactly as in the Hausdorff case  $C^* (V,F)$ is defined as the
$C^*$ completion of $C_c^{\infty} (G , \Omega^{1/2})$ with norm
$\sup_{x \in V} \, \Vert \pi_x (\varphi) \Vert$.

\medskip

The obtained functor from foliations to $C^*$-algebras makes it
possible first
of all to translate from basic geometric properties to algebraic
ones and the simplest examples of foliations already exhibit
remarkable $C^*$-algebras. For instance the horocycle foliation of
the unit sphere bundle of a Riemann surface of genus $>1$ gives a
simple $C^*$-algebra without idempotent. The Kronecker foliation
gives rise to the noncommutative torus, which we   describe
in more detail in Section \ref{nctori}.

In the type II situation \ie in the presence of a holonomy invariant
transverse measure $\Lambda$ the basic result of the theory is the
longitudinal index theorem which computes the $L^2$-index of
differential operators $D$ on the foliated manifold $(V,F)$ which
are elliptic in the longitudinal direction (\ie $D$ restrict to
the leaves $L$ as elliptic operators $D_L$). One starts with a pair of
smooth vector bundles $E_1$, $E_2$ on $V$ together with a
differential operator $D$ on $V$ from sections of $E_1$ to sections
of $E_2$ such that:
\begin{enumerate}
\item[1)] $D$ restricts to leaves, \ie $(D \xi)_x$
only depends on the restriction of $\xi$ to a neighborhood of
$x$ in the leaf of $x$ ({\it i.e.} $D$ only uses partial
differentiation in the leaf direction).
\item[2)] $D$ is elliptic when restricted to any leaf.
\end{enumerate}

\begin{thm}\label{longmeasindex} \cite{Co78}

{\rm a)} There exists a Borel transversal $B$ (resp. $B'$) such that
the bundle $(\ell^2 (L \cap B))_{L \in V/F}$ is measurably
isomorphic to the bundle $({\rm Ker} \, D_L)_{L \in V/F}$ (resp. to
$({\rm Ker} \, D_L^*)_{L \in V/F})$.

\smallskip

\noindent {\rm b)} The scalar $\Lambda (B) < \infty$ is independent
of the choice of $B$ and noted $\dim_{\Lambda} ({\rm Ker} \, (D))$.

\smallskip

\noindent  ${\rm c)}\qquad\dim_{\Lambda} ({\rm Ker} \, (D)) -
\dim_{\Lambda} ({\rm Ker} \, (D^*)) = \,\varepsilon \, \langle {\rm
ch} \, \sigma_D \, {\rm Td} \, (F_{\mathbb C}) , [C] \rangle$

 ($\varepsilon = (-1) \, \frac{k(k+1)}{2}$, $k = \dim F$,
 ${\rm Td} \, (F_{\mathbb C})=$  Todd genus,
$\sigma_D=$ symbol of $D$) .

\end{thm}

Here $[C]\in H_k(V,\C)$ is the homology class of the Ruelle-Sullivan
current, a closed de-Rham current of dimension $k = \dim F$ which
encodes the transverse measure $\Lambda$ by integration of a
$k$-dimensional differential form $\omega$ on $V$ along the plaques
of foliation charts.

\bigskip

In particular the Betti numbers $\beta_j$ of a measured foliation
were defined in \cite{Co78} and give the $L^2$-dimension of the
space of $L^2$-harmonic forms along the leaves, more precisely one
has the following result.

{\rm a)} For each $j = 0,1,2,\ldots ,\dim F$, there exists a Borel
transversal $B_j$ such that the bundle $(H^j (L,{\mathbb C}))_{L \in
V/F}$ of $j$-th square integrable harmonic forms on $L$ is
measurably isomorphic to $(\ell^2 (L \cap B))_{L \in V/F}$.

\smallskip

\noindent {\rm b)} The scalar $\beta_j = \Lambda (B_j)$ is finite,
independent of the choice of $B_j$, of the choice of the Euclidean
structure on $F$.

\smallskip

\noindent {\rm c)} One has $\Sigma \, (-1)^j \, \beta_j = \chi
(F,\Lambda)$.

Here the Euler characteristic is simply given by the pairing of the
Ruelle-Sullivan current with the Euler class $e(F)$ of the oriented
bundle $F$ on $V$.

Extending ideas of Cheeger and Gromov \cite{ChGr} in the case of discrete groups,
D.~Gaboriau has shown in a remarkable recent work
(\cf \cite{Gabo1}, \cite{Gabo2}) that the Betti
numbers $\beta_j(F,\Lambda)$ of a foliation with contractible leaves
are invariants of the measured equivalence relation $\mathcal
R=\,\{(x,y)\,|\,y \in {\rm leaf}(x)\}$.

\medskip

In the general case one cannot expect to have a holonomy invariant
transverse measure and in fact the simplest foliations are of type
III from the measure theoretic point of view. Obtaining an analogue
in general of Theorem \ref{longmeasindex} was the basic motivation
for the construction of the assembly map (the second step of section
\ref{road}). Let us now briefly state the longitudinal index
theorem.

Let $D$ be as above an elliptic differential operator along the
leaves of the foliation $(V,F)$. Since $D$ is elliptic it has an
inverse modulo $C^* (V,F)$ hence it gives an element ${\rm Ind}_a
(D)$ of $K_0 (C^* (V,F))$ which is the {\em analytic index} of $D$.
The {\em topological index} is obtained as follows. Let $i$ be an
auxiliary imbedding of the manifold $V$ in ${\mathbb R}^{2n}$. Let
$N$ be the total space of the normal bundle to the leaves: $N_x =
(i_* (F_x))^{\perp} \subset {\mathbb R}^{2n}$. Let us foliate
$\tilde V = V \times \R^{2n}$ by $\tilde F$, $\tilde
F_{(x,t)} = F_x \times \{ 0 \}$, so that the leaves of $(\tilde V ,
\tilde F)$ are just $\tilde L = L \times \{ t \}$, where $L$ is a
leaf of $(V,F)$ and $t \in {\mathbb R}^{2n}$. The map $(x,\xi) \to
(x,i(x) + \xi)$ turns an open neighborhood of the $0$-section in $N$
into an open transversal $T$ of the foliation $(\tilde V , \tilde
F)$. For a suitable open neighborhood $\Omega$ of $T$ in $\tilde V$,
the $C^*$-algebra $C^* (\Omega , \tilde F)$ of the restriction of
$\tilde F$ to $\Omega$ is (Morita) equivalent to $C_0 (T)$, hence
the inclusion $C^* (\Omega , \tilde F) \subset C^* (\tilde V ,
\tilde F)$ yields a $K$-theory map: $K^0 (N) \to K_0 (C^* (\tilde V
, \tilde F))$. Since $C^* (\tilde V , \tilde F) = C^* (V,F) \otimes
C_0 (\R^{2n})$, one has, by Bott periodicity, the equality
$K_0 (C^* (\tilde V , \tilde F)) = K_0 (C^* (V,F))$.

Using the Thom isomorphism, $K^0 (F^*)$ is identified with $K^0 (N)$
so that one gets by the above construction the topological index:
$$
{\rm Ind}_t : K^0 (F^*) \to K_0 (C^* (V,F)) \, .
$$
The {\em longitudinal index theorem} \cite{CoSk}  is the equality
\begin{equation}\label{longind}
{\rm Ind}_a (D) = {\rm Ind}_t ([\sigma_D]),
\end{equation}
where $\sigma_D$ is the
longitudinal symbol of $D$ and $[\sigma_D]$ is its class in $K^0
(F^*)$.

Since the group $K_0 (C^* (V,F))$ is still fairly hard to compute
one needs computable invariants of its elements and this is where
cyclic cohomology enters the scene. In fact its early development
was already fully completed in 1981 for that precise goal  (\cf
\cite{khal}). The role of the trace on $C^* (V,F)$ associated to the
transverse measure $\Lambda$ is now played by cyclic cocycles on a
dense subalgebra of $C^* (V,F)$. The hard analytic problem is to
show that these cocycles have enough semi-continuity properties to
define invariants of $K_0 (C^* (V,F))$. This was achieved for some
of them in \cite{Co1} and makes it possible to formulate corollaries whose
statements are independent of the general theory, such as

\smallskip
\begin{thm}\label{transahat} \cite{Co1}
Let $M$ be a compact oriented manifold and assume that the $\hat
A$-genus
 $\hat A (M)$ is non-zero (since $M$ is not
assumed to be a spin manifold $\hat A (M)$ need not be an integer).
Let then $F$ be an integrable Spin sub-bundle of $TM$. There exists
no metric on $F$ for which the scalar curvature (of the leaves) is
strictly positive $(\geq \varepsilon > 0)$ on $M$.
\end{thm}

There is a very rich interplay between the theory of foliations and
their characteristic classes and operator algebras even at the
purely measure theoretic level \ie the classification of factors.

In a remarkable series of papers (see \cite{Hu} for references),
J.~Heitsch and S.~Hurder have analyzed the interplay between the
vanishing of the Godbillon-Vey invariant of a compact foliated
manifold $(V,F)$ and the type of the von Neumann algebra of the
foliation. Their work culminates in the following beautiful result
of S.~Hurder (\cite{Hu}). If the von Neumann algebra is {\it
semi-finite}, then the Godbillon-Vey invariant {\it vanishes}. We
have shown, in fact, that cyclic cohomology yields a stronger
result, proving that, if ${\rm GV} \ne 0$, then the central
decomposition of $M$ contains necessarily factors $M$, whose virtual
modular spectrum is of finite covolume in ${\mathbb R}_+^*$.

\begin{thm}\label{th3} \cite{Co1}
Let $(V,F)$ be an oriented, transversally oriented, compact,
foliated manifold, $({\rm codim} \, F = 1)$. Let $M$ be the
associated von Neumann algebra, and ${\rm Mod}(M)$ be its flow of
weights.
 Then, if the Godbillon-Vey class of
$(V,F)$ is different from $0$, there exists an invariant probability
measure for the flow ${\rm Mod}(M)$.
\end{thm}

 One actually constructs an invariant measure for the flow ${\rm
Mod}(M)$, exploiting the following remarkable property of the
natural cyclic 1-cocycle $\tau$ on the algebra $\cA$ of the
transverse 1-jet bundle for the foliation. When viewed as a linear
map $\delta$ from $\cA$ to its dual, $\delta$ is an unbounded
derivation, which is {\em closable}, and whose domain extends to the
center $Z$ of the von-Neumann algebra generated by $\cA$. Moreover,
$\delta$ vanishes on this center, whose elements $h\in Z$ can then
be used to obtain new cyclic cocycles $\tau_h$ on $\cA$. The pairing
$$
L(h)=\,\langle \tau_h,\, \mu(x)\rangle
$$
with the K-theory classes $\mu(x)$ obtained from the assembly map
$\mu$, which we had constructed with P. Baum \cite{BauCo},
then gives a measure on $Z$, whose invariance under the
flow of weights follows from the discreteness of the K-group. To
show that it is non-zero, one uses an index formula that evaluates
the cyclic cocycles, associated as above to the Gelfand-Fuchs
classes, on the range of  the assembly map $\mu$.

\medskip

The central question in the analysis of the noncommutative leaf
space of a foliation is step 3) (of section \ref{road}), namely the
metric aspect which entails in particular constructing a spectral
triple describing the transverse geometry. The reason why the
problem is really difficult is that it essentially amounts to doing
``metric" geometry on manifolds in a way which is ``background
independent" to use the terminology of physicists \ie which is {\em
invariant} under diffeomorphisms rather than covariant as in
traditional Riemannian geometry. Indeed the transverse space of a
foliation is a manifold endowed with the action of a large pseudo
group of partial diffeomorphisms implementing the holonomy. Thus in
particular no invariant metric exists in the general case and the
situation is very similar to trying to develop gravity without
making use of any particular ``background" metric that automatically
destroys the invariance under the action of diffeomeorphisms (\cf
\cite{CoMo-back}). Using the theory of hypoelliptic differential
operators and the basic technique of reduction from type III to type
II, a general construction of a spectral triple was done by
Connes-Moscovici in \cite{CoMo}. The remaining problem of the
computation of the local index formula in cyclic cohomology was
solved in \cite{CoMoHopf} and led in particular to the discovery of
new symmetries given by an action of a Hopf algebra which only
depends upon the transverse dimension of the foliation.

This also led to the development of the noncommutative analogue of
the Chern-Weil theory of characteristic classes \cite{CoMoHopf2} in
the general context of Hopf algebra actions on noncommutative spaces
and cyclic cohomology, a subject which  is undergoing rapid progress,
in particular thanks to the recent works of M. Khalkhali and
collaborators  \cite{kr1}, \cite{kr2}, \cite{kr3}, \cite{Hajac}.

\bigskip

\section{The noncommutative tori}\label{nctori}

This is perhaps considered as the prototype example of a
noncommutative space, since it illustrates very clearly
the properties and structures of noncommutative geometries.
Noncommutative tori played a key role in
the early developments of the theory in the 1980's (\cf
\cite{C3}), giving rise to noncommutative analogues of
vector bundles, connections, curvature, etc.

\smallskip

\begin{figure}
\begin{center}
\includegraphics[scale=0.8]{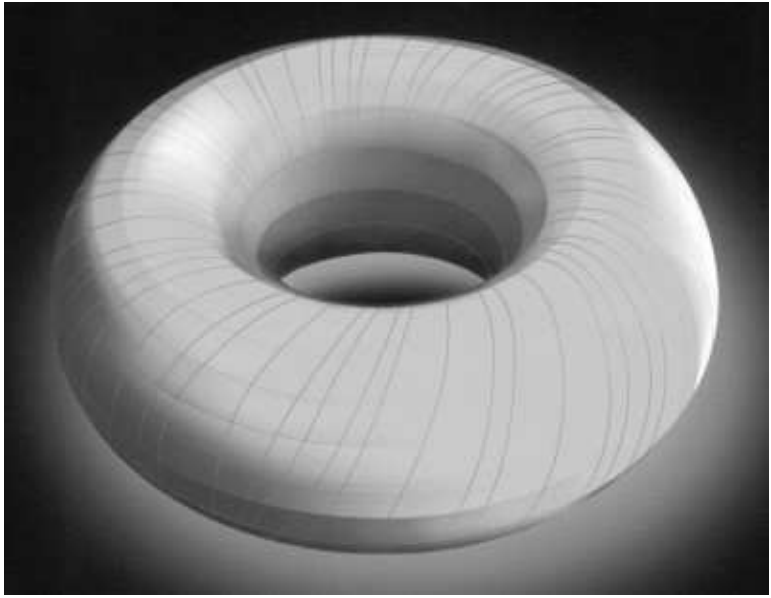}
\end{center}
\caption{The Kronecker foliation and the noncommutative
torus \label{FigNCtori}}
\end{figure}

One can regard noncommutative tori as a special case of
noncommutative spaces arising from foliations. In this case, one
considers certain vector fields on the ordinary 2-dimensional
real torus $T^2=\R^2/\Z^2$. In fact, one considers on $T^2$ the
Kronecker foliation $dx =\theta dy$, where $\theta$ is a given real
number. We are especially interested in the case where $\theta$ is
irrational. That is, we consider the
space of solutions of the differential
equation,
\begin{equation}
dx = \theta dy, \ \ \forall x,y \in \R/\Z \label{Krfol}
\end{equation}
for $\theta \in ]0,1[$ is a fixed irrational number. In other words, we
are considering the space of
leaves of the Kronecker foliation on the torus (\cf Figure \ref{FigNCtori}).

We can choose a transversal $T$ to the foliation, given by
$$ T=\{ y=0 \}, \ \ \ \  T\cong S^1\cong\R/\Z. $$
Two points of the transversal which differ by an integer multiple
of $\theta$ give rise to the same leaf.
We want to describe the further quotient
\begin{equation}
S^1 / \theta \Z \label{quottheta}
\end{equation}
by the equivalence relation which identifies any two points on the
orbits of the irrational rotation
\begin{equation}
R_{\theta} x = x + \theta  \mod 1 \, .
\label{Rtheta}
\end{equation}

\smallskip

We can regard the circle $S^1$ and the quotient space
\eqref{quottheta} at various levels of
regularity (smooth, topological, measurable). This corresponds to
different algebras of functions on the circle,
\begin{equation}
C^{\infty} (S^1) \subset C(S^1) \subset L^{\infty} (S^1) \, .
\label{S1}
\end{equation}

When passing to the quotient \eqref{quottheta}, if we just consider
invariant functions we obtain a very poor algebra of functions, since,
even at the measurable level, we would only have constant functions.
If instead we consider the noncommutative algebra of functions
obtained by the general recipe of ``noncommutative quotients''
(functions on the graph of the equivalence relation with the
convolution product), we obtain a very interesting and highly
non--trivial algebra of functions describing the space of leaves
of the foliation. This is given (in the topological category) by the
``irrational rotation algebra'', \ie the $C^*$-algebra
\begin{equation} \label{algNCT}
\cA_\theta: = \{ \left( a_{ij} \right)
\, \, \, \, i,j\in T \, \, \text{ in the same leaf } \}.
\end{equation}
Namely, elements in the algebra $\cA_\theta$ associated to the
transversal $T\simeq S^1$ are just matrices $(a_{ij})$ where the
indices are arbitrary pairs of elements $i,j$ of $S^1$ belonging
to the same leaf.

The algebraic rules are the same as for ordinary matrices. In the
above situation, since the equivalence is given by a group action,
the construction coincides with the crossed product. For
instance, in the topological category,
$\cA_\theta$ is identified with the crossed product $C^*$-algebra
\begin{equation}\label{irr-rot-alg}
\cA_\theta = C(S^1) \rtimes_{R_\theta} \Z.
\end{equation}

\smallskip

The algebra (\ref{algNCT}) has two natural generators:
\begin{equation} \label{U} U=\left\{ \begin{array}{ll}
1 & n=1 \\ 0 & \text{otherwise} \end{array}\right.  \end{equation}
and
\begin{equation} \label{VV} V=\left\{ \begin{array}{ll}
e^{2\pi ia} & n=0 \\ 0 & \text{otherwise}
\end{array} \right. \end{equation}
In fact, an element $b=(a_{ij})$ of $\cA_\theta$ can be written as
power series
\begin{equation}
b = \sum_{n \in \Z} b_n U^n ,
\label{powerU}
\end{equation}
where each $b_n$ is an element of the algebra (\ref{S1}), with
the multiplication rule given by
\begin{equation}
U h U^{-1} = h \circ R_{\theta}^{-1} \, . \label{RthetaU}
\end{equation}
The algebra \eqref{S1} is generated by the function $V$ on
$S^1$,
\begin{equation}
V(\alpha) = \exp (2\pi i \alpha) \ \ \ \forall \alpha \in S^1
\label{Valpha}
\end{equation}
and it follows that $\cA_\theta$ is generated by two unitaries
$(U,V)$ with presentation given by the relation
\begin{equation}
VU = \lambda \, U V, \ \ \text{ with } \ \  \lambda = \exp (2\pi i \theta) \, .
\label{[U,V]}
\end{equation}

\smallskip

If we work in the smooth category, then a generic element
$b$ in \eqref{S1} is given by a power series
\begin{equation}
b = \sum_{\Z^2} b_{nm} U^n V^m \ \in \cS (\Z^2)
\label{smoothNCt}
\end{equation}
where $\cS (\Z^2)$ is the Schwartz space of sequences of rapid
decay on $\Z^2$. We refer to the algebra of smooth functions \eqref{smoothNCt}
as $\cC^\infty (\bT^2_\theta)$, where we think of
$\bT^2_\theta$ as the (smooth) non--commutative torus.

\smallskip

Notice that in the definition (\ref{algNCT}) it is not necessary
to restrict to the condition that $i,j$ lie on the transversal
$T$. It is possible to also form an algebra
\begin{equation} \cB_\theta=\{ \left( a_{ij}
\right) \, \, \, \, i,j\in T^2 \, \, \text{ in the same leaf } \},
\end{equation}
where now the parameter of integration is no longer discrete. This
ought to correspond to the same non--commutative space, and in
fact the algebras are related by
$$ \cB_\theta = \cA_\theta \otimes \cK, $$
where $\cK$ is the algebra of all compact operators.

\medskip

The tangent space to the ordinary torus $T^2$ is spanned by the
tangent directions $\frac{\partial}{\partial x}$ and
$\frac{\partial}{\partial y}$. By choosing coordinates $U,V$, with
$U=e^{2\pi i x}$ and $V=e^{2\pi i y}$, the tangent vectors are
given by $\frac{\partial}{\partial x}=2\pi i
U\frac{\partial}{\partial U}$ and $\frac{\partial}{\partial
y}=2\pi i V\frac{\partial}{\partial V}$. These have analogs in
terms of derivations of the algebra of the non--commutative torus.
The two commuting vector fields which span the tangent space for
an ordinary (commutative) 2-torus correspond algebraically to two
commuting derivations of the algebra of smooth functions.

These derivations continue to make sense when we replace the
generators $U$ and $V$ of $C^{\infty} (\bT^2)$ by the generators of
the algebra $\cC^{\infty} (\bT_{\theta}^2)$, which no longer
commute, as shown in \eqref{[U,V]}. The derivations are still
given by the same formulas as in the commutative case,
\begin{equation}
\delta_1 = 2\pi i U \, \frac{\partial}{\partial U} \ \ \ \
\delta_2 = 2\pi i V \, \frac{\partial}{\partial V}
\label{derivationsUV}
\end{equation}
so that $\delta_1 \left( \sum b_{nm} U^n V^m \right) = 2\pi i \sum
n b_{nm} U^n V^m$, and similarly for $\delta_2$.

The operators (\ref{derivationsUV}) are commuting derivations of
the algebra $\cC^{\infty}(\bT_{\theta}^2)$.

In fact, it is straightforward to verify that $\delta_1$ and
$\delta_2$ satisfy
\begin{equation}
\delta_1 \delta_2 = \delta_2 \delta_1
\label{der-commuting}
\end{equation}
and
\begin{equation}
\delta_j (bb') = \delta_j (b) b' + b \delta_j (b') \quad \forall b
, b' \in \cA_\theta. \label{der-rule}
\end{equation}

Just as in the classical case of a (commutative)
manifold, what ensures that the derivations considered are enough
to span the whole tangent space is the condition of ellipticity
for the Laplacian
$$ \Delta = \delta_1^2 + \delta_2^2. $$
In Fourier modes the Laplacian is of the form $n^2 + m^2$, hence
$\Delta^{-1}$ is a compact operator.

\medskip

The geometry of the Kronecker foliation is closely related to the
structure of the algebra. In fact, a choice of a {\it closed
transversal} $T$ of the foliation corresponds
canonically to a {\it finite projective module} over the algebra
$\cA_\theta$.

In fact, the main result on finite projective module over the
non--commutative tori $\bT_{\theta}^2$ is the following
classification, which is obtained by combining the results of
\cite{PV}, \cite{C3}, \cite{Rieffel}.

\begin{thm}\label{NCTmodules}
Finite projective modules over $\cA_{\theta}$ are classified up to
isomorphism by a pair of integers $(p,q)$ such that $p+q \theta
\geq 0$. For a choice of such pair, the corresponding module
${\mathcal H}_{p,q}^{\theta}$ is obtained from the transversal
$T_{p,q}$ given by the closed geodesic of the torus $T^2$
specified by $(p,q)$, via the following construction. Elements of
the module associated to the transversal $T_{p,q}$ are rectangular
matrices, $(\xi_{i,j})$ with $(i,j) \in T \times S^1$, and with
$i$ and $j$ belonging to the same leaf. The right action of
$(a_{i,j}) \in \cA_{\theta}$ is by matrix multiplication.
\end{thm}

For instance, from the transversal $x=0$ one obtains the following
right module over $\cA_\theta$. The underlying linear space is the
usual Schwartz space
\begin{equation}
\cS (\R) = \{ \xi : \, \xi (s) \in \C, \quad \forall s \in \R
\} \label{SRmodule}
\end{equation}
of complex valued smooth functions on $\R$, all of whose
derivatives are of rapid decay.
The right module structure is given by the action of
the generators $U,V$
\begin{equation}
(\xi U) (s) = \xi (s+\theta) \ \ \ \ (\xi V) (s) = e^{2\pi is} \xi
(s) \ \ \ \forall s \in \R \, .
\label{actionUV}
\end{equation}
One of course checks that the relation \eqref{[U,V]} is satisfied,
and that, as a right module over
$\cA_\theta$, the space $\cS (\R)$ is {\it finitely generated}
and {\it projective} (\ie it complements to a free module).

\smallskip

Finitely generated projective modules play an important role in
noncommutative geometry, as they replace {\em vector bundles} in the
commutative setting. In fact, in ordinary commutative geometry,
one can equivalently describe vector bundles through their sections,
which in turn form a finite projective module over the algebra of
smooth functions. The
notion of finite projective module continues to make sense in the
non--commutative setting, and provides this way a good notion of
``non--commutative vector bundles''.

\smallskip

Suppose given a vector bundle $E$, described algebraically through its
space of smooth sections $C^{\infty} (X,E)$. One can
compute the dimension of $E$ by computing the trace of the identity
endomorphism. In terms of the space of smooth sections, hence of finite
projective modules, it is possible to recover the dimension of the
vector bundle as a limit
\begin{equation} \label{dimE}
\dim_{\cA}(\cE)= \lim_{N\to \infty}\frac{1}{N}\left( \#
\text{Generators of} \, \, \underbrace{\cE \oplus \cdots \oplus
\cE}_{N \text{times}} \right).
\end{equation}
This method applies to the noncommutative setting. In the case
of noncommutative tori, one finds that the Schwartz space $\cS(\R)$
has dimension the real number
\begin{equation}
\dim_{\cB} (\cS) = \theta \, .
\label{dimtheta}
\end{equation}
One similarly finds values $p+q\theta$ for the more general
case of Theorem \ref{NCTmodules}.

\smallskip

The appearance of a real values dimension is related to the {\em
density} of transversals in the leaves, that is, the limit of
$$ \frac{\# B_R \cap S}{\text{size of}\, \, B_R}, $$
for a ball $B_R$ of radius $R$ in the leaf. In this sense, the
dimension $\theta$ of the Schwartz space measure the relative
densities of the two transversals $S=\{ x=0 \}$ and $T=\{ y=0 \}$.

In general, the appearance of non integral dimension is a basic
feature of von
Neumann algebras of type II. The dimension of a vector bundle
is the only invariant that remains when one looks from the measure
theoretic point of view, using the algebra of measurable functions
$L^\infty(S^1)$ in (\ref{S1}). The von Neumann algebra which
describes the quotient space from the measure theoretic point of
view is the crossed product
\begin{equation}
R = L^{\infty} (S^1) \rtimes_{R_{\theta}} \Z.
\label{measNCT}
\end{equation}
This is the well known hyperfinite factor of type II$_1$. In
particular the classification of finite projective modules $\cE$
over $R$ is given by a positive real number, the Murray and von
Neumann dimension
\begin{equation}
\dim_R (\cE) \in \R_+ \, . \label{MvNR}
\end{equation}

\medskip

The simplest way to describe the phenomenon of Morita equivalence
for non--commutative tori is in terms of the Kronecker foliation,
where it corresponds to reparameterizing the leaves space in terms
of a different closed transversal. Thus, Morita equivalence of the
algebras $\cA_{\theta}$ and $\cA_{\theta'}$ for $\theta$ and
$\theta'$ in the same $\PGL(2,\Z)$ orbit becomes simply a
statement that the leaf--space of the original foliation is
independent of the transversal used to parameterize it. For
instance, Morita equivalence between $\cA_{\theta}$ and
$\cA_{-1/\theta}$ corresponds to changing the parameterization of
the space of leaves from the transversal $T=\{ y=0 \}$ to the
transversal $S=\{ x=0 \}$.

More generally, an explicit construction of bimodules
$\cM_{\theta,\theta'}$ was obtained in \cite{C3}. These are
given by the Schwartz space ${\mathcal S}(\R \times \Z/c)$, with
the right action of $\cA_\theta$ given by
$$ Uf\, (x,u)=f\left( x-\frac{c\theta+d}{c}, u-1\right) $$
$$ Vf\, (x,u)=\exp(2\pi i(x-u d/c)) f(x,u) $$
and the left action of $\cA_{\theta'}$
$$ U'f\, (x,u)=f\left( x-\frac{1}{c}, u-a \right) $$
$$ V'f\, (x,u)=\exp\left( 2\pi i \left(\frac{x}{c\theta+d}
-\frac{u}{c}\right)\right) f(x,u). $$
The bimodule $\cM_{\theta,\theta'}$ realizes the
Morita equivalences between $\cA_\theta$ and
$\cA_{\theta'}$ for
$$ \theta'=\frac{a\theta + b}{c\theta + d}=g\theta $$
with $g\in \PGL(2,\Z)$, \cf \cite{C3}, \cite{Ri2}.

\bigskip

\section{Duals of discrete groups}\label{discgr}

Noncommutative geometry provides naturally
a generalization of Pontrjagin duality for
discrete groups. While the Pontrjagin dual $\hat\Gamma$ of a
finitely generated discrete abelian group is a compact abelian group,
the dual of a more general finitely generated discrete group
is a noncommutative space.

To see this, recall that the usual Pontrjagin duality assigns to
a finitely generated discrete abelian group $\Gamma$ its group of
characters $\hat\Gamma=\Hom(\Gamma,U(1))$. The duality
is given by Fourier transform $e^{i\langle k, \gamma\rangle}$, for
$\gamma\in \Gamma$ and $k\in \hat\Gamma$.

In particular, Fourier transform gives an identification between
the algebra of functions on
$\hat\Gamma$ and the (reduced) $C^*$-algebra of
the group $\Gamma$,
\begin{equation}\label{CGammadual}
C(\hat\Gamma) \cong C^*_r(\Gamma),
\end{equation}
where the reduced $C^*$-algebra $C^*_r(\Gamma)$ is the
$C^*$-algebra generated by $\Gamma$ in the regular representation
on $\ell^2(\Gamma)$.

When $\Gamma$ is non-abelian Pontrjagin duality no
longer applies in the classical sense. However, the left hand side of
\eqref{CGammadual} still makes sense and it behaves ``like''
the algebra of functions on the dual group. One can then
say that, for a non-abelian group, the Pontrjagin dual
$\hat\Gamma$ still exists as a noncommutative space whose algebra
of coordinates is the $C^*$-algebra $C^*_r(\Gamma)$.

\medskip

As an example that illustrates this general philosophy we
give a different version of Example \ref{ex2}.

\begin{ex}{\em
The algebra (\ref{alg1}) of Example \ref{ex2}, is the group ring of
the dihedral group $\Z \rtimes \Z/2 \cong \Z/2 * \Z/2$. }
\label{dihedral}
\end{ex}

In fact, first notice that to a representation of the group $\Z/2 *
\Z/2$ (free product of two copies of the group with two elements) is
the same thing as a pair of subspaces in the Hilbert space,
$E,F\subset \cH$. The corresponding operators are $U=I-2P_E$,
$V=I-2P_V$, with $P_E,P_V$ the projections. The operators $U,V$
represent reflections, since $U=U^*$, $U^2=I$, $V=V^*$, $V^2=I$. The
group $\Gamma =\Z/2 * \Z/2$ realized as words in the generators $U$
and $V$ can equivalently be described as the semi-direct product
$\Gamma =\Z \rtimes \Z/2$, by setting $X=UV$, with the action
$UXU^{-1}=X^{-1}$. The regular representation of $\Gamma$ is
analyzed using Mackey's theory for semi-direct products. 0ne
considers first representations of the normal subgroup, and then
orbits of the action of $\Z/2$. The irreducible representations of
the normal subgroup $\Z$ are labeled by $S^1=\{z\in\C\,|z|=1\}$ and
given by $X^n\mapsto z^n$. The action of $\Z/2$ is the involution given
by conjugation $z \mapsto \bar z$. The quotient of $S^1$ by the
$\Z/2$ action is identified with the interval $[-1,1]$ by the map
$z\mapsto \Re( z)$. For points inside the interval the corresponding
irreducible representation of $\Gamma$ is two dimensional. At each
of the two endpoints $\pm 1$ one gets two inequivalent irreducible
representations of $\Gamma$. Thus we recover the picture of Example
\ref{ex2} and an isomorphism $C^*(\Gamma)\sim A$ where $A$ is the
algebra \eqref{alg1}.

\medskip
 The first two basic steps of the general theory are known
for arbitrary discrete groups $\Gamma$, namely

1) The
resolution of the diagonal and computation of the cyclic cohomology
are provided by the geometric model (due to Burghelea \cite{Burgh}) given by the
free loop space of the classifying space $B\Gamma$.

2) The assembly map (BC-map) of \cite{BauCo}
from the K-homology of the classifying
space $B\Gamma$ to the K-theory of the reduced $C^*$-algebra
$C^*_r(\Gamma)$ is refined in \cite{BCH} to take care of torsion
in the group
$\Gamma$ and gives a pretty good approximation to the K-theory of
$C^*_r(\Gamma)$ (\cf \cite{Skan1} and references therein).

In the presence of a natural
smooth subalgebra of $C^*_r(\Gamma)$ containing the group ring and
stable under holomorphic functional calculus, the combination of
the two steps described above makes it possible  to prove an index
theorem which is an higher dimensional form of Atiyah's $L^2$-index
theorem for coverings. This gave the first proof of the Novikov
conjecture for hyperbolic groups (\cite{CoMoind}). Since then the
analysis of dense smooth subalgebras has played a key role, in
particular in the ground breaking work of Vincent Lafforgue.
See \cite{BauCo}, \cite{Julg}, \cite{Laff}, \cite{Skan1}, \cite{Skan2}.

The next step, \ie the construction of a spectral geometry, is
directly related to geometric group theory. In general one cannot
expect to get a finite dimensional spectral triple since the growth
properties of the group give (except for groups of polynomial
growth) a basic obstruction (\cf \cite{Co-fr}). A general
construction of a theta summable spectral triple was given in
\cite{Co94} Section IV.9. Basically the transition from finitely
summable spectral triples to the theta summable ones is the
transition from finite dimensional geometry to the infinite
dimensional case. In the theta summable case the Chern character is
no longer a finite dimensional cyclic cocycle and one needs to
extend cyclic cohomology using cocycles with infinite support in the
$(b,B)$ bicomplex fulfilling a subtle growth condition. The general
theory of {\em entire} cyclic cohomology was developed in
\cite{Cotheta}. It is in general quite difficult to compute the
Chern character in the theta summable case and one had to wait quite
a long time until it was done for the basic example of discrete
subgroups of semi-simple Lie groups. This has been achieved in a
remarkable recent paper of M. Puschnigg \cite{Pusch} in the case of
real rank one.

The fourth step \ie the thermodynamics might seem irrelevant in the
type II context of discrete groups. However a small variant of the
construction of the group ring, namely the Hecke algebra associated
to an almost normal inclusion of discrete groups (in the sense
considered in \cite{BC})
suffices to meet the type III world. One of the open fields is to
extend the above steps 1), 2) and 3) in the general context of
almost normal inclusions of discrete groups, and to perform the
thermodynamical analysis in the spirit of \cite{CCM} in that
context.

\bigskip

\section{Brillouin zone and the quantum Hall effect}\label{qhe}

An important application to physics of the theory of
non--commutative tori was the development of a rigorous
mathematical model for the Integer Quantum Hall Effect (IQHE)
obtained by Bellissard and collaborators \cite{Bellissard1},
\cite{Bellissard2}, \cite{Co94}.

The classical Hall effect is a physical phenomenon first observed
in the XIX century \cite{Hall}. A very thin metal sample is
immersed in a constant uniform strong magnetic field orthogonal to
the surface of the sample. By forcing a constant current to flow
through the sample, the flow of charge carriers in the metal is
subject to a Lorentz force perpendicular to the current and the
magnetic field. The equation for the equilibrium of forces in the
sample is of the form
\begin{equation} \label{clHall} N e {\bf E} + {\bf j} \wedge {\bf B}
=0, \end{equation} where ${\bf E}$ is the electric field, $e$ and
$N$ the charge and number of the charge carriers in the metal,
${\bf B}$ the magnetic field, and ${\bf j}$ the current.

\begin{center}
\begin{figure}
\includegraphics[scale=1.1]{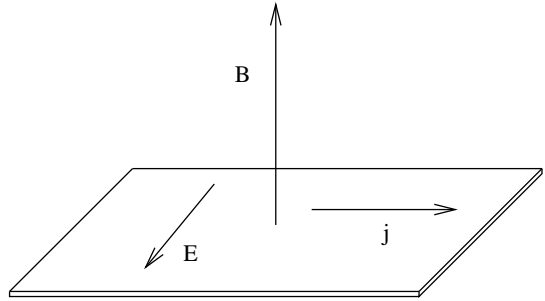}
\caption{The quantum Hall effect experiment.
\label{FigQHE}}
\end{figure}
\end{center}

The equation \eqref{clHall} defines a linear relation: the ratio
of the intensity of the Hall current to the intensity of the
electric field defines the Hall conductance,
\begin{equation} \label{condHall} \sigma_H = \frac{N e \delta}{B},
\end{equation}
with $B=|{\bf B}|$ the intensity of the magnetic field and
$\delta$ the sample width. The dimensionless quantity
\begin{equation} \label{filling} \nu_H = \frac{N\delta h}{B e} =
\sigma_H R_H \end{equation} is called the filling factor, while
the quantity $R_H = h/e^2$ is the Hall resistance. The filling
factor measures the fraction of Landau level filled by conducting
electrons in the sample. Thus, classically, the Hall conductance,
measured in units of $e^2/h$, equals the filling factor.

In 1980, about a century after the classical Hall effect was
observed, von Klitzing's experiment showed that, lowering the
temperature below 1 K, the relation of Hall conductance to filling
factor shows plateaux at integer values, \cite{vonKlitz}. The
integer values of the Hall conductance are observed with a
surprising experimental accuracy of the order of $10^{-8}$. This
phenomenon of quantization of the Hall conductance is known as
Integer Quantum Hall Effect (IQHE).

Laughlin first suggested that IQHE should be of a geometric origin
\cite{Lau}. A detailed mathematical model of the IQHE, which
accounts for all the important features of the experiment
(quantization, localization, insensitivity to the presence of
disorder, vanishing of direct conductance at plateaux levels)
improving over the earlier Laughlin model, was developed by
Bellissard and collaborators \cite{Bellissard1},
\cite{Bellissard2}.

Bellissard's approach to the IQHE is based on non--commutative geometry.
The quantization of the Hall conductance at integer values is
indeed geometric in nature: it resembles another well known
``quantization'' phenomenon that happens in the more familiar
setting of the geometry of compact 2--dimensional manifolds,
namely the Gauss--Bonnet theorem, where the integral of the
curvature is an integer multiple of $2\pi$, a property that is
stable under deformations. In the same spirit, the values of the
Hall conductance are related to the evaluation of a certain
characteristic class, or, in other words, to an index theorem for
a Fredholm operator.

More precisely, in the physical model one makes the simplifying
assumption that the IQHE can be described by non--interacting
particles. The Hamiltonian then describes the motion of a single
electron subject to the magnetic field and an additional potential
representing the lattice of ions in the conductor. In a perfect
crystal and in the absence of a magnetic field, there is a group
of translational symmetries. This corresponds to a group of
unitary operators $U(a)$, $a\in G$, where $G$ is the locally
compact group of symmetries. Turning on the magnetic field breaks
this symmetry, in the sense that translates of the Hamiltonian
$H_a = U(a) H U(a)^{-1}$ no longer commute with the Hamiltonian
$H$. Since there is no preferred choice of one translate over the
others, the algebra of observables must include all translates of
the Hamiltonian, or better their resolvents, namely the bounded
operators
\begin{equation} \label{translR} R_a(z) = U(a) (zI - H)^{-1}
U(a)^{-1}. \end{equation}

For a particle of (effective) mass $m$ and charge $e$ confined to
the plane, subject to a magnetic field of vector potential ${\bf
A}$ and to a bounded potential $V$, the Hamiltonian is of the form
\begin{equation} \label{magnLapl} H= \frac{1}{2m} \sum_{j=1,2} (p_j -e
A_j)^2 + V = H_0 + V, \end{equation} where the unperturbed part
$H_0$ is invariant under the magnetic translations, namely the
unitary representation of the translation group $\R^2$ given by
$$  U(a) \psi(x) = \exp\left( \frac{-ieB}{2\hbar} \omega(x,a) \right)
\psi(x-a), $$ with $\omega$ the standard symplectic form in the
plane. The hull (strong closure) of the translates (\ref{translR})
yields a topological space, whose homeomorphism type is
independent of the point $z$ in the resolvent of $H$. This
provides a non--commutative version of the Brillouin zone.

Recall that the Brillouin zones of a crystals are fundamental
domains for the reciprocal lattice $\Gamma^\sharp$ obtained via
the following inductive procedure. The {\em Bragg hyperplanes} of
a crystal are the hyperplanes along which a pattern of diffraction
of maximal intensity is observed when a beam of radiation (X-rays
for instance) is shone at the crystal. The $N$-th Brillouin zone
consists of all the points in (the dual) $\R^d$ such that the line
from that point to the origin crosses exactly $(n-1)$ Bragg
hyperplanes of the crystal.

\begin{figure}
\begin{center}
\includegraphics[scale=0.5]{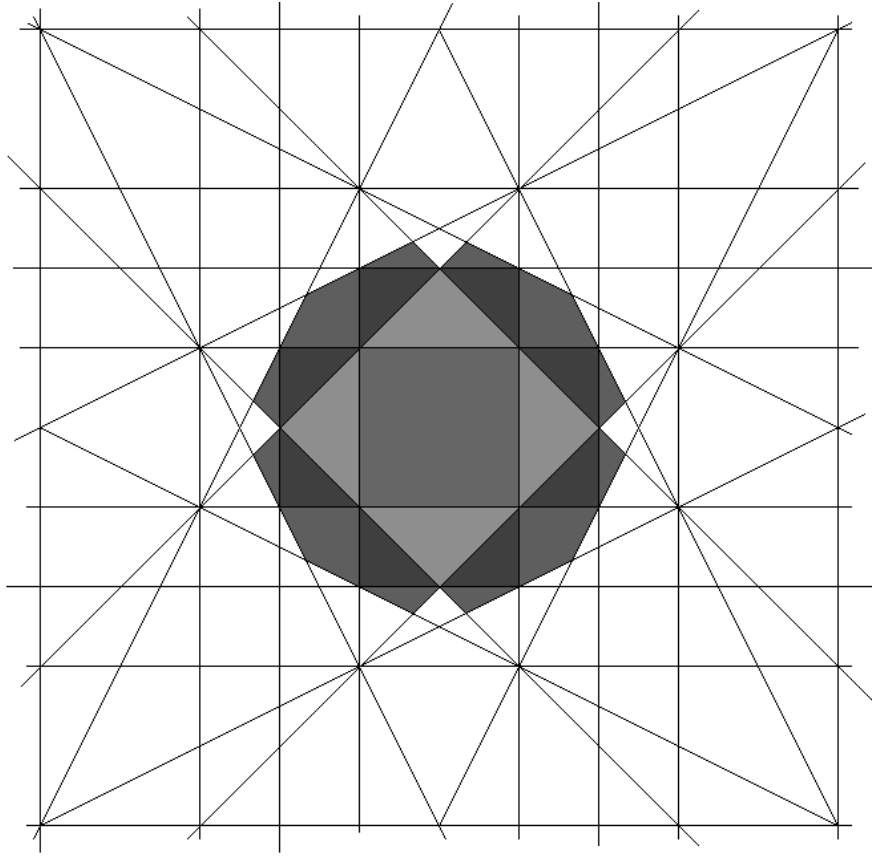}
\end{center}
\caption{Brillouin zones for a 2-dimensional crystal
\label{FigBrill}}
\end{figure}

More precisely, in our case, if $e_1$ and $e_2$ are generators of the
periodic lattice, we obtain a commutation relation
$$ U(e_1) U(e_2) = e^{2\pi i \theta} U(e_2) U(e_1), $$
where $\theta$ is the flux of the magnetic field through a
fundamental domain for the lattice, in dimensionless units, hence
the non-commutative Brillouin zone is described by a
non-commutative torus.

This can also be seen easily in a discrete model, where the
Hamiltonian is given by an operator
\begin{equation}\label{discr-Ham}
\begin{array}{rll}
(H_a \, f)(m,n) = & e^{-ia_1 n} f(m+1,n) + & e^{ia_2 n} f(m-1,n)
\\[3mm] + &
 e^{-ia_2 m} f(m, n+1) + & e^{ia_2 m } f(m, n-1),
\end{array}
\end{equation}
for $f\in L^2(\Z^2)$. This is a discrete version of the magnetic
Laplacian. Notice then that \eqref{discr-Ham} can be written in
the form
$$ H_a = U+V + U^* + V^*, $$
for
$$ (U\, f)(m,n)= e^{-ia_2m} f(m,n+1) \ \ \ \ (V\, f)(m,n)=
e^{-ia_1 n} f(m+1,n).
$$
These clearly satisfy the commutation relation \eqref{[U,V]} of
$\bT^2_\theta$ with $\theta= a_2-a_1$.

In the zero-temperature limit, the Hall conductance satisfies the
Kubo formula
\begin{equation}\label{kubo} \sigma_H = \frac{1}{2\pi i R_H}
\tau(P_\mu [\delta_1 P_\mu, \delta_2 P_\mu]),
\end{equation}
where $P_\mu$ is a spectral projection of the Hamiltonian on
energies smaller or equal to the Fermi level $E_\mu$, $\tau$ is
the trace on $\cA_\theta$ given by
\begin{equation} \label{trace-theta}
\tau\left( \sum a_{n,m} U^n V^m \right) = a_{0,0}.
\end{equation}
and $\delta_1$, $\delta_2$ are as in
\eqref{derivationsUV}. Here we assume that the Fermi level $\mu$
is in a gap in the spectrum of the Hamiltonian. Then the spectral
projections $P_\mu$ belong to the $C^*$-algebra of observables.

The Kubo formula \eqref{kubo} can be derived from purely physical
considerations, such as transport theory and the quantum adiabatic
limit.

The main result then is the fact that the integrality of the
conductance observed in the Integer Quantum Hall Effect is explained
topologically, that is, in terms of the integrality of the {\em
cyclic cocycle} $\tau(a^0 (\delta_1 a^1 \delta_2 a^2 - \delta_2 a^1
\delta_1 a^2))$ (\cf \cite{C3}).

\medskip

\begin{center}
\begin{figure}
\includegraphics[scale=0.6]{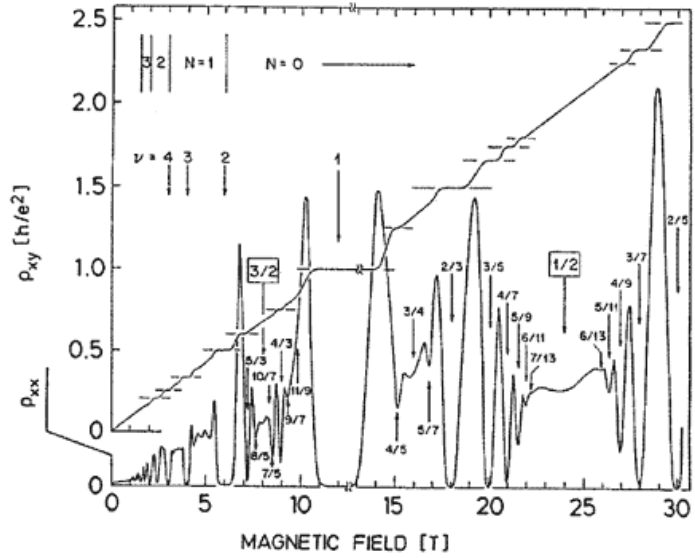}
\caption{Observed fractions in the quantum Hall effect.
\label{FigFQHE}}
\end{figure}
\end{center}

The fractional QHE was discovered by Stormer and Tsui in 1982. The
setup is as in the quantum Hall effect: in a high quality
semi-conductor interface, which will be modelled by an infinite
2-dimensional surface, with low carrier concentration and
extremely low temperatures $\sim 10 mK$, in the presence of a very
strong magnetic field, the experiment shows that the graph of
$\frac{h}{e^2}\sigma_H$ against the filling factor $\nu$ exhibits
plateaux at certain fractional values.

The independent electron approximation that, in the case of
the integer quantum Hall effect, reduces the problem to a single
electron wavefunction is no longer viable in this case and one has to
incorporate the Coulomb interaction between the electrons in a
many-electron theory. Nonetheless, it is possible to use a crude
approximation, whereby one alters the underlying geometry to
account for an average effect of the multi-electron interactions.
One can obtain this way a model of the fractional quantum Hall effect
via noncommutative geometry (\cf \cite{MM1}, \cite{MM2}, \cite{MM3}),
where one uses hyperbolic geometry to simulate the interactions.

The noncommutative geometry approach to the quantum Hall effect
described above was extended to hyperbolic geometry in \cite{CHMM}.
The analog of the operator \eqref{discr-Ham} is given by the Harper
operator on the Cayley graph of a finitely generated discrete subgroup
$\Gamma$ of $\PSL_2(\R)$. Given $\sigma: \Gamma \times
\Gamma \to U(1)$ satisfying $\sigma(\gamma,1)=\sigma(1,\gamma)=1$
and
$$ \sigma(\gamma_1,\gamma_2)\sigma(\gamma_1\gamma_2,\gamma_3)=
\sigma(\gamma_1,\gamma_2\gamma_3)\sigma(\gamma_2,\gamma_3), $$
one considers the right $\sigma$-regular representation on
$\ell^2(\Gamma)$ of the form
\begin{equation}\label{rightreg}
R^\sigma_\gamma
\psi(\gamma') = \psi(\gamma' \gamma) \sigma(\gamma', \gamma)
\end{equation}
satisfying $$R^\sigma_\gamma
R^\sigma_{\gamma'} = \sigma(\gamma,\gamma')
R^\sigma_{\gamma\gamma'}.$$
For $\{\gamma_i\}_{i=1}^r$ a symmetric set of generators of
$\Gamma$, the Harper operator is of the form
\begin{equation}\label{HarperGamma}
 \cR_\sigma = \sum_{i=1}^r  R^\sigma_{\gamma_i},
\end{equation}
and the operator $r-\cR_\sigma$ is the discrete analog of the magnetic
Laplacian (\cf \cite{Sun}).

The idea is that, by effect of the strong
interaction with the other electrons, a single electron ``sees''
the surrounding geometry as hyperbolic, with lattice sites that
appear (as a multiple image effect) as the points in a lattice
$\Gamma\subset \PSL_2(\R)$. Thus, one considers the general form
of such a lattice
\begin{equation}\label{FuchsianG}
\Gamma = \Gamma(g; \nu_1, \ldots, \nu_n),
\end{equation}
with generators $a_i, b_i, c_j$, with $i=1,\ldots,g$ and
$j=1,\ldots,n$ and a presentation of the form
\begin{equation}\label{presentation}
\Gamma(g; \nu_1, \ldots, \nu_n)=\langle a_i, b_i, c_j\,\,
\left|\,\, \prod_{i=1}^g [a_i,b_i]c_1\cdots c_n =1, \,\,\,
c_j^{\nu_j} =1 \right. \rangle.
\end{equation}
The quotient of the action of $\Gamma$ by isometrieson $\H$,
\begin{equation}\label{hyporb}
\Sigma(g; \nu_1, \ldots, \nu_n):= \Gamma \backslash \H,
\end{equation}
is a hyperbolic orbifold.

\begin{center}
\begin{figure}
\includegraphics[scale=0.6]{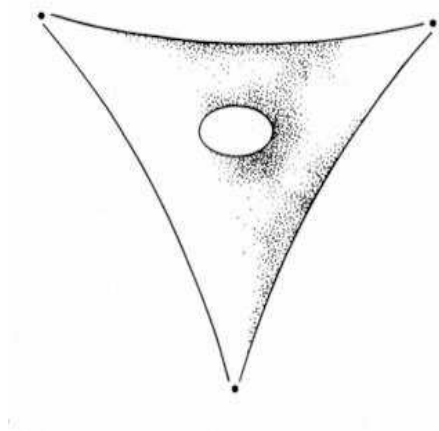}
\caption{Hyperbolic orbifolds.
\label{FigOrbi}}
\end{figure}
\end{center}

Let $P_E$ denote denote the spectral projection associated to the
Fermi level, \ie $P_E = \chi_{(-\infty, E]}(H)$. Then, in the zero
temperature limit, the Hall conductance is given by
\begin{equation}\label{sigmaE}
\sigma_E = \tr_K (P_E, P_E, P_E),
\end{equation}
where $\tr_K$ denotes the conductance 2-cocycle. It is a
cyclic 2-cocycle on the twisted group algebra $\C(\Gamma, \sigma)$
of the form
\begin{equation}\label{trK}
 \tr_K(f_0, f_1, f_2) =
\sum_{j=1}^g \tr(f_0  (\delta_j(f_1)\delta_{j+g}(f_2) -
  \delta_{j+g}(f_1)\delta_j(f_2))),
\end{equation}
where the $\delta_j$ are derivations associated to the
1-cocycles $a_j$ associated to a symplectic basis
$\{ a_j, b_j\}_{j=1,\ldots, g}$ of $H^1(\Gamma, \R)$ (\cf
\cite{MM3}).

Within this model, one obtains the fractional values of the Hall
conductance as integer multiples of orbifold Euler characteristics
\begin{equation}\label{orbEulch}
\chi_{orb}(\Sigma(g;\nu_1,\ldots,\nu_n))=2-2g +\nu -n  \in \Q.
\end{equation}
In fact, one shows (\cf \cite{MM2}, \cite{MM3}) that the conductance
2-cocycle is cohomologous to another cocycle, the area 2-cocycle, for
which one can compute the values on $K$-theory (hence the value of
\eqref{sigmaE}) by applying a twisted version of the Connes--Moscovici
higher index theorem \cite{CoMoind}.

\smallskip

While in the case of the integer quantum Hall effect the
noncommutative geometry model is completely satisfactory and explains
all the physical properties of the system, in the fractional case
the orbifold model can be considered
as a first rough approximation to the quantum field theory that
coverns the fractional quantum Hall effect. For instance, the geometry of
2-dimensional hyperbolic orbifolds is related to Chern--Simons theory
through the moduli spaces of vortex equations. This remains an interesting
open question.

\bigskip

\section{Tilings}\label{tiles}

In general, by a {\em tiling} $\cT$ in $\R^d$ one means the following. One
considers a finite collection $\{ \tau_1, \ldots, \tau_N \}$ of closed
bounded subsets of $\R^d$ homeomorphic to the unit ball. These are
called the {\em prototiles}. One usually assumes that the prototiles
are polytopes in $\R^d$ with a single $d$-dimensional cell which is
the interior of the prototile, but this assumption can be relaxed.
A tiling $\cT$ of $\R^d$ is then a covering of $\R^d$ by sets with disjoint
interior, each of which is a {\em tile}, that is, a translate of one
of the prototiles.

\begin{center}
\begin{figure}
\includegraphics[scale=0.8]{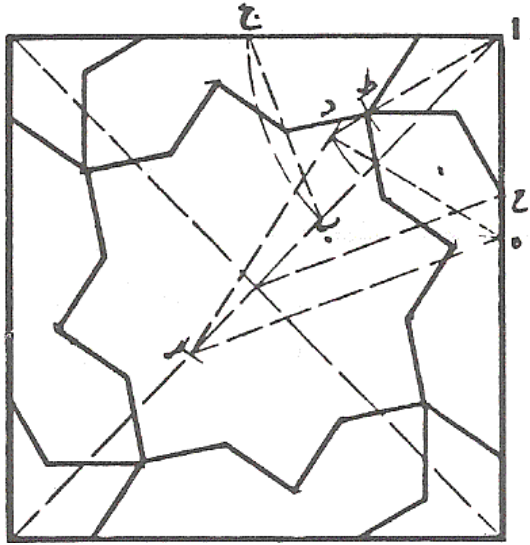}
\caption{Prototiles for tilings.
\label{FigTiles}}
\end{figure}
\end{center}

Given a tiling $\cT$ of $\R^d$ one can form its orbit closure
under translations. The metric on tilings is defined by saying that
two tilings are close if they almost agree on a large ball centered
at the origin in $\R^d$ (for more details and equivalent definitions
see \eg \cite{AndPut}, \cite{BellC}).

Tilings can be periodic or aperiodic. There are many familiar examples
of periodic tilings, while the best known examples of aperiodic
tilings are the Penrose tilings \cite{Penrose}. Similar types of
aperiodic tilings have been widely studied in the physics of
quasicrystals (\cf \eg \cite{BaMo}, \cite{BellC}).

It was understood very early on in the development of noncommutative
geometry (\cf \cite{Co-fr} and pp.5--7, pp.88--93, and pp.175--178 of
\cite{Co94}) that Penrose tilings provide an interesting class of
noncommutative spaces.

In fact, one can consider on the set $\Omega$ of tilings $\cT$ with
given prototiles $\{ \tau_1, \ldots, \tau_N \}$ the equivalence relation
given by the action of $\R^d$ by translations, \ie one identifies
tilings that can be obtained from one another by translations. In the
case of aperiodic tilings, this yields the type of quotient
construction described in Section \ref{quotients}, which leads
naturally to noncommutative spaces. An explicit description of this
noncommutative space for the case of Penrose tilings can be found in
\S II.3 of \cite{Co94}.

\begin{center}
\begin{figure}
\includegraphics[scale=0.6]{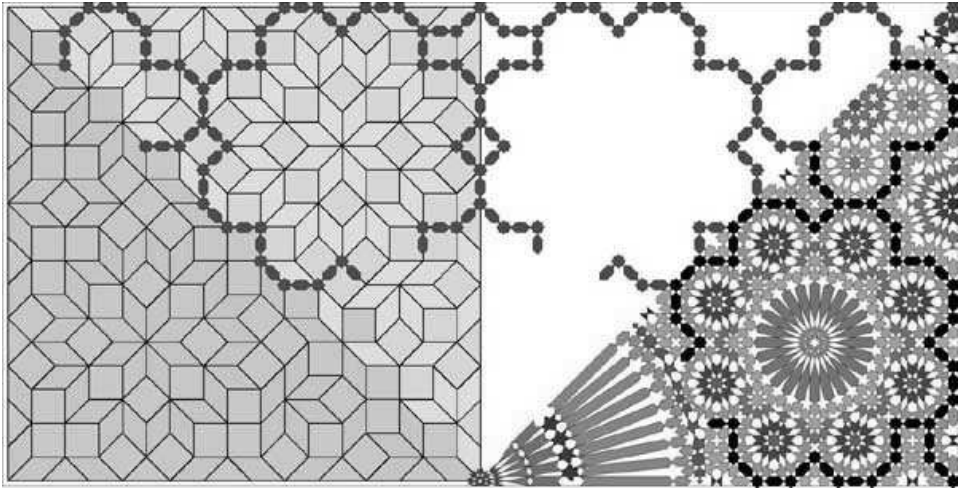}
\caption{Quasiperiodic tilings and zellijs.
\label{FigZell}}
\end{figure}
\end{center}

To simplify the picture slightly, we can consider the similar problem
(dually) with arrangements of points in $\R^d$ instead of
tilings. This is the formulation used in the theory of aperiodic
solids and quasicrystals (\cf \cite{BellC}). Then, instead of tilings
$\cT$, we consider discrete subsets of points $\cL\subset \R^d$. Such
$\cL$ is a Delauney set if there are radii $r,R>0$ such that every
open ball of radius $r$ meets $\cL$ in at most one point and every
closed ball of radius $R$ meets $\cL$ in at least one point. One can
describe $\cL$ by the counting measure
$$ \mu_{\cL}(f) =\sum_{x\in\cL} f(x), $$
and one can take the orbit closure $\Omega$ of the action of $\R^d$ by
translations $$\mu_{\cL}\mapsto T_{-a}\mu_{\cL}=\mu_{\cL}\circ T_a, \
\ \ \text{ for } a\in \R^d, $$
in the space $\cM(\R^d)$ of Radon measures with the weak$^*$
topology. The Hull of $\cL$ is the dynamical system $(\Omega,T)$,
where $T$ denotes the action of $\R^d$ by translations.

This dynamical system determines a corresponding noncommutative space,
describing the quotient of $\Omega$ by translations, namely the
crossed product $C^*$-algebra
\begin{equation}\label{crossOmega}
\cA=C(\Omega)\rtimes_T \R^d.
\end{equation}
In fact, one can also consider the groupoid with set of units the
transversal
\begin{equation}\label{transvOmega}
X=\{ \omega\in \Omega: 0\in {\rm Support}(\omega)\},
\end{equation}
arrows of the form $(\omega,a)\in \Omega\times \R^d$, with
source and range maps $s(\omega,a)=T_{-a}\omega$, $r(\omega,a)=\omega$
and $(\omega,a)\circ (T_{-a}\omega,b)=(\omega, a+b)$ (\cf
\cite{BellC}). This defines a locally compact groupoid $\cG(\cL,X)$.
The $C^*$-algebras $C^*(\cG(\cL,X))$ and $C(\Omega)\rtimes_T \R^d$ are
Morita equivalent.

In the case where $\cL$ is a periodic arrangement of points with
cocompact symmetry group $\Gamma\subset \R^d$, the space
$\Omega$ is an ordinary commutative space, which is topologically
a torus $\Omega=\R^d/\Gamma$. The $C^*$-algebra $\cA$ is in this case
isomorphic to $C(\hat\Gamma)\otimes \cK$, where $\cK$ is the algebra
of compact operators and $\hat\Gamma$ is the Pontrjagin dual of the
abelian group $\Gamma\cong \Z^d$, isomorphic to $T^d$, obtained by taking
the dual of $\R^d$ modulo the reciprocal lattice
\begin{equation}\label{recLat}
 \Gamma^\sharp = \{ k \in \R^d : \langle k, \gamma\rangle\in 2
\pi \Z, \forall \gamma \in \Gamma \}.
\end{equation}
Thus, in physical language, $\hat\Gamma$ is identified with the
Brillouin zone $B=\R^d/\gamma^\sharp$ of the periodic
crystal $\cL$ (\cf Section \ref{qhe}). In this periodic case, the
transversal $X=\cL/\Gamma$ is a finite set of points. The groupoid
$C^*$-algebra $C^*(\cG(\cL,X))$ is in this case isomorphic to
$C(\hat\Gamma)\otimes M_k(\C)$, where $k$ is the cardinality of the
transversal $X$. Thus, the periodic case falls back into the realm of
commutative spaces, while the aperiodic patterns give rise to truly
noncommutative spaces, which are highly nontrivial and interesting.

\begin{center}
\begin{figure}
\includegraphics[scale=0.6]{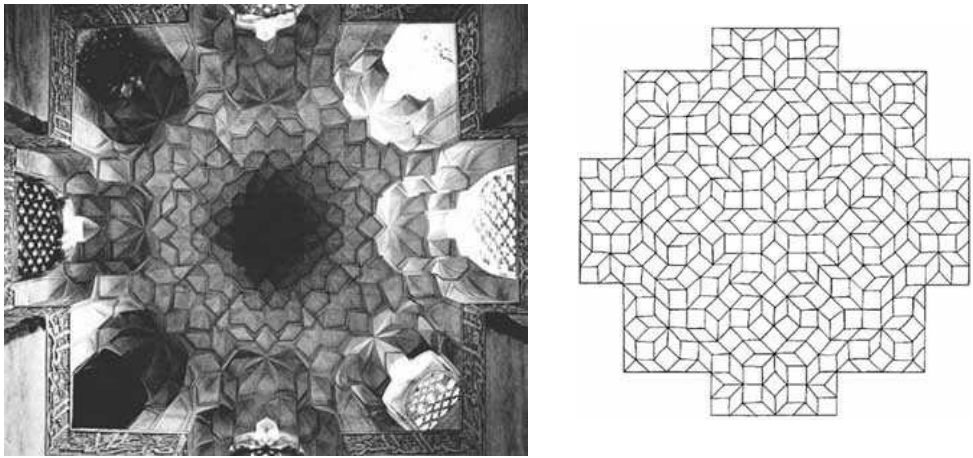}
\caption{Quasiperiodic tilings and muqarnas.
\label{FigMuq}}
\end{figure}
\end{center}

One of the richest sources of interesting tilings are the
{\em zellijs} and {\em muqarnas} widely used in ancient architecture.
Also collectively defined as ``arabesques'', not only these patterns
exhibit highly nontrivial geometries, but they reflect
the intricate interplay between philosophy, mathematics,
and aesthetics (\cf \cite{ArdBak}, \cite{Bulatov}).
Some of the best studies on
zellijs and muqarnas concentrate on 2-dimensional {\em periodic}
patterns. For instance we find in \cite{ArdBak}, p.43:
\begin{quote}
``As Nature is based on rhythm, so the arabesque is rhythmic
in concept. It reflects movement marked by the regular recurrence
of features, elements, phenomena; hence it has periodicity.''
\end{quote}
It seems from this viewpoint that only the theory of periodic
tilings (\ie commutative geometry) should be relevant in this
context. However, more recent
studies (\cf \cite{Bulatov}, \cite{Castera}, \cite{Castera2},
\cite{Nagy}) suggest that the design of zellijs and muqarnas was not
limited to 2-dimensional crystallographic groups, but, especially
during the Timurid period, it
involved also {\em aperiodic} patterns with fivefold symmetry, analogous to
those observed in quasi-crystals. This is no accident and was
certainly the result of a highly developed geometric theory:
already in the historic textbook of Abu'l-Wafa' al-Buzjani (940-998)
on geometric constructions \cite{Wafa} there is explicit mention of
meetings and discussions where mathematicians were directly involved
alongside artisans in the design of arabesque patterns.

The appearance of aperiodic
tilings is documented in the anonymous Persian manuscript
\cite{Anonym} ``On interlocking similar and congruent figures'', which
dates back to the 11th-13th century. Some of these aperiodic aspects of
zellijs and muqarnas were studied by Bulatov in the
book \cite{Bulatov}, which also contains Vil'danova's Russian
translation of the ancient Persian text.

\bigskip

\section{Noncommutative spaces from dynamical systems}\label{dynsys}

We will look at some examples of noncommutative spaces associated to
a discrete dynamical system $T$, for instance given by a self
mapping of a Cantor set. Such noncommutative spaces
 have been extensively studied in a
series of papers (\cf \cite{GPS} and \cite{skau} for a survey) where
C. Skau and his coworkers have obtained  remarkable results on the
classification of minimal actions of $\Z$ on Cantor
  sets using the  $K$-theory of the associated $C^*$-algebra.

It was found recently (\cf \cite{CM1}, \cite{CM2}, \S 4 of
\cite{Mar} and \S 8 of \cite{MaPa})  that the mapping torus of such
systems can be used to model the ``dual graph'' of the fibers at the
archimedean primes of arithmetic surfaces, in Arakelov geometry, in
the particular case in which the dynamical system $T$ is a subshift
of finite type encoding the action of a Schottky group
$\Gamma\subset \SL_2(\C)$ on its limit set $\Lambda_\Gamma\subset
\P^1(\C)$. In fact, the results of \cite{CM1} were motivated by
earlier results of Manin \cite{Man-hyp} that provided a geometric
model for such dual graphs in terms of hyperbolic geometry and
Schottky uniformizations.

\smallskip

More generally, given an alphabet with letters $\{ \ell_1, \ldots, \ell_N \}$,
the space $\cS^+_A$ of a subshift of finite type consists
of all right-infinite {\em admissible} sequences
\begin{equation}\label{admseqplus}
w=a_0 a_1 a_2 \ldots a_n \ldots
\end{equation}
in the letters of the alphabet. Namely, $a_i\in \{ \ell_1, \ldots,
\ell_N \}$ subject to an admissibility condition specified by an
$N\times N$ matrix $A$ with entries in $\{0,1\}$. Two letters
$\ell_i$ and $\ell_j$ in the alphabet can appear as consecutive digits
$a_k$, $a_{k+1}$ in the word $w$ if and only if the entry $A_{ij}$ of
the admissibility matrix $A$ is equal to $1$.
One defines similarly the space $\cS_A$ as the set of doubly-infinite
admissible sequences
\begin{equation}\label{admseq}
w=\ldots a_{-m}\ldots a_{-2} a_{-1} a_0 a_1 a_2 \ldots a_n \ldots
\end{equation}
The sets $\cS^+_A$ and $\cS_A$ have a natural choice of a topology.
In fact, on $\cS_A$ we can put the topology generated by the
sets $W^s(x,\ell) =\{ y\in {\mathcal S}_A | x_k = y_k, k\geq
\ell \}$, and the $W^u(x,\ell)=\{ y\in {\mathcal S}_A | x_k = y_k, k\leq
\ell \}$ for $x\in {\mathcal S}_A$ and $\ell \in \Z$. This induces a
topology with analogous properties on ${\mathcal S}^+_A$ by realizing it as a
subset of ${\mathcal S}_A$, for instance, by extending each sequence to
the left as a constant sequence. One then considers on $\cS_A$ (or on
$\cS_A^+$) the action of the two-sided (resp.~one-sided) shift
$T$ defined by $(Tw)_k = a_{k+1}$, where the $a_k$ are the digits of
the word $w$. Namely, the one-sided shift on $\cS_A^+$ is of the form
\begin{equation}\label{shift+}
T( a_0 a_1 a_2 \ldots a_\ell \ldots) = a_1 a_2 \ldots
a_\ell \ldots
\end{equation}
while the two-sided shift on $\cS_A$ acts as
\begin{equation}\label{shift}
\begin{array}{rcccccccccl} T(& \ldots & a_{-m} & \ldots & a_{-1}& a_0
& a_1 &
\ldots & a_{\ell} & \ldots &) = \\ & \ldots & a_{-m+1} & \ldots & a_{0} &
a_1 & a_2 & \ldots & a_{\ell+1} & \ldots & \end{array}
\end{equation}
Tipically spaces ${\mathcal S}^+_A$ and ${\mathcal S}_A$ are topologically
Cantor sets. The one-sided shift $T$ of \eqref{shift+} is a
continuous surjective map on ${\mathcal S}^+_A$, while the two-sided
shift $T$ of \eqref{shift} is a homeomorphism of ${\mathcal S}_A$.

For example, let $\Gamma$ be a free group in $g$ generators $\{
\gamma_1,\ldots, \gamma_g \}$. Consider the alphabet
$\{ \gamma_1,\ldots, \gamma_g, \gamma_1^{-1}, \ldots, \gamma_g^{-1}
\}$. Then one can consider the right-infinite, or doubly-infinite
words in these letters, without cancellations, that is, subject to the
admissibility rule that $a_{k+1}\neq a_k^{-1}$. This defines a
subshift of finite type where the matrix $A$ is the symmetric
$2g\times 2g$ matrix with $A_{ij}=0$ for $|i-j|=g$ and $A_{ij}=1$
otherwise.  Suppose that $\Gamma$ is
a Schottky group of genus $g$, \ie a finitely generated discrete
subgroup $\Gamma\subset \SL_2(\C)$, isomorphic to a free group in $g$
generators, where all nontrivial elements are hyperbolic.
Then the points in $\cS_A^+$ parameterize points in the
limit set $\Lambda_\Gamma\subset \P^1(\C)$ (the set of accumulation
points of orbits of $\Gamma$).
The points in $\cS_A$ parameterize geodesics in the
three dimensional real hyperbolic space $\H^3$ with ends at points on
the limit set $\Lambda_\Gamma$.

The pair $({\mathcal S}_A,T)$ is a typical example of an interesting
class of dynamical systems, namely it is a {\em Smale
space}. This means that locally ${\mathcal S}_A$ can be decomposed
as the product of expanding and contracting directions for $T$.
Namely, the following
properties are satisfied.
\begin{itemize}
\item For every point $x\in {\mathcal S}_A$ there exist subsets $W^s(x)$
and $W^u(x)$ of ${\mathcal S}_A$, such that $W^s(x) \times  W^u(x)$ is
homeomorphic to a neighborhood of $x$.
\item The map $T$ is contracting on $W^s(x)$ and expanding on
$W^u(x)$, and $W^s(Tx)$ and $T(W^s(x))$ agree in some neighborhood of
$x$, and so do $W^u(Tx)$ and $T(W^u(x))$.
\end{itemize}

A construction of Ruelle shows that one can
associate different $C^*$--algebras to Smale spaces (\cf \cite{Rue2},
\cite{Putn}, \cite{PutSpi}). For Smale spaces like $({\mathcal
S}_A,T)$ there are four basic possibilities: the crossed product algebra
$C({\mathcal S}_A)\rtimes_T \Z$ and the ${\rm C}^*$--algebras ${\rm
C}^*({\mathcal G}^s)\rtimes_T \Z$, ${\rm C}^*({\mathcal G}^u)\rtimes_T
\Z$, ${\rm C}^*({\mathcal G}^a)\rtimes_T \Z$ obtained by considering
the action of the shift $T$ on the groupoid ${\rm C}^*$--algebra
associated to the groupoids ${\mathcal G}^s$, ${\mathcal G}^u$,
${\mathcal G}^a$ of the  stable, unstable, and asymptotic equivalence
relations on $({\mathcal S}_A,T)$.

The first choice, $C({\mathcal S}_A)\rtimes_T \Z$, is
closely related to the continuous dynamical system given by the
mapping torus of $T$, while a choice like ${\rm C}^*({\mathcal
G}^u)\rtimes_T \Z$ is related to the ``bad quotient'' of $\cS_A^+$ by
the action of $T$. In the example of the Schottky group this
corresponds to the action of $\Gamma$ on its limit set.

One can consider the suspension flow $\cS_T$ of a dynamical
system $T$, that is, the mapping torus of the dynamical
system $( {\mathcal S}_A, T)$, which is defined as
\begin{equation} \label{suspensionT}
{\mathcal S}_T := {\mathcal S}_A
\times [0,1] / (x,0)\sim (Tx,1).
\end{equation}

\begin{figure}
\begin{center}
\includegraphics[scale=0.65]{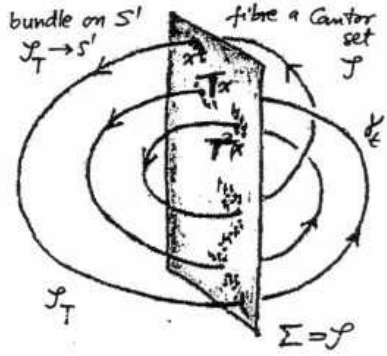}
\end{center}
\caption{Mapping Torus \label{FigMapTorus}}
\end{figure}

The first cohomology group of $\cS_T$ is
the ``ordered cohomology'' of the dynamical system $T$, in the sense of
\cite{BoHa} \cite{PaTu}.
There is an identification of $H^1({\mathcal S}_T,\Z)$ with the
$K_0$-group of the crossed product ${\rm C}^*$-algebra for the
action of $T$ on ${\mathcal S}_A$,
\begin{equation}
H^1({\mathcal S}_T,\Z) \cong K_0({\rm C}({\mathcal
S}_A) \rtimes_T \Z).  \label{H1K0}
\end{equation}
This can be seen from the Pimsner--Voiculescu exact sequence (\cf
\cite{PV}) for the $K$-theory of a crossed product by $\Z$, which in
this case reduces to
\begin{equation}\label{PV} 0 \to K_1({\rm C}({\mathcal S})\rtimes_T
\Z) \to {\rm C}({\mathcal S},\Z) \stackrel{I-T_*}{\to} {\rm
C}({\mathcal S},\Z) \to K_0({\rm C}({\mathcal S})\rtimes_T \Z) \to
0, \end{equation}
It can also be seen in terms of the Thom isomorphism
of \cite{Co-Thom}, \cite{Co1}.

In fact, as we discussed in Section \ref{road}, one of the
fundamental construction of noncommutative geometry (\cf
\cite{Co1}) is that of {\em homotopy
quotients}. These are commutative spaces, which provide, up to
homotopy, geometric models for the corresponding noncommutative
spaces. The noncommutative spaces themselves, as we are going to show
in our case, appear as quotient spaces of foliations on the homotopy
quotients with contractible leaves.

For the noncommutative space ${\mathcal S}_A/\Z$, with $\Z$ acting as
powers of the invertible two-sided shift, the homotopy quotient
is precisely the mapping torus \eqref{suspensionT},
\begin{equation}\label{htpyquot2}
{\mathcal S}_T = {\mathcal S}\times_\Z \R.
\end{equation}
The noncommutative space ${\mathcal S}/\Z$ can be
identified with the quotient space of the natural foliation on
\eqref{htpyquot2} whose generic leaf is contractible (a copy of
$\R$).

Another noncommutative space associated to a subshift
of finite type $T$ (which, up to Morita equivalence, correspond to
another choice of the $C^*$-algebra of a Smale space, as mentioned
above) is the Cuntz--Krieger algebra $\cO_A$,
where $A$ is the admissibility matrix of the subshift finite type
(\cf \cite{Cu} \cite{CuKrie}).

A partial isometry is a linear operator $S$ satisfying the relation
$S= S S^* S$. The Cuntz--Krieger algebra ${\mathcal O}_A$ is defined
as the universal $C^*$--algebra generated by partial isometries
$S_1, \ldots, S_{N}$, satisfying the relations
\begin{equation} \label{CK1rel} \sum_j S_j S_j^* =I \end{equation}
\begin{equation} \label{CK2rel} S_i^* S_i =\sum_j A_{ij} \, S_j
S_j^*. \end{equation}

In the case of a Schottky group $\Gamma\subset \PSL_2(\C)$
of genus $g$, the Cuntz--Krieger
algebra ${\mathcal O}_A$ can be described in terms
of the action of the free group $\Gamma$ on its limit set
$\Lambda_\Gamma\subset \P^1(\C)$ (\cf \cite{Rob}, \cite{Spi}), so that we can regard
${\mathcal O}_A$ as a noncommutative space replacing the classical
quotient $\Lambda_\Gamma / \Gamma$,
\begin{equation} \label{treefree}
{\mathcal O}_A \cong
C(\Lambda_\Gamma) \rtimes \Gamma.
\end{equation}
The quotient space
\begin{equation}\label{htpyquot1}
\Lambda_\Gamma \times_\Gamma \H^3 = \Lambda_\Gamma \times_\Gamma
\underline{E}\Gamma,
\end{equation}
is precisely the homotopy quotient of $\Lambda_\Gamma$ with respect to
the action of $\Gamma$, with $\underline{E}\Gamma =\H^3$ and the
classifying space $\underline{B}\Gamma = \H^3/\Gamma$. Here
$\H^3/\Gamma$ is a hyperbolic 3-manifold of infinite volume, which
is topologically a handlebody of genus $g$.
In this case also we find
that the noncommutative space $\Lambda_\Gamma /\Gamma$ is the quotient
space of a foliation on the homotopy quotient \eqref{htpyquot1} with
contractible leaves $\H^3$.

\bigskip

\section{Noncommutative spaces from string theory}\label{string}

The main aspects of string and D-brane theory that
involve noncommutative geometry are the bound states of
configurations of parallel D-branes \cite{Wit1}, the matrix models
for M-theory \cite{BFSS} and the strong coupling limit of
string theory (\cf \eg \cite{Ard1}, \cite{Ard2}). It also plays an
important role in the M-theory compactifications \cite{CDS}.
We shall not discuss all these aspects in detail here. Since the
focus of this review is on examples we   only mention a
couple of examples of noncommutative spaces arising from
string and D-brane theory.

\medskip

The noncommutative tori and the
components of the Yang-Mills connections appear in the
classification of the BPS states in M-theory \cite{CDS}.

Recall first that Yang--Mills theory on noncommutative tori
can be formulated (\cf \cite{CoRi}) using suitable notions of
connections and curvature for noncommutative spaces. In fact,
the analogs of connection
and curvature of vector bundles are straightforward to obtain
(\cite{C3}): a connection is just given by the associated
covariant differentiation $\nabla$ on the space of smooth
sections. Thus here it is given by a pair of linear operators,
\begin{equation}
\nabla_j : \cS (\R) \rightarrow \cS (\R)
\label{nabla}
\end{equation}
such that
\begin{equation}
\nabla_j (\xi b) = (\nabla_j \xi)b + \xi \delta_j (b) \quad
\forall \xi \in \cS \ , b \in \cA_\theta \, .
\label{connection}
\end{equation}
One checks that, as in the usual case, the trace of the curvature
$$\Omega = \nabla_1 \nabla_2 - \nabla_2 \nabla_1, $$
is independent of the choice of the connection.

We can make the following choice for the connection:
\begin{equation}
(\nabla_1 \xi) (s) = - \frac{2\pi i s}{\theta} \, \xi (s) \
\,\,\,\,\,\, (\nabla_2 \xi)(s) = \xi' (s) \, . \label{connection2}
\end{equation}
Notice that, up to the correct powers of
$2\pi i$, the total curvature of $\cS$ is an {\em integer}.
In fact, the curvature $\Omega$ is constant, equal to $\frac{1}{\theta}$,
so that the irrational number $\theta$ disappears in the total
curvature, $\theta \times \frac{1}{\theta}$.
This integrality phenomenon, where the pairing of dimension and
curvature (both of which are non--integral) yields an integer:
$$ \dim \times \Omega \sim \theta \times \frac{1}{\theta} =
\text{integer}, $$ is the basis for the development of a theory of
characteristic classes for non--commutative spaces. In the general
case, this requires the development of more sophisticated tools,
since analogs of the derivations $\delta_i$ used in the
case of the noncommutative tori are not there in general.
The general theory is obtained through cyclic
homology, as developed in \cite{ConnesCH}.

\smallskip

Consider then the modules ${\mathcal H}^\theta_{p,q}$ described in
Section \ref{nctori}. It is possible to define an $\cA_\theta$ valued
inner product $\langle \cdot, \cdot \rangle_{\cA}$ on ${\mathcal
H}^\theta_{p,q}$, as in \cite{Rieffel}, which is used to show that
${\mathcal H}^\theta_{p,q}$ is a projective module.
Connections are required to be compatible with the metric,
\begin{equation}\label{metr-conn}
 \delta_j \langle\xi,\eta \rangle_{\cA} =\langle \nabla_j \xi,\eta
\rangle_{\cA} +
\langle \xi, \nabla_j \eta \rangle_{\cA}.
\end{equation}
It is proved in \cite{C3} that such connections always exist. The
curvature $\Omega$ has values in $\cE=\End_{\cA}({\mathcal H})$.
An $\cE$--valued inner product on ${\mathcal H}$ is given by
$$ \langle \xi,\eta \rangle_{\cE} \zeta = \xi \langle \eta, \zeta
\rangle_{\cA}, $$ and a canonical faithful trace $\tau_E$ is
defined as
$$ \tau_E(\langle \xi,\eta \rangle_{\cE}) = \tau(\langle \eta,\xi
\rangle_{\cA}), $$ where $\tau$ is the trace on the algebra
$\cA_\theta$, given by \eqref{trace-theta}.

The Yang--Mills action is defined (\cf \cite{CoRi}) as
\begin{equation} \label{YM-action}
\tau(\langle \Omega, \Omega\rangle_{\cE}).
\end{equation}
One seeks for minima of the Yang--Mills action among metric
compatible connections \eqref{metr-conn}. The main result is
that this recovers the classical moduli spaces
of Yang--Mills connections on the ordinary torus (\cite{CoRi}):

\begin{thm} \label{NCT-YM}
For a choice of a pair of integers $(p,q)$ with $p + q\theta \geq
0$, the moduli space of Yang--Mills connections on the
$\cA_\theta$ module ${\mathcal H}^\theta_{pq}$ is a classical
space given by the symmetric product
$$ s^N(T^2)= (T^2)^N / \Sigma_N, $$
where $\Sigma_N$ is the group of permutations in $N$--elements,
for $N=\gcd (p,q)$.
\end{thm}

\smallskip

In the matrix formulation of M-theory the basic equations to
obtain periodicity of two of the basic coordinates $X_i$ turn out
to be the
\begin{equation}
U_i X_j U_i^{-1} =X_j + a \delta_i^j , \ \ \  i= 1,2,  \label{Mtheory}
\end{equation}
where the $U_i$ are unitary gauge transformations.
The multiplicative commutator $U_1 U_2 U_1^{-1}U_2^{-1}$ is then
central and in the irreducible case its scalar value $ \lambda
=\exp 2\pi i \theta $ brings in the algebra of coordinates on the
noncommutative torus. The $X_j$ are then the components of the
Yang-Mills connections. The same
picture emerged from the other information one has about M-theory
concerning its relation with 11 dimensional supergravity and that
string theory dualities can then be interpreted using Morita
equivalence, relating the values of $\theta$ on an orbit of $SL_2
(\Z)$.

\smallskip

Nekrasov and Schwarz \cite{NS} showed that Yang-Mills gauge
theory on noncommutative $\R^4$ gives a conceptual understanding of the
nonzero B-field desingularization of the moduli space of
instantons obtained by perturbing the ADHM equations.
In \cite{Witten}, Seiberg and Witten exhibited the unexpected
relation between the standard gauge theory and the noncommutative
one, and clarified the limit in which the entire string dynamics
is described by a gauge theory on a noncommutative space.
Techniques from noncommutative differential and Riemannian geometry,
in the sense discussed in Section \ref{road} were applied to
string theory, for instance in \cite{Ard1}. The role of noncommutative
geometry in the context of $T$-duality was considered in very
interesting recent work of Mathai and collaborators, \cite{BMa1},
\cite{BMa2}, \cite{RoMa}.

\medskip

Recently, in the context of the holographic description of
type IIB string theory on the plane-wave background,
Shahin M.M.~Sheikh-Jabbari obtained (\cf \cite{Shahin}) an interesting
class of noncommutative spaces from the quantization of Nambu $d$-brackets.
The classical Nambu brackets
\begin{equation}\label{NambuC}
\{ f_1,\ldots, f_k \}= \epsilon^{i_1\cdots i_k} \frac{\partial
f_1}{\partial x^{i_1}}\cdots \frac{\partial
f_k}{\partial x^{i_k}}
\end{equation}
of $k$ real valued functions of variables $(x^1,\ldots, x^k)$
is quantized in the even case to the expression in $2k$
operators
\begin{equation}\label{NambuQev}
\frac{1}{i^k} [ F_1,\ldots, F_{2k}]=\frac{1}{i^k (2k)!}
\epsilon^{i_1\cdots i_{2k}} F_{i_1}\cdots F_{i_{2k}}.
\end{equation}
This generalizes the Poisson bracket quantization
$\{f_1,f_2\} \mapsto \frac{-i}{\hbar} [F_1,F_2]$.
The odd case is more subtle and it involves an additional operator
related to the chirality operator $\gamma_5$. One sets
\begin{equation}\label{NambuQodd}
\frac{1}{i^k} [ F_1,\ldots, F_{2k-1}, \gamma]=\frac{1}{i^k (2k)!}
\epsilon^{i_1\cdots i_{2k}} F_{i_1}\cdots F_{i_{2k-1}} \gamma,
\end{equation}
where $\gamma$ is the chirality operator in $2k$-dimensions.
For example, for $k=2$ one gets
$$ \begin{array}{rl}
[F_1,F_1,F_3,\gamma]=& \frac{1}{24}([F_1,F_2][F_3,\gamma]
+[F_3,\gamma][F_1,F_2] \\[2mm] & -([F_1,F_3][F_2,\gamma]+[F_2,\gamma][F_1,F_3])
\\[2mm] & +[F_2,F_3][F_1,\gamma] +[F_1,\gamma][F_2,F_3]).
\end{array} $$
If one describes the ordinary $d$-dimensional sphere of radius $R$
by the equation
\begin{equation}\label{Sd1}
\sum_{i=1}^{d+1} (x^i)^2=R^2,
\end{equation}
the coordinates satisfy
\begin{equation}\label{Sd2}
\{ x^{i_1},\ldots,x^{i_d}\}=R^{d-1} \epsilon^{i_1\cdots i_{d+1}}
x^{i_{d+1}}.
\end{equation}
The equation \eqref{Sd1} and \eqref{Sd2} are then replaced by their
quantized version, using the quantization of the Nambu bracket and
the introduction of a quantization parameter $\ell$. This defines
algebras generated by unitaries $X^i$ subject to the relations given
by the quantization of \eqref{Sd1} and \eqref{Sd2}. Matrix
representations of these algebra correspond to certain fuzzy spheres.
It would be interesting to study the general structure of these
noncommutative spaces from the point of view of the various steps
introduced in Section \ref{road}, \cf also the discussion in Section
\ref{deform}.

\bigskip

\section{Groupoids and the index theorem }
\label{cotangent}

Since the construction of the $C^*$-algebra of foliations based on
the holonomy groupoid (section \ref{foliations}), groupoids have
played a major role in noncommutative geometry. In fact the original
construction of matrix mechanics by Heisenberg (section
\ref{heisenberg}) is exactly that of the convolution algebra of the
groupoid of transitions imposed by experimental results. The
convolution algebra of groupoids can be defined in the context of
von-Neumann algebras and of $C^*$-algebras (\cf \cite{Co78} and
\cite{Ren}). It is particularly simple and canonical in the context
of smooth groupoids (\cf \cite{Co94} section II.5). One virtue of
the general construction is that it provides a geometric mental
picture of complicated analytical constructions. The prototype
example is given by the tangent groupoid of a manifold (\cf
\cite{Co94} section II.5). It is obtained by blowing up the diagonal
in the square $V\times V$ of the manifold and as a set is given by
$$
G_V=\,V\times V\times ]0,1]\cup TV
$$
where $TV$ is the (total space of the) tangent bundle of $V$. A
tangent vector $X\in T_x(V)$ appears as the limit of nearby triples
$(x_1, x_2, \varepsilon)$ provided in any chart the ratios
$(x_1-x_2)/\varepsilon$ converge to $X$. When $\varepsilon \to 0$
the Heisenberg (matrix) law of composition :
$$
(x_1, x_2, \varepsilon)\circ (x_2, x_3, \varepsilon)=(x_1, x_3,
\varepsilon)
$$
converges to the addition of tangent vectors, so that $G_V$ becomes
a smooth groupoid. The functoriality of the construction $G\to
C^*(G)$ of the convolution algebra for smooth groupoids $G$ is then
enough to define  the Atiyah-Singer analytic index of
pseudo-differential operators. It is simply given by the map in
K-theory for the exact sequence of $C^*$-algebras associated to the
geometric sequence
$$
V\times V\times ]0,1] \to G_V \supset TV
$$
where $TV$ is viewed as a closed subgroupoid of $G_V$. The
corresponding exact sequence of $C^*$-algebras can be written as
$$
0\to C_0(]0,1])\times \mathcal K\to C^*(G_V)\to C_0(T^*V)\to 0
$$
and is a geometric form of the extension of pseudo-differential
operators. By construction the algebra $C_0(]0,1])$ is {\em
contractible} and the same holds for the tensor product
$C_0(]0,1])\times \mathcal K$ by the algebra $\mathcal K$ of compact
operators. This shows that the restriction map $C^*(G_V)\to
C_0(T^*V)$ is an isomorphism in K-theory :
\begin{equation}\label{tstariso}
K_0(C_0(T^*V))\sim K_0(C^*(G_V))
\end{equation}
 Since the K-theory of
$\mathcal K$ is  $\Z$ for $K_0$, one gets the analytic index by the
evaluation map
$$
C^*(G_V)\to \mathcal K\,,\quad K_0(C^*(G_V))\to K_0(\mathcal K)=\Z
$$
composed with the isomorphism \eqref{tstariso}. Using the Thom
isomorphism yields a geometric proof (\cf \cite{Co94}) of the
Atiyah-Singer index theorem, where all the analysis has been taken
care of once and for all by the functor $G\to C^*(G)$.

This paradigm for a geometric set-up of the index theorem has been
successfully extended to manifolds with singularities (\cf
\cite{Month} \cite{Month1} and references there) and to manifolds
with boundary \cite{ANS}.

\bigskip

\section{Riemannian manifolds, conical singularities}\label{singconic}

A main property of the homotopy type of a compact
oriented manifold is that it satisfies Poincar\'e duality not
just in ordinary homology but also in $K$-homology.
In fact, while
Poincar\'e duality in ordinary homology is not sufficient to describe
homotopy type of manifolds (\cf \cite{Mi-S}), Sullivan proved (\cf
\cite{Sull}) that for simply connected PL
manifolds of dimension at least $5$, ignoring 2-torsion,
the same property in $KO$-homology suffices and the
Chern character of the $KO$-homology
fundamental class carries all the rational information on the
Pontrjagin classes.

For an ordinary manifold the choice of the fundamental cycle in
$K$-homology is a refinement of the choice of orientation of the
manifold and, in its simplest form, it is a choice of
Spin-structure. Of course the role of a spin structure is to allow
for the construction of the corresponding Dirac operator which
gives a corresponding Fredholm representation of the algebra of
smooth functions. The choice of a square root involved in the
Dirac operator $D$ corresponds to a choice of $K$-orientation.

$K$-homology admits a
fairly simple definition in terms of Hilbert spaces and Fredholm
representations of algebras. In fact,
we have the following notion of Fredholm module (\cf \cite{Co94}):

\begin{defn}\label{Fred-mod-act}
Let $\cA$ be an algebra, an {\em odd} Fredholm module over $\cA$
is given by:
\begin{enumerate}
\item a representation of $\cA$ in a Hilbert space $\cH$.
\item an operator $F = F^*$, $F^2 = 1$, on $\cH$ such that
$$
[F , a] \ \hbox{is a compact operator for any} \ a \in \cA \, .
$$
\end{enumerate}
An {\it even} Fredholm module is given by an odd Fredholm module
$(\cH , F)$ together with a $\Z / 2$ grading
$\gamma$, $\gamma = \gamma^*$,
$\gamma^2 = 1$ of the Hilbert space $\cH$ satisfying:
\begin{enumerate}
\item $\gamma a = a \gamma$, for all $a \in \cA$
\item $\gamma F = -F \gamma$.
\end{enumerate}
\end{defn}

This definition is derived from
Atiyah's definition \cite{AT} of abstract elliptic operators,
and agrees with Kasparov's definition \cite{Kasp} for the cycles
in $K$-homology, $KK(A,\C)$, when $A$ is a $C^*$-algebra.

The notion of Fredholm module can be illustrated by the following
examples (\cf \cite{Co94}).

\begin{ex}\label{Fred-manifold} {\em If $X$ is a manifold, an
elliptic operator on $X$ can be twisted with vector bundles, so as
to give rise to an index map ${\rm Ind}: K_0(C(X)) \to \Z$. If $P$
is an elliptic operator (the symbol is invertible) and a
pseudodifferential operator of order zero, $P: L^2(X,E_+)\to
L^2(X,E_-)$, then there exists a parametrix $Q$ for $P$. This is
also an operator of order zero, and a quasi-inverse for $P$, in
the sense that it is an inverse at the symbol level, namely $PQ-I$
and $QP-I$ are compact operators. Consider then the operator
$$ F=\left( \begin{array}{cc} 0 & Q \\ P & 0  \end{array} \right) $$
on $\cH=L^2(X,E_+)\oplus L^2(X,E_-)$. The algebra $C(X)$ acts on
$\cH$ and $[F,f]$ is a compact operator for all $f\in C(X)$. Since
$F^2-I$ is compact, it is possible to add to $\cH$ a finite
dimensional space to obtain $F^2=I$. Notice that the functions of
$C(X)$ act differently on this modified space. In particular the
function $f\equiv 1$ no longer acts as the identity: one recovers
the index of $P$ this way.}
\end{ex}

\begin{ex}\label{Fred-freeG} {\em Let $\Gamma=\Z * \Z$ a free
group, and let $\cA =\C \Gamma$. Let $\cT$ be the tree of $\Gamma$
with $\cT^0$ the set of vertices and $\cT^1$ the set of edges. Fix
an origin $x_0$ in $\cT^0$. For any vertex $v\in \cT^0$ there
exists a unique path connecting it to the origin $x_0$. This
defines a bijection $\phi: \cT^0 \backslash \{ x_0 \} \to \cT^1$
that assigns $v\mapsto \phi(v)$ with $\phi(v)$ this unique edge.
Let $U_\phi$ be the unitary operator implementing $\phi$, and
consider the operator
$$ F = \left( \begin{array}{cc} 0 & U_\phi \\ U_\phi^* & 0
\end{array} \right) $$
acting on $\cH = \ell^2(\cT^0) \oplus \ell^2(\cT^1) \oplus \C$. By
construction $\Gamma$ acts naturally on $\cT_j$ which gives a
corresponding action of $\cA$ in $\cH$. The pair $(\cH,F)$ is a
Fredholm module over $\cA$.}
\end{ex}

\begin{ex}\label{Hilb-transform} {\em On $S^1\simeq
\bP^1(\R)$, consider the algebra of functions ${\rm
C}(\bP^1(\R))$, acting on the Hilbert space $L^2(\R)$, as
multiplication operators $(f\, \xi)(s) = f(s) \xi(s)$. Let $F$ be
the Hilbert transform
$$ (F\, \xi)(s)=\frac{1}{\pi i} \int \frac{f(t)}{s-t} \, dt. $$
This multiplies by $+1$ the positive Fourier modes and by $-1$ the
negative Fourier modes. A function $f\in {\rm C}(\bP^1(\R))$ has
the property that $[F,f]$ is of finite rank if and only if $f$ is
a rational function $f(s) =P(s)/Q(s)$. This is Kronecker's
characterization of rational functions.}
\end{ex}

\medskip

Besides the $K$-homology class, specified by a Fredholm module,
one also wants to generalize to the noncommutative setting the
infinitesimal line element $ds$ of a Riemannian manifold. In ordinary
Riemannian geometry one deals rather with the $ds^2$ given by
the usual local expression $g_{\mu \nu} \, dx^{\mu} \, dx^{\nu}$.
However, in order to extend the notion of metric space to the
noncommutative setting it is more natural to deal with $ds$, for
which the ansatz is
\begin{equation}
ds = \times\!\!\!\!-\!\!\!-\!\!\!-\!\!\!\!\!\times \, ,
\label{dsx-x}
\end{equation}
where the right hand side has the meaning usually attributed to it in
physics, namely the Fermion propagator
\begin{equation}
\times\!\!\!\!-\!\!\!-\!\!\!-\!\!\!\!\!\times = D^{-1} ,
\label{Fermprop}
\end{equation}
where $D$ is the Dirac operator. In other words, the presence of a
spin (or spin$^c$) structure makes it possible to extract the square
root of $ds^2$, using the Dirac operator as a
differential square root of a Laplacian.

This prescription recovers the usual geodesic distance on a Riemannian
manifold, by the following result (\cf
\cite{Co-fr}).

\begin{lem}\label{geod-distance}
On a Riemannian spin manifold the geodesic distance $d(x,y)$
between two points is computed by the formula
\begin{equation}
d(x,y) = \sup  \{ \vert f(x) - f(y) \vert \, ; \ f \in \cA
\, , \ \Vert [D,f] \Vert \leq 1 \} \label{geoddist}
\end{equation}
where $D$ is the Dirac operator, $D=ds^{-1}$, and $\cA$ is the
algebra of smooth functions.
\end{lem}

This essentially follows from the fact that the quantity $\Vert
[D,f] \Vert$ can be identified with the Lipschitz norm of the
function $f$,
$$ \Vert [D,f] \Vert = {\rm ess} \sup_{x\in M} \| (\nabla f)_x \|
= \sup_{x\neq y}\frac{ | f(x)-f(y) | }{d(x,y)}. $$

Notice that, if $ds$ has the dimension of a length $L$, then $D$
has dimension $L^{-1}$ and the above expression for $d(x,y)$ also
has the dimension of a length. On a Riemannian spin manifold $X$,
the condition $\| [D,f] \|\leq 1$, for $D$ the Dirac operator, is
equivalent to the condition that $f$ is a Lipschitz function with
Lipschitz constant $c\leq 1$.

The advantage of the definition \eqref{dsx-x}, \eqref{Fermprop} of the
line element is that it is of a {\em spectral} and operator theoretic
nature, hence it extends to the noncommutative setting.

The structure that combines the $K$-homology fundamental cycle with
the spectral definition of the line element $ds$ is the notion of {\em
spectral triple} $(\cA,\cH,D)$ (\cf \cite{ConnesS3}, \cite{CoMo}).

\begin{defn}\label{def-triple}
A (compact) noncommutative geometry is a triple
\begin{equation}
(\cA ,\cH ,D)
\label{triple}
\end{equation}
where $\cA$ is a unital algebra
represented concretely as an
algebra of bounded operators on the Hilbert space $\cH$. The unbounded
operator $D$ is the inverse of the line  element
\begin{equation}
ds=1/D.\label{dsD}
\end{equation}
Such a triple $(\cA ,\cH ,D)$ is requires to satisfy the properties:
\begin{enumerate}
\item $[D,a]$ is bounded for any $a \in A^\infty$, a dense subalgebra of
the $C^*$-algebra $A=\bar\cA$.
\item $D=D^*$ and $(D+\lambda)^{-1}$ is a compact operator, for
all $\lambda \not\in \R$.
\end{enumerate}
We say that a spectral triple $(\cA ,\cH ,D)$ is even if the
Hilbert space $\cH$ has a $\Z/2$-grading by an operator $\gamma$
satisfying
\begin{equation}\label{gammaZ2}
\gamma =\gamma^*, \ \ \ \gamma^2 =1, \ \ \ \gamma\, D =- D \,
\gamma, \ \ \ \gamma \, a = a\, \gamma \ \ \forall a \in A.
\end{equation}
\end{defn}

This definition is entirely spectral. The elements of the algebra (in
general noncommutative) are operators and the line element is also an
operator. The polar decomposition $D=|D|F$ recovers the Fredholm
module $F$ defining the fundamental class in $K$-homology.
The formula for the geodesic distance extends to this
context as follows.

\begin{defn}\label{distanceNC}
Let $\varphi_i: A \to \C$, for $i=1,2$, be states on $A$, \ie
normalized positive linear functionals on $A$ with $\varphi_i (1) =
1$ and $\varphi_i (a^* a) \geq 0$ for all $a \in A$. Then the
distance between them is given by the formula
\begin{equation}
d(\varphi_1 ,\varphi_2) = \sup  \{ \vert \varphi_1 (a) -
\varphi_2 (a) \vert \ ;
\ a\in A \ , \ \Vert [D,a]\Vert \leq 1 \} \, .
\label{geodNC}
\end{equation}
\end{defn}

\smallskip

A spectral triple $(\cA,\cH,D)$ is of {\em metric dimension} $p$,
or $p$-summable, if $|D|^{-1}$ is an infinitesimal of order $1/p$
(\ie $|D|^{-p}$ is an infinitesimal of order one). Here $p<
\infty$ is a positive real number. A spectral triple $(\cA,\cH,D)$ is
$\theta$-summable if $\Tr(e^{-t D^2}) <\infty$ for all $t>0$.
The latter case corresponds to an infinite dimensional geometry.

Spectral triples also provide a more refined notion of dimension
besides the metric dimension (summability). It is given by the
{\em dimension spectrum}, which is not a number but a subset of
the complex plane.

More precisely, let  $(\cA,\cH,D)$ be
a spectral triple satisfying
the regularity hypothesis
\begin{equation}
a \ \hbox{and } \ [D,a] \ \in \ \cap_k {\rm Dom}( \delta^k ) , \ \forall \,
a\in A^\infty ,
\label{regS3}
\end{equation}
where $\delta$ is the derivation $\delta(T) = [\vert D \vert,T]$, for any
operator $T$.
Let $\cB$ denote the algebra generated by  $\delta^k (a)$ and $\delta^k
([D,a])$.
The dimension spectrum of the triple $(\cA,\cH,D)$
is the subset $\Sigma \subset \C$ consisting of all the singularities of
the analytic functions $\zeta_b (z)$ obtained by continuation of
\begin{equation}
\zeta_b (z) = \Tr  (b \vert D \vert^{-z}), \ \ \ \Re(z)
> p \ , \ \ \ b\in \cB \, . \label{zeta-b}
\end{equation}

\begin{ex} {\em
Let $M$ be a smooth compact Riemannian spin manifold, and
$(\cA,\cH,D)$ is the corresponding spectral triple given by the
algebra of smooth functions, the space of spinors, and the Dirac
operator. Then the metric dimension agrees with the
usual dimension $n$ of $M$. The dimension spectrum of $M$ is the
set $\{ 0,1,\ldots ,n\}$, where $n=\dim M$, and it is simple.
(Multiplicities appear for singular manifolds.)}
\end{ex}

It is interesting in the case of an ordinary Riemannian manifold $M$
to see the meaning of the points in the dimension spectrum that are
smaller than $n=\dim M$. These are dimensions in which the space
``manifests itself nontrivially'' with some interesting geometry.

For instance, at the point $n=\dim M$ of the dimension spectrum one
can recover the
volume form of the Riemannian metric by the equality (valid up to
a normalization constant \cf \cite{Co94})
\begin{equation}
{\int \!\!\!\!\!\! -}f \,\vert ds\vert^n =  \, \int_{M_n} f \,
\sqrt{g} \ d^n x \, ,                       \label{nvolume}
\end{equation}
where the integral $\cutint T$ is given (\cf \cite{Co94}) by the
Dixmier trace (\cf \cite{Dix}) generalizing the Wodzicki residue of
pseudodifferential operators (\cf \cite{Wo}).

One can also consider integration $\cutint ds^k$ in any other
dimension in the
dimension spectrum, with $ds=D^{-1}$ the line element.
In the case of a Riemannian manifold one finds
other important curvature expressions. For instance, if $M$ is a
manifold of dimension $\dim M=4$, when one considers integration in
dimension $2$ one finds the Einstein--Hilbert action. In fact,
a direct computation
yields the following result (\cf \cite{Kast} \cite{K-W}):

\begin{prop}
Let $dv=\sqrt{g} \ d^4 x$ denote the volume form, $ds=D^{-1}$  the
length element, \ie the inverse of the Dirac operator, and $r$ the
scalar curvature. We obtain:
\begin{equation}
{\int \!\!\!\!\!\! -} \, ds^2 = \frac{-1}{48\pi^2} \, \int_{M_4} r
\, \, dv   \, .                      \label{EHaction}
\end{equation}
\end{prop}

In general, one obtains the scalar curvature of an $n$-dimensional manifold
from the integral $\cutint ds^{n-2}$.

\medskip

Many interesting examples of spectral triples just satisfy the
conditions stated in Definition \ref{def-triple}. However, there are
significant case where more refined properties of manifolds carry over
to the noncommutative case, such as the presence of a real structure
(which makes it possible to distinguish between $K$-homology and
$KO$-homology) and the ``order one condition'' for the Dirac
operator. These properties are described as follows (\cf \cite{Co3}
and \cite{Coo2}).

\begin{defn}\label{realstr}
A real structure on an $n$-dimensional spectral triple
$(\cA,\cH,D)$ is an antilinear isometry $J: \cH \to \cH$, with the
property that
\begin{equation}\label{per8}
J^2 = \varepsilon, \ \ \ \ JD = \varepsilon' DJ, \ \ \text{and} \ \
J\gamma = \varepsilon'' \gamma
J \, \text{(even case)}.
\end{equation}
The numbers $\varepsilon ,\varepsilon' ,\varepsilon'' \in \{ -1,1\}$
are a function of $n \mod 8$ given by

\begin{center}
\begin{tabular}
{|c| r r r r r r r r|} \hline {\bf n }&0 &1 &2 &3 &4 &5 &6 &7 \\
\hline \hline
$\varepsilon$  &1 & 1&-1&-1&-1&-1& 1&1 \\
$\varepsilon'$ &1 &-1&1 &1 &1 &-1& 1&1 \\
$\varepsilon''$&1 &{}&-1&{}&1 &{}&-1&{} \\  \hline
\end{tabular}
\end{center}

Moreover, the action of $\cA$ satisfies the commutation rule
\begin{equation}\label{comm-rule}
[a,b^0] = 0 \quad \forall \, a,b \in \cA,
\end{equation}
where
\begin{equation}\label{b0}
b^0 = J b^* J^{-1} \qquad \forall b \in \cA,
\end{equation}
and the operator $D$ satisfies
\begin{equation}\label{order1}
[[D,a],b^0] = 0 \qquad \forall \, a,b \in \cA \, .
\end{equation}
\end{defn}

The anti-linear isometry $J$ is given, in ordinary Riemannian
geometry, by the charge conjugation operator acting on spinors. In
the noncommutative case, this is replaced by the Tomita
antilinear conjugation operator (\cf \cite{Tak}).

\smallskip

In \cite{Co3} and \cite{FGV} Theorem 11.2, necessary and sufficient
conditions are given that a spectral triple $(\cA,\cH,D)$ (with real
structure $J$) should fulfill in order to come from an ordinary
compact Riemannian spin manifold:

\begin{enumerate}
\item $ds=D^{-1}$ is an infinitesimal of order $1/n$.
\item There is a real structure in the sense of Definition
\ref{realstr}.
\item The commutation relation \eqref{order1} holds (this is
$[[D,a],b]=0$, for all $a,b\in \cA$, when $\cA$ is commutative).
\item The regularity hypothesis of \eqref{regS3} holds:
$a$ and $[D,a]$ are in $\cap_k {\rm Dom}( \delta^k)$ for all $a\in \cA_\infty$.
\item There exists a Hochschild cycle $c\in Z_n(\cA_\infty,\cA_\infty)$ such
that its representation $\pi(c)$ on $\cH$ induced by
$$ \pi(a^0\otimes\cdots\otimes a^n)=a^0 [D,a^1]\cdots [D,a^n] $$
satisfies $\pi(c)=\gamma$, for $\gamma$ as in \eqref{gammaZ2}, in
the even case, and $\pi(c)=1$ in the odd case.
\item The space $\cH^\infty= \cap_k {\rm Dom}( D^k)$ is a finite
projective $\cA$-module, endowed with a $\cA$-valued inner product
$\langle \xi, \eta \rangle_{\cA}$ defined by
$$ \langle a \xi, \eta \rangle = \cutint a\, \langle \xi, \eta
\rangle_{\cA}\, ds^n. $$
\item The intersection form
\begin{equation}
K_* (\cA) \times K_* (\cA) \to \Z
\label{intform}
\end{equation}
obtained from the Fredholm index of $D$ with
coefficients in $K_* (\cA \otimes \cA^0)$
is invertible.
\end{enumerate}

When $\cA$ is commutative, the above conditions characterize a
smooth Riemannian manifold $M$, with $\cA_\infty=\cC^\infty(M)$ (we
refer to \cite{FGV} for the precise statement). However, the
conditions can be stated without any commutativity assumption on
$\cA$. They are satisfied, for instance, by the isospectral
deformations of \cite{CLa}, which we discuss in Section
\ref{isodef}. Another very significant noncommutative example is the
standard model of elementary particles (\cf \cite{Co3}), which we
discuss in Section \ref{stmod}.

\bigskip

Another example of spectral triple associated to a classical
space, which is not classically a smooth manifold, is the case of
manifolds with singularities. In particular, one can consider the
case of an isolated conical singularity. This case was studied by
Lescure \cite{Lescure}.

\smallskip

Let $X$ be a manifold with an isolated conical singularity. The
cone point $c\in X$ has the property that there is a neighborhood
$U$ of $c$ in $X$, such that $U\smallsetminus \{ c \}$ is of the
form $(0, 1] \times N$, with $N$ a smooth compact manifold, and
metric $g|_U = dr^2 + r^2 g_N$, where $g_N$ is the metric on $N$.

\smallskip

A natural class of differential operators on manifolds with
isolated conical singularities is given by the elliptic operators
of Fuchs type, acting on sections of a bundle $E$. These are
operators whose restriction over $(0, 1] \times N$ takes the form
$$ r^{-\nu} \sum_{k=0}^d a_k(r) (-r \partial_r)^k, $$
for $\nu \in \R$ and $a_k \in C^\infty([0,1], {\rm Diff}^{d-k}(N,
E|_N))$, which are elliptic with symbol $\sigma_M(D)=\sum_{k=0}^d
a_k(0) z^k$ that is an elliptic family parameterized by $Im(z)$. In
particular, operators of Dirac type are elliptic of Fuchs type. For
such an operator $D$, which is of first order and symmetric, results
of Chou \cite{chou}, Br\"uning, Seeley \cite{BruSee} and Lesch
\cite{Lesch} show that its self-adjoint extension has discrete
spectrum, with $(n+1)$-summable resolvent, for $\dim X= n$.

\smallskip

The algebra that is used in the construction of the spectral
triple is $\cA= C^\infty_c(X) \oplus \C$, the algebra of functions
that are smooth on $X\smallsetminus \{ c \}$ and constant near the
singularity. The Hilbert space on which $D$ acts is chosen from a
family of weighted Sobolev spaces. Roughly, one defines weighted
Sobolev spaces that look like the standard Sobolev space on the
smooth part and on the cone are defined locally by norms
$$
\| f \|_{s,\gamma}^2 = \int_{\R^*_+ \times \R^{m-1}} \left(1 + (\log
t)^2 +\xi^2\right)^s \left| \widehat{(r^{-\gamma +1/2} f)}
(t,\xi)\right|^2 \frac{dt}{t} d\xi,
$$
where $\hat f$ denotes Fourier transform on the group $\R_+^*
\times \R^{m-1}$.

\smallskip

Then one obtains the following result (\cf \cite{Lescure})

\begin{thm}
The data $(\cA, \cH, D)$ given above define a spectral triple. In
particular, the zeta functions $\Tr(a |D|^{-z})$ admit analytic
continuation to $\C \smallsetminus \Sigma$, where the dimension
spectrum is of the form
$$ \Sigma=\{ \dim X - k, k\in \N \}, $$
with multiplicities $\leq 2$.
\end{thm}

\smallskip

The analysis of the zeta functions uses the heat kernel
$$ \Tr (| D |^{-z}) = \frac{1}{\Gamma(z/2)} \int_0^\infty t^{z/2
-1} \Tr (e^{-t D^2}) dt, $$ for which one can rely on the results of
\cite{chou} and \cite{Lesch}. The case of $\Tr(a |D|^{-z})$ of the
form $\Tr(Q |D|^{-z})$ with $Q\in \Psi^\ell_c(E)$, is treated by
splitting $Q|D|^{-z}$ as a sum of a contribution from the smooth
part and one from the singularity.

\smallskip

The Chern character for this spectral triple gives a map
$$ Ch: K_*(X) \to H_*(X,\C), $$
where we have $K^*(X)\cong K_*(\cA)$ and $H_*(X,\C) \cong
PHC^*(\cA)$, the periodic cyclic homology of the algebra $\cA$.

The cocycles $\varphi_n$ in the $(b,B)$-bicomplex for the algebra
$\cA$ have also been computed explicitly and are of the form
$$ \varphi_n (a_0,\ldots,a_n) =\nu_n \int_X a_0 da_1 \wedge \cdots
\wedge da_n \wedge \hat A(X) \wedge Ch(E),
$$
for $n \geq 1$, while for $n=0$, $\lambda\in \C$,
$\varphi_0(a+\lambda)=\int_X a \hat A(X)\wedge Ch(E) + \lambda
\Ind(D_+)$.

\bigskip

\section{Cantor sets and fractals}\label{fractals}

An important class of $C^*$-algebras are those obtained as
direct limits of a sequence of finite dimensional
subalgebras and embeddings. These are called {\em approximately finite
dimensional}, or simply AF-algebras.

An AF algebra $\cA$ is determined by a diagram of finite
dimensional algebras and inclusions, its  Bratteli diagram
\cite{Bratteli}, and from the diagram itself it is possible to
read a lot of the structure of the algebra, for instance its ideal
structure. Some simple examples of algebras that belong to this
class are:

\begin{ex}\label{Cantorset} {\it
An example of a commutative AF is the algebra of complex valued
continuous functions on a Cantor set, where a Bratteli diagram is
determined by a decreasing family of disjoint intervals covering
the Cantor set. } \end{ex}

A non--commutative example of AF algebra is given by the algebra
of the canonical anticommutation relations of quantum mechanics.

\begin{ex} {\it Consider a real Hilbert space $\cE$ and a
linear map $\cE \to \cB(\cH)$, $f \mapsto T_f$, to bounded operators
in a Hilbert space $\cH$, satisfying
$$ T_f T_g +T_g T_f =0 $$
$$ T_f^* T_g+ T_g T_f^* = \langle g,f \rangle I, $$
and the algebra $\cA$ generated by all the operators $T_f$
satisfying these relations. } \label{exCAR}
\end{ex}

A survey with many examples of AF algebras and their properties is
given for instance in \cite{Davidson}.

Let $\cA$ be a commutative AF ${\rm C}^*$-algebra. A
commutative AF algebra $\cA$ is spanned by its projections,
since finite dimensional commutative algebras are generated by
orthogonal projections. This condition is equivalent to the
spectrum $\Lambda=\Sp(\cA)$ of the algebra being a totally
disconnected compact Hausdorff space, typically a Cantor set.
Realizing such Cantor set as the intersection of a decreasing
family of disjoint intervals covering $\Lambda$ also provides a
Bratteli diagram for the AF algebra $\cA=C(\Lambda)$.

As described in \cite{Co94},
in order to construct the Hilbert space $\cH$ for a Cantor set
$\Lambda\subset \R$, let $J_k$ be the collection of bounded open
intervals in $\R\setminus \Lambda$. We denote by $L=\{ \ell_k
\}_{k\geq 1}$ the countable collection of the lengths of the
intervals $J_k$. We can assume that the lengths are ordered
\begin{equation}\label{ord-ell}
 \ell_1 \geq \ell_2 \geq \ell_3 \geq \cdots \geq \ell_k \cdots
>0.
\end{equation}
We also denote by $E=\{ x_{k,\pm} \}$ the set of the endpoints of
the intervals $J_k$, with $x_{k,+}>x_{k,-}$. Consider the Hilbert
space
\begin{equation}\label{H-endpts}
\cH := \ell^2 (E)
\end{equation}

Since the endpoints of the $J_k$ are points of $\Lambda$, there is
an action of ${\rm C}(\Lambda)$ on $\cH$ given by
\begin{equation}\label{act-CS}
f\cdot \xi (x) = f(x) \xi(x), \ \ \ \forall f\in {\rm
C}(\Lambda),\, \, \forall \xi \in \cH,\, \, \forall x \in E.
\end{equation}

A sign operator $F$ can be obtained (\cf \cite{Co94})
by choosing the closed subspace
$\hat\cH\subset \cH$ given by
\begin{equation}\label{PH}
\hat\cH=\{ \xi \in \cH : \xi(x_{k,-})=\xi(x_{k,+}), \, \, \forall
k \}.
\end{equation}
Then $F$ has eigenspaces $\hat\cH$ with eigenvalue $+1$ and
$\hat\cH^\perp$ with eigenvalue $-1$, so that, when restricted to
the subspace $\cH_k$ of coordinates $\xi(x_{k,+})$ and
$\xi(x_{k,-})$, the sign $F$ is given by
$$ F|_{\cH_k} = \left(\begin{array}{cc} 0&1\\1&0
\end{array}\right). $$

Finally, a Dirac operator $D=|D| F$ is obtained as
\begin{equation}\label{fractalDirac}
 D|_{\cH_k} \left(\begin{array}{c} \xi(x_{k,+})\\ \xi(x_{k,-})
\end{array}\right) = \ell_k^{-1} \cdot \, \left(\begin{array}{c}
\xi(x_{k,-})\\ \xi(x_{k,+}) \end{array}\right).
\end{equation}

We then obtain the following result.

\begin{prop}
Let $\Lambda\subset \R$ be a Cantor set. Let $\cA_\infty\subset
C(\Lambda)$ be the dense subalgebra of locally constant functions on
the Cantor set. Then the data $(\cA,\cH,D)$ form a spectral triple,
with $\cH$ as in
(\ref{H-endpts}), the action (\ref{act-CS}), and $D$ as in
(\ref{fractalDirac}). The zeta function satisfies
$$ \Tr( |D|^{-s} ) = 2 \zeta_L (s), $$
where $\zeta_L(s)$ is the geometric zeta function of
$L=\{ \ell_k \}_{k\geq 1}$, defined as
\begin{equation}\label{zeta-fractal}
\zeta_L(s):= \sum_k \ell_k^s.
\end{equation}
\label{propZCS}
\end{prop}

\smallskip

These zeta functions are related to the theory of Dirichlet series
and to other arithmetic zeta functions, and also to Ruelle's
dynamical zeta functions (\cf \cite{LapFra}).

\smallskip

For example, for the classical middle-third Cantor set, we have
set of lengths $\ell_k = 3^{-k}$ and multiplicities $m_k=2^{k-1}$,
for $k\geq 1$, so that we obtain
\begin{equation}\label{zeta-CS}
\Tr(|D|^{-s})= 2\zeta_L(s)=\sum_{k\geq 1} 2^k 3^{-sk} =
\frac{2\cdot 3^{-s}}{1-2\cdot 3^{-s}}.
\end{equation}
This shows that the dimension spectrum of the spectral triple of a
Cantor set has points off the real line. In fact, the set of poles
of (\ref{zeta-CS}) is
\begin{equation}\label{poles-CS}
\left\{ \frac{\log 2}{\log 3} + \frac{2\pi i n}{\log 3}
\right\}_{n\in \Z}.
\end{equation}

In this case the dimension spectrum lies on a vertical line and it
intersects the real axis in the point $D=\frac{\log 2}{\log 3}$
which is the Hausdoff dimension of the ternary Cantor set. The
same is true for other Cantor sets, as long as the self-similarity
is given by a unique contraction (in the ternary case the original
interval is replaced by two intervals of lengths scaled by 1/3).

If one considers slightly more complicated fractals in $\R$, where
the self-similarity requires more than one scaling map, the
dimension spectrum may be correspondingly more complicated. This
can be seen in the case of the Fibonacci Cantor set, for instance (\cf
\cite{LapFra}).

\begin{center}
\begin{figure}
\includegraphics[scale=0.35]{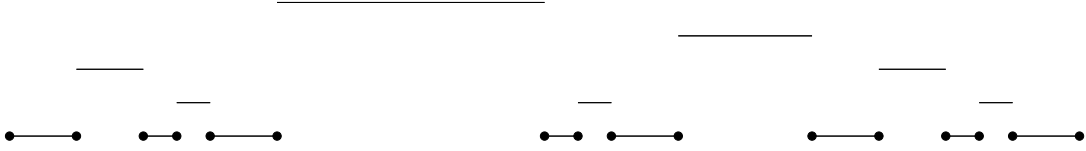}
\caption{The Fibonacci Cantor set.
\label{FigFC}}
\end{figure}
\end{center}

The Fibonacci Cantor set $\Lambda$ is obtained from the interval
$I=[0,4]$ by successively removing $F_{n+1}$ open intervals
$J_{n,j}$ of lengths $\ell_n =1/2^n$
according to the rule of Figure \ref{FigFC}.
We can associate to this Cantor set the
commutative AF algebra $\cA=C(\Lambda)$.

To obtain the Hilbert space we consider again the set $E$ of
endpoints $x_{n,j,\pm}$ of the intervals $J_{n,j}$ and we take
$\cH=\ell^2(E)$. We define the Dirac operator as in the previous
case, and we again consider the dense involutive subalgebra
$\cA_\infty$ of locally constant functions.

The data $(\cA,\cH,D)$ give a spectral triple. The Dirac
operator is related to the geometric zeta function of the
Fibonacci Cantor set by
$$ Tr(|D|^{-s})= 2 \zeta_F(s) =\frac{2}{1-2^{-s}-4^{-s}}, $$
where the geometric zeta function is $\zeta_F(s)=\sum_n F_{n+1}
2^{-ns}$, with $F_n$ the Fibonacci numbers.

A simple argument shows that the dimension
spectrum is given by the set
$$ \Sigma= \left\{
\frac{\log\phi}{\log 2} + \frac{2\pi i n}{\log 2}  \right\}_{n\in
\Z} \cup \left\{ -\frac{\log\phi}{\log 2} + \frac{2\pi i
(n+1/2)}{\log 2} \right\}_{n\in \Z}, $$ where
$\phi=\frac{1+\sqrt{5}}{2}$ is the golden ratio.

\medskip

Recent results on the noncommutative geometry of fractals and Cantor
sets and spectral triple constructions for AF algebras can be found
in \cite{AntChris}, \cite{GuIso1}, \cite{GuIso2}. The construction 
in \cite{AntChris} is in fact a spectral triple for the dual group
of the Cantor set seen as the product of countably many copies of 
the group $\Z/2$. The recent work
\cite{IvChris} shows that it is easy to describe a compact metric
space exactly (\ie recovering the metric) via a spectral triple,
which is a sum of two-dimensional modules, but spectral triples
carry much more information than just the one regarding the metric.

\bigskip

\section{Spaces of dimension $z$ and DimReg}\label{DimReg}

In perturbative quantum field theory, one computes expectation values
of observables via a formal series, where the terms are parameterized
by Feynman graphs and reduce to ordinary finite dimensional integrals
in momentum space of expressions assigned to the graphs by the
Feynman rules. These expressions typically produce divergent
integrals. For example, in the example of the scalar $\phi^3$ theory
in dimension $D=4$ or $D=4+2\N$, one encounters a divergence already
in the simplest one loop diagram, with corresponding integral (in
Euclidean signature)
\begin{equation}\label{Fint}
  \hbox{\psfig{figure=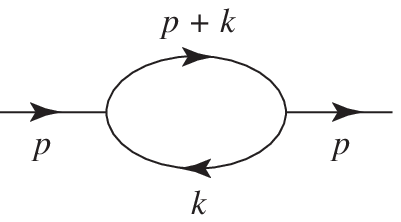}\quad =}  \int \frac{1}{k^2+m^2}\,
 \frac{1}{((p+k)^2+m^2)} \, d^D k.
\end{equation}

One needs therefore a regularization procedure for these divergent
integrals. The regularization most commonly adopted in quantum field
theory computation is ``Dimensional Regularization (DimReg) and Minimal
Subtraction (MS)''. The method was introduced in the '70s in
\cite{BoGia} and \cite{tHV} and it has the advantage of preserving
basic symmetries.

The regularization procedure of DimReg is essentially based on the use
of the formula
 \begin{equation}\label{gaussienne}
 \int\,e^{-\lambda\,k^2}\,d^{d}
k\,=\,{\pi^{d/2}}\,{\lambda^{-d/2}}\,,
 \end{equation}
to {\em define} the meaning of the integral in
$d=(D-z)$-dimensions, for $z\in \C$ in a neighborhood of zero.
For instance, in the case of \eqref{Fint}, the procedure of
dimensional regularization yields the result
$$ \pi^{(D-z)/2} \Gamma\left(\frac{4-D+z}{2}\right) \int_0^1 \left( (
x-x^2)p^2 +m^2\right)^{\frac{D-z-4}{2}} \, dx .$$

\medskip

In the recent survey \cite{ManDim}, Yuri Manin refers to DimReg as
``dimensions in search of a space''\footnote{Nicely reminiscent of
Pirandello's play ``Six characters in search of an author''.}.
Indeed, in the usual approach in perturbative quantum field theory, the
dimensional regularization procedure is just regarded as a formal rule
of analytic continuation of formal (divergent) expressions in integral
dimensions $D$ to complex values of the variable $D$.

However, using noncommutative geometry, it is possible to construct
actual spaces (in the sense of noncommutative Riemannian geometry)
$X_z$ whose dimension (in the sense of dimension spectrum) is a
point $z\in \C$ (\cf \cite{CoMar-an}).

It is well known in the physics literature that there are problems
related to using dimensionl regularization in chiral theory, because
of how to give a consistent prescription on how to extend the
$\gamma_5$ (the product of the matrices $\gamma^i$ when $D=4$) to
noninteger dimension $D-z$. It turns out that a prescription
known as Breitenlohner-Maison (\cite{BM}, \cite{Collins}) admits an
interpretation in terms of the cup product of spectral triples, where
one takes the product of the spectral triple associated to the ordinary
geometry in the integer dimension $D$ by
a spectral triple $X_z$ whose dimension spectrum is reduced
to the complex number $z$ (\cf \cite{CoMar-an}).

We illustrate here the construction for the case where $z\in
\R^*_+$. The more general case of $z\in \C$ is more delicate.

One needs to work in a slightly modified setting for spectral triples,
which is given by the type II spectral triples (\cf \cite{BF},
\cite{CPRS}, \cite{CPRS1}). In this setting the usual type I  trace
of operators in $\cL(\cH)$ is replaced by the
trace on a type II$_{\infty}$ von-Neumann algebra.

One considers a self-adjoint operator $Y$, affiliated to
a type II$_{\infty}$ factor $N$, with spectral measure given by
\begin{equation}\label{Sp-measure}
\Tr_N (\chi_E(Y))=\,\frac{1}{2}\,
\int_E\,dy
\end{equation}
for any interval $E\subset \R$, with characteristic function
$\chi_E$.

If $Y=F\, |Y|$ is the polar decomposition of $Y$, one sets
\begin{equation}\label{newdz}
D_z=\, \rho(z) \,F\, |Y|^{1/z}
 \end{equation}
with the complex
power $|Y|^{1/z}$ defined by the
functional calculus.
The normalization constant $ \rho(z)$
is chosen to be
\begin{equation}\label{rhoz}
 \rho(z)=\,\pi^{-\frac{1}{2}}\,
\left(\Gamma(\frac{z}{2}+1)\right)^{\frac{1}{z}}
\end{equation}
so that one obtains
\begin{equation}\label{basic}
\Tr \left(e^{-\lambda D^2}\right)=\,
\,{\pi^{z/2}}\,{\lambda^{-z/2}} \ \ \ \forall \lambda \in \R^*_+\,.
\end{equation}

This gives a geometric meaning to the basic formula \eqref{gaussienne}
of DimReg. The algebra $\cA$ of the spectral triple $X_z$ can be made
to contain any operator $a$ such that $[D_z,a]$ is bounded
and both $a$
and $[D_z,a]$ are smooth for the ``geodesic flow"
\begin{equation}\label{geod}
T\mapsto \,e^{it|D_z|}\,T\,e^{-it|D_z|}.
\end{equation}
The dimension spectrum of $X_z$ is reduced to the
single point $z$, since
\begin{equation}\label{trace'}
{\rm Trace}'_N(
(D^2_z)^{-s/2})=\,\rho^{-s}\, \int_1^\infty\,u^{-s/z}\,du=
\,\rho^{-s}\,\frac{z}{s-z}
\end{equation}
has a single (simple) pole
at $s=z$ and is absolutely convergent in the half space
Re$(s/z)>1$. Here ${\rm Trace}'_N$ denotes the trace with an
infrared cutoff (\ie integrating outside $|y|<1$).

\bigskip

\section{Local algebras in supersymmetric QFT}\label{susy}

It is quite striking that the general framework of noncommutative
geometry is suitable not only for handling  finite dimensional
spaces (commutative or not, of non-integer dimension etc.) but is
also compatible with infinite dimensional spaces. We already saw in
Section \ref{discgr} that discrete groups of exponential growth
naturally give rise to noncommutative spaces which are described by
a $\theta$-summable spectral triple, but not by a finitely summable
spectral triple. This is characteristic of an infinite dimensional
space and in that case, as we saw for discrete groups, cyclic
cohomology needs to be extended to {\em entire} cyclic cohomology. A
very similar kind of noncommutative spaces arises from Quantum Field
Theory in the supersymmetric context \cite{Co94} Section
IV.9.$\beta$. We   briefly recall this below and then explain
open questions also in the context of  supersymmetric theories.

The simplest example to understand the framework is that of the free
Wess-Zumino model in two dimensions, a supersymmetric free field
theory in a two dimensional space-time where space is compact
(\cite{Co94}). Thus space is a circle $S^1$ and space-time is a
cylinder $C=\,S^1\times \R$ endowed with the Lorentzian metric.  The
fields are given by a complex scalar bosonic field $\phi$ of mass
$m$ and a spinor field $\psi$ of the same mass. The Lagrangian of
the theory is of the form $\cL=\,\cL_b +\cL_f$ where,
$$
\cL_b=\,\frac{1}{2}(|\partial_0 \phi|^2- |\partial_1 \phi|^2 - m^2
|\phi|^2)
$$
and for the fermions,
$$
\cL_f=\,i\,\bar \psi\,\gamma^\mu\,\partial_\mu\,\psi-\,m\,\bar
\psi\,\psi
$$
where the spinor field is given by a column matrix, with $\bar
\psi=\,\gamma^0\,\psi^*$ and the $\gamma^\mu$ are two by two Pauli
matrices, anticommuting, self-adjoint and of square $1$.

The Hilbert space of the quantum theory is the tensor product
$\cH=\,\cH_b \otimes \cH_f$ of the bosonic one $\cH_b$ by the
fermionic one $\cH_f$. The quantum field $\phi(x)$ and its conjugate
momentum $\pi(x)$ are operator valued distributions in $\cH_b$ and
the bosonic Hamiltonian is of the form
$$
H_b=\,\int_{S^1} \,:\,|\pi(x)|^2 + |\partial_1 \phi(x)|^2 + m^2
|\phi(x)|^2  \,:\, dx
$$
where the Wick ordering takes care of an irrelevant additive
constant. The fermionic Hilbert space $\cH_f$ is given by the Dirac
sea representation which simply corresponds to a suitable spin
representation of the infinite dimensional Clifford algebra
containing the fermionic quantum fields $\psi_j(x)$. The fermionic
Hamiltonian is then the positive operator in $\cH_f$ given by
$$
H_f=\,\int_{S^1}\,:\,\bar
\psi\,\gamma^1\,i\,\partial\,\psi-\,m\,\bar \psi\,\psi \,:\,
$$
The full Hamiltonian of the non-interacting theory is acting in the
Hilbert space $\cH=\,\cH_b \otimes \cH_f$ and is the positive
operator
$$
H=\,H_b\otimes 1 +\,1\otimes H_f \,.
$$
This is were supersymmetry enters the scene in finding a
self-adjoint square root of $H$ in the same way as the Dirac
operator is a square root of the Laplacian in the case of finite
dimensional manifolds. This square root, called the {\em
supercharge} operator, is given by
$$
Q=\,\frac{1}{\sqrt{2}}\,\int_{S^1}(
\psi_1(x)(\pi(x)-\partial\phi^*(x)-i m \phi(x))+
\psi_2(x)(\pi^*(x)-\partial\phi(x)-i m \phi^*(x)) + {\rm h.c.})dx
$$
where the symbol $+{\rm h.c.}$ means that one adds the hermitian
conjugate.

The basic relation with spectral triples is then given by the
following result (\cite{Co94} Section IV).

\begin{thm} For any local region $\mathcal O\subset C$ let $\cA(\mathcal
O)$ be the algebra of functions of quantum fields with support in
$\mathcal O$ acting in the Hilbert space $\cH$. Then the triple
$$
(\cA(\mathcal O),\cH,Q)
$$
is an even $\theta$-summable spectral triple, with $\Z_2$-grading
given by the operator $\gamma=\,(-1)^{N_f}$ counting the parity of
the fermion number operator $N_f$.
\end{thm}

To be more specific the algebra $\cA(\mathcal O)$ is generated by
the imaginary exponentials $e^{i(\phi(f)+\phi(f)^*)}$ and
$e^{i(\pi(f)+\pi(f)^*)}$ for $f\in C_c^{\infty}(\mathcal O)$. As
shown in \cite{Co94} Section IV.9.$\beta$, and exactly as in the
case of discrete groups with exponential growth, one needs the
entire cyclic cohomology rather than its finite dimensional version
in order to obtain the Chern character of $\theta$-summable spectral
triples. Indeed the index map is non polynomial in the above example
of the Wess-Zumino model in two dimensions and the $K$-theory of the
above local algebras is highly non-trivial. In fact it is in that
framework that the JLO-cocycle was discovered by Jaffe-Lesniewski
and Osterwalder \cite{JLO}.

It is an open problem to extend the above result to interacting
theories in higher dimension and give a full computation of the
$K$-theory of the local algebras as well as of the Chern character
in entire cyclic cohomology. The results of Jaffe and his
collaborators on constructive quantum field theory yield many
interacting non-trivial examples of supersymmetric two dimensional
models. Moreover the recent breakthrough of Puschnigg in the case of
lattices  of semi-simple Lie groups of rank one opens the way to the
computation of the Chern character in entire cyclic cohomology.

\bigskip

\section{Spacetime and the standard model of elementary
particles}\label{stmod}

The standard model of elementary particle physics provides a
surprising example of a spectral triple in the noncommutative
setting, which in addition to the conditions of Definition
\ref{def-triple} also has a real structure satisfying all the
additional conditions of Definition \ref{realstr}.

The noncommutative geometry of the standard model developed in
\cite{Co3} (\cf also \cite{C-C}, \cite{C-C2}, \cite{CL},
\cite{Kast2}) gives a concise conceptual way to describe, through
a simple mathematical structure, the full complexity of the input
from physics. As we   recall here, the model also allows for
predictions.

The physics of the standard model can be described by a
Lagrangian. We consider here the standard model minimally coupled
to gravity, so that the Lagrangian we shall be concerned with is
the sum
\begin{equation}\label{Lagr-EHSM}
  \cL = \cL_{EH} + \cL_{SM}
\end{equation}
of the Einstein--Hilbert Lagrangian $\cL_{EH}$ and the standard
model Lagrangian $\cL_{SM}$.

The standard model Lagrangian $\cL_{SM}$ has a very complicated
expression, which, if written in full, might take a full page (\cf
\eg \cite{Velt2}). It comprises five types of terms,
\begin{equation}\label{Lagr-SM}
\cL_{SM} = \cL_G + \cL_{GH} + \cL_H + \cL_{Gf} + \cL_{Hf},
\end{equation}
where the various terms involve:
\begin{itemize}
\item spin $1$ bosons $G$: the eight gluons, $\gamma$, $W^\pm$, $Z$;
\item spin $0$ bosons $H$ such as the Higgs fields;
\item spin $1/2$ fermions $f$: quarks and leptons.
\end{itemize}
The term $\cL_G$ is the pure gauge boson part, $\cL_{GH}$ for the
minimal coupling with the Higgs fields, and $\cL_H$ gives the
quartic Higgs self interaction. In addition to the coupling
constants for the gauge fields, the fermion kinetic term
$\cL_{Gf}$ contains the hypercharges $Y_L$, $Y_R$. These numbers,
which are constant over generations, are assigned
phenomenologically, so as to obtain the correct values of the
electromagnetic charges. The term $\cL_{Hf}$ contains the Yukawa
coupling of the Higgs fields with fermions. A more detailed and
explicit description of the various terms of \eqref{Lagr-SM} is
given in \cite{Co94} \S VI.5.$\beta$. See also \cite{Velt2}.

\begin{figure}
\begin{center}
\includegraphics[scale=0.65]{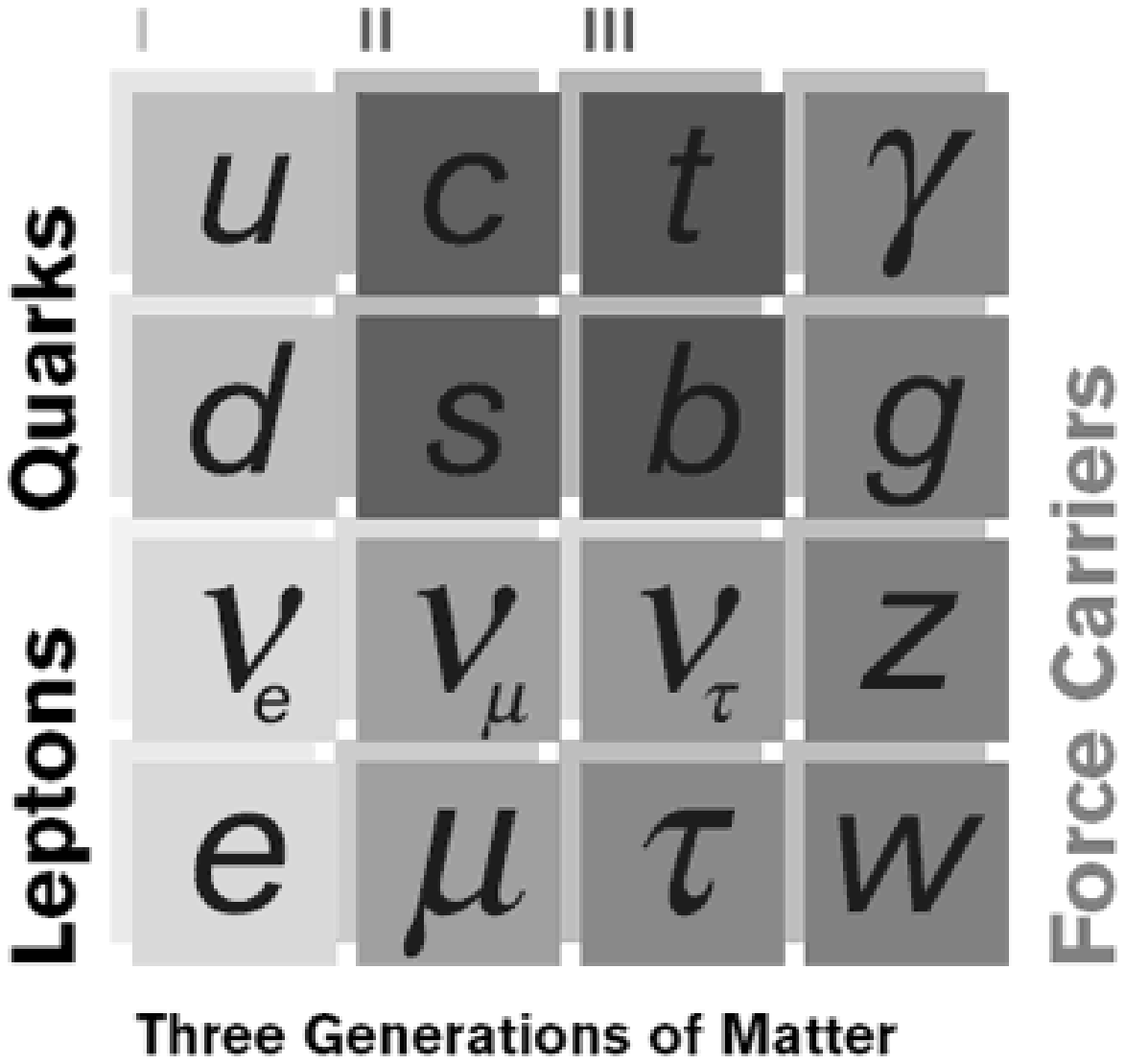}
\end{center}
\caption{Elementary particles \label{FigSM}}
\end{figure}

The symmetry group of the Einstein--Hilbert Lagrangian $\cL_{EH}$
by itself would be, by the equivalence principle, the
diffeomorphism group ${\rm Diff}(X)$ of the space-time manifold.
In the standard model Lagrangian $\cL_{SM}$, on the other hand,
the gauge theory has another huge symmetry group which is the
group of local gauge transformations. According to our current
understanding of elementary particle physics, this is given by
\begin{equation}\label{gauge-gr}
G_{SM}(X)= {\rm C}^\infty(X,U(1)\times SU(2)\times SU(3)).
\end{equation}
(At least in the case of a trivial principal bundle, \eg when the
spacetime manifold $X$ is contractible.)

Thus, when one considers the Lagrangian $\cL$ of
\eqref{Lagr-EHSM}, the full symmetry group $G$ will be a
semidirect product
\begin{equation}\label{symm-gr-SM}
G(X) = G_{SM}(X) \rtimes {\rm Diff}(X).
\end{equation}
In fact, a diffeomorphism of the manifold relabels the gauge
parameters.

To achieve a {\em geometrization} of the standard model, one would
like to be able to exhibit a space $X$ for which
\begin{equation}\label{diffG}
G(X)={\rm Diff}(X).
\end{equation}
If such a space existed, then we would be able to say that the whole
theory is pure gravity on $X$. However, it is impossible to find
such a space $X$ among ordinary manifolds. In fact, a result of W.
Thurston, D. Epstein and J. Mather (\cf \cite{mather}) shows that
the connected component of the identity in the diffeomorphism group
of a (connected) manifold is a simple group (see \cite{mather} for
the precise statement). A simple group cannot have a nontrivial
normal subgroup, so it cannot have the structure of semi-direct
product like $G(X)$ in \eqref{symm-gr-SM}.

However, it is possible to obtain a space with the desired
properties among {\em noncommutative spaces}. What plays the role of
the  connected component of the identity in the diffeomorphism group
${\rm Diff}(X)$ in the noncommutative setting is the group
$\Aut^+(\cA)$ of automorphisms of the (noncommutative) algebra that
preserve the fundamental class in $K$-homology \ie that can be
implemented by a unitary compatible with the grading and real
structure.

When the algebra $\cA$ is not commutative, among its automorphisms
there are, in particular, inner ones. They associate to an element
$x$ of the algebra the element $uxu^{-1}$, for some $u\in \cA$. Of
course $uxu^{-1}$ is not, in general, equal to $x$ because the
algebra is not commutative. The inner automorphisms form a normal
subgroup of the group of automorphisms. Thus, we see that the
group $\Aut^+ (\cA)$ has in general the same type of structure as
our desired group of symmetries $G(X)$, namely, it has a normal
subgroup of inner automorphisms and it has a quotient. It is
amusing how the physical and the mathematical vocabularies agree
here: in physics one talks about internal symmetries and in
mathematics one talks about inner automorphisms (one might as well
call them internal automorphisms).

There is a very simple non commutative algebra $\cA$ whose group
of inner automorphisms corresponds to the group of gauge
transformations $G_{SM}(X)$, and such that the quotient $\Aut^+
(\cA) / {\rm Inn} (\cA)$ corresponds exactly to diffeomorphisms
(\cf \cite{Schu}). The noncommutative space is a product $X\times
F$ of an ordinary spacetime manifold $X$ by a ``finite
noncommutative space'' $F$. The noncommutative algebra $\cA_F$ is
a direct sum of the algebras $\bC$, $\bH$ (here denoting the
quaternions), and $M_3 (\bC)$ (the algebra of $3\times 3$ complex
matrices).

The algebra $\cA_F$ corresponds to a {\it finite} space where the
standard model fermions and the Yukawa parameters (masses of
fermions and mixing matrix of Kobayashi Maskawa) determine the
spectral geometry in the following manner. The Hilbert space
$\cH_F$ is finite-dimensional and admits the set of elementary
fermions as a basis. This comprises the generations of quarks
(down--up, strange--charmed, bottom--top),
\begin{equation}
\begin{array}{ccccccccc}
u_L & u_R & d_L & d_R& & \bar{u}_L & \bar{u}_R & \bar{d}_L & \bar{d}_R \\[2mm]
c_L & c_R & s_L & s_R& & \bar{c}_L & \bar{c}_R & \bar{s}_L & \bar{s}_R \\[2mm]
t_L & t_R & b_L & b_R& & \bar{t}_L & \bar{t}_R & \bar{b}_L &
\bar{b}_R,
\end{array}
\label{quarks}
\end{equation}
with the additional color index $(y,r,b)$, and the generations of
leptons (electron, muon, tau, and corresponding neutrinos)
\begin{equation}
\begin{array}{ccccccc}
e_L& e_R& \nu^e_L& & \bar e_L& \bar e_R& \bar \nu^e_L
\\[2mm]
\mu_L& \mu_R& \nu^\mu_L& & \bar \mu_L& \bar \mu_R& \bar \nu^\mu_L
\\[2mm]
\tau_L& \tau_R& \nu^\tau_L& & \bar \tau_L& \bar \tau_R &
\bar\nu^\tau_L
\end{array}
\label{leptons}
\end{equation}
(We discuss here only the minimal standard model with no right
handed neutrinos.)

The $\Z/2$ grading $\gamma_F$ on the Hilbert space $\cH_F$ has
sign $+1$ on left handed particles (\eg the $u_L$, $d_L$, etc.)
and sign $-1$ on the right handed particles. The involution $J_F$
giving the real structure is the charge conjugation, namely, if we
write $\cH_F = \cE \oplus \bar \cE$, then $J_F$ acts on the
fermion basis as $J_F(f,\bar h) = (h,\bar f)$. This satisfies
$J_F^2=1$ and $J_F\gamma_F=\gamma_F J_F$, as should be for
dimension $n=0$.

The algebra $\cA_F$ admits a natural representation in $\cH_F$
(see \cite{Coo2}). An element $(z,q,m)\in \C \oplus \H \oplus
M_3(\C)$ acts as
$$ (z,q,m) \cdot \left(\begin{array}{c} u_R \\ d_R \end{array}\right)
= \left(\begin{array}{c} z\, u_R \\ \bar z\, d_R
\end{array}\right) \ \ \ \ (z,q,m) \cdot e_R = \bar z\,  e_R $$
$$ (z,q,m) \cdot \left(\begin{array}{c} u_L \\ d_L \end{array}\right)
= q \left(\begin{array}{c} u_L \\ d_L \end{array}\right) \ \ \ \
(z,q,m) \cdot \left(\begin{array}{c} \nu^e_L \\ e_L
\end{array}\right) = q \left(\begin{array}{c} \nu^e_L \\ e_L
\end{array}\right), $$
$$ (z,q,m) \cdot \left(\begin{array}{c} \bar e_L \\ \bar e_R
\end{array}\right) = \left(\begin{array}{c} z\, \bar e_L \\ z\, \bar e_R
\end{array}\right) $$
$$ (z,q,m) \cdot \bar u_R  = m\,  \bar u_R \ \ \ \ (z,q,m) \cdot \bar
d_R =m  \, \bar  d_R $$ and similarly for the other generations.
Here $q\in \H$ acts as multiplication by the matrix
$$ q= \left(\begin{array}{cc} \alpha & \beta \\
-\bar \beta & \bar \alpha \end{array}\right),  $$ where $q=\alpha+
\beta\, j$, with $\alpha,\beta\in \C$. The matrix $m\in M_3(\C)$
acts on the color indices $(y,r,b)$.

The data $(\cA_F, \cH_F)$ can be completed to a spectral triple
$(\cA_F,\cH_F,D_F)$ where the Dirac operator (in this finite
dimensional case a matrix) is given by
\begin{equation}\label{YukawaD}
D_F = \left( \begin{array}{cc} Y & 0 \\ 0 & \bar Y
\end{array}\right)
\end{equation}
on $\cH_F=\cE\oplus \bar \cE$, where $Y$ is the Yukawa coupling
matrix, which combines the masses of the elementary fermions
together with the Cabibbo--Kobayashi--Maskawa (CKM) quark mixing matrix.

The fermionic fields acquire mass through the spontaneous symmetry
breaking produced by the Higgs fields. The Yukawa coupling matrix
takes the form $Y=Y_q \otimes 1 \oplus Y_f$, where the matrix
$Y_f$ is of the form
$$ \left(\begin{array}{ccc} 0&0& M_e\\
0&0&0\\ M_e^* &0&0 \end{array}\right), $$ in the basis
$(e_R,\nu_L,e_L)$ and successive generations, while $Y_q$ is of
the form
$$ \left(\begin{array}{cccc} 0&0& M_u & 0\\
0&0&0& M_d \\M_u^* &0&0&0 \\ 0& M_d^* & 0&0
\end{array}\right), $$
in the basis given by $(u_R,d_R,u_L,d_L)$ and successive
generations. In the case of the lepton masses, up to rotating the
fields to mass eigenstates, one obtains a mass term for each
fermion, and the off diagonal terms in $M_e$ can be reabsorbed in
the definition of the fields. In the quark case, the situation is
more complicated and the Yukawa coupling matrix can be reduced to
the mass eigenvalues and the CKM quark mixing. By rotating the
fields, it is possible to eliminate the off diagonal terms in
$M_u$. Then $M_d$ satisfies $V M_d V^*=M_u$, where $V$ is the CKM
quark mixing, given by a $3\times 3$ unitary matrix
$$ V= \left( \begin{array}{ccc} V_{ud} & V_{us} & V_{ub} \\
V_{cd} & V_{cs} & V_{cb} \\
V_{td} & V_{ts} & V_{tb} \end{array} \right) $$ acting on the
charge $-e/3$ quarks (down, strange, bottom). The entries of this
matrix can be expressed in terms of three angles
$\theta_{12},\theta_{23},\theta_{13}$ and a phase, and can be
determined experimentally from weak decays and deep inelastic
neutrino scatterings.

The detailed structure of the Yukawa coupling matrix $Y$ (in
particular the fact that color is not broken) allows one to check
that the finite geometry $(\cA_F,\cH_F,D_F)$ satisfies all the
axioms of Definition \ref{realstr} for a noncommutative spectral
manifold. The key point is that elements $a \in \cA_F$ and
$[D_f,a]$ commute with $J_F\, \cA_F\, J_F$. These operators
preserve the subspace $\cE \subset \cH_F$. On this subspace, for
$b=(z,q,m)$, the action of $J_F \, b^* \, J_F$ is by
multiplication by $z$ or by the transpose $m^t$. It is then not
hard to check explicitly the commutation with $a$ or $[D,a]$ (\cf
\cite{Co94} \S VI.5.$\delta$). By exchanging the roles of $a$ and
$b$, one sees analogously that $a$ commutes with $J_F b J_F$ and
$[D,J_F b J_F]$ on $\bar \cE$, hence the desired commutation
relations hold on all of $\cH_F$.

We can then consider the product $X\times F$, where $X$ is an
ordinary $4$--dimensional Riemannian spin manifold and $F$ is the
finite geometry described above. This product geometry is a
spectral triple $(\cA,\cH,D)$ obtained as the cup product of a
triple $(\cC^\infty(X),L^2(X,S), D_1)$, where $D_1$ is the Dirac
operator on $X$ acting on square integrable spinors in $L^2(X,S)$,
with the spectral triple $(\cA_F,\cH_F,D_F)$ described above.
Namely, the resulting (smooth) algebra and Hilbert space are of
the form
\begin{equation}\label{algXF}
\cA_\infty=\cC^\infty(X,\cA_F) \ \ \ \  \cH= L^2(X,S)\otimes
\cH_F,
\end{equation}
and the Dirac operator is given by
\begin{equation}\label{diracXF}
D = D_1 \otimes 1 + \gamma \otimes D_F,
\end{equation}
where $\gamma$ is the usual $\Z/2$ grading on the spinor bundle
$S$. The induced $\Z/2$ grading on $\cH$ is the tensor product
$\gamma\otimes \gamma_F$, and the real structure is given by
$J=C\otimes J_F$, where $C$ is the charge conjugation operator on
spinors.

Notice that, so far, we have only used the information on the
fermions of the standard model. We'll see now that the bosons,
with the correct quantum numbers, are {\em deduced} as inner
fluctuations of the metric of the spectral triple $(\cA,\cH,D)$.

It is a general fact that, for noncommutative geometries
$(\cA,\cH,D)$, one can consider inner fluctuations of the metric
of the form
$$
D\mapsto D+A+JAJ^{-1}
$$
where $A$ is of the form
\begin{equation}\label{Aab}
A= \sum a_i \, [D, a_i' ]  \ \ \ \  a_i,a_i' \in \cA.
\end{equation}

In the case of the standard model, a direct computation of the
inner fluctuations gives the standard model gauge bosons $\gamma ,
W^{\pm} ,Z$, the eight gluons and the Higgs fields $\varphi$ with
accurate quantum numbers (\cf \cite{Co3}). In fact, a field $A$ of
the form \eqref{Aab} can be separated in a ``discrete part''
$A^{(0,1)} = \sum a_i \, [\gamma\otimes D_F, a_i' ]$ and a
continuous part $A^{(1,0)} = \sum a_i \, [D_1\otimes 1, a_i' ]$,
with $a_i=(z_i,q_i,m_i)$ and $a_i'=(z_i',q_i',m_i')$,
$q_i=\alpha_i + \beta_i j$ and $q_i'=\alpha_i' + \beta_i' j$. The
discrete part gives a quaternion valued function
$$ q(x) = \sum z_i \left( (\alpha_i'-z_i') + z_i \beta_i'\, j \right)
= \varphi_1 + \varphi_2 \, j  $$ which provides the Higgs doublet.
The continuous part gives three types of fields:
\begin{itemize}
\item A $U(1)$ gauge field $U =\sum z_i \, d z_i'$
\item An $SU(2)$ gauge field $Q= \sum q_i \, d q_i'$
\item A $U(3)$ gauge field $M=\sum m_i \, d m_i'$,
which can be reduced to an $SU(3)$ gauge field $M'$ by subtracting
the scalar part of the overall gauge field which eliminates
inessential fluctuations that do not change the metric.
\end{itemize}
The resulting internal fluctuation of the metric $A+J A J^{-1}$ is
then of the form (\cf \cite{Co3})
$$ \left(\begin{array}{ccc} -2U & 0 & 0 \\
0 & Q_{11}-U & Q_{12} \\ 0& Q_{21} & Q_{22}-U
\end{array} \right) $$ on the basis of leptons $(e_R,\nu_L,e_L)$
and successive generations, and
$$ \left(\begin{array}{cccc} \frac{4}{3} U + M' & 0 & 0 & 0 \\[2mm]
0 & \frac{-2}{3} U + M' & 0 & 0 \\[2mm]
0 & 0 & Q_{11} + \frac{1}{3} U + M' & Q_{12} \\[2mm]
0 & 0 & Q_{21} & Q_{22} + \frac{1}{3} U + M'
\end{array}\right), $$
on the basis of quarks given by $(u_R,d_R,u_L,d_L)$ and successive
generations. A striking feature that these internal fluctuations
exhibit is the fact that the expressions above recover all the
exact values of the hypercharges $Y_L$, $Y_R$ that appear in the
fermion kinetic term of the standard model Lagrangian.

Finally, one can also recover the bosonic part of the standard model
Lagragian from a very general principle, the {\em spectral action
principle} of Chamseddine--Connes (\cf \cite{C-C}, \cite{C-C2},
\cite{C-C3}). The result is that the Hilbert--Einstein action
functional for the Riemannian metric, the Yang--Mills action for the
vector potentials, and the self interaction and the minimal coupling
for the Higgs fields all appear with the correct signs in the
asymptotic expansion for large $\Lambda$ of the number $N(\Lambda)$
of eigenvalues of $D$ which are $\leq \Lambda$ (\cf \cite{C-C}),
\begin{equation}
N(\Lambda) = \# \ \hbox{eigenvalues of $D$ in} \
[-\Lambda,\Lambda].
\end{equation}

The spectral action principle, applied to a spectral triple
$(\cA,\cH,D)$, can be stated as saying that the physical action
depends only on $\Sp(D)\subset \R$. This spectral datum corresponds
to the data $(\cH,D)$ of the spectral triple, independent of the
action of $\cA$. Different $\cA$ that correspond to the same
spectral data can be thought of as the noncommutative analog of
isospectral Riemannian manifolds (\cf the discussion of isospectral
deformations in Section \ref{isodef}). A natural expression for an
action that depends only on $\Sp(D)$ and is {\em additive} for
direct sums of spaces is of the form
\begin{equation}\label{Sp-action}
\Tr\, \chi (\frac{D}{\Lambda} )+ \langle \psi, D \, \psi \rangle,
\end{equation}
where $\chi$ is a positive even function and $\Lambda$ is a scale.

In the case of the standard model, this formula \eqref{Sp-action} is
applied to the full ``metric" including the internal fluctuations
and gives the full standard model action minimally coupled with
gravity. The Fermionic part of the action \eqref{Sp-action} gives
(\cf \cite{C-C}, \cite{C-C2})
\begin{equation}\label{fermionic}
\langle \psi, D \, \psi \rangle= \int_X \left( \cL_{Gf} + \cL_{Hf}
\right) \sqrt{|g|} d^4 x .
\end{equation}

The bosonic part of  the action \eqref{Sp-action} evaluated via heat
kernel invariants gives the standard model Lagrangian minimally
coupled with gravity. Namely, one writes the function $\chi
(\frac{D}{\Lambda})$  as the superposition of exponentials. One then
computes the trace by a semiclassical approximation from local
expressions involving the familiar heat equation expansion. This
delivers all the correct terms in the action (\cf \cite{C-C2} for an
explict calculation of all the terms involved).

Notice that here one treats the spacetime manifold $X$ in the
Euclidean signature. The formalism of spectral triple can be
extended in various ways to the Lorentzian signature (\cf \eg
\cite{EHaw}). Perhaps the most convenient choice is to drop the
self-adjointness condition for $D$ while still requiring $D^2$ to
be self-adjoint.

\bigskip

\section{Isospectral deformations}\label{isodef}

A very rich class of examples of noncommutative manifolds is
obtained by considering isospectral deformations of a classical
Riemannian manifold. These examples satisfy all the axioms of
ordinary Riemannian geometry (\cf \cite{Co3}) except
commutativity. They are obtained by the following result
(Connes--Landi \cite{CLa}):

\begin{thm}\label{Theorem6}
Let $M$ be a compact Riemannian spin manifold. Then if the
isometry group of $M$ has rank $r \geq 2$, $M$ admits a
non-trivial one parameter isospectral deformation to
noncommutative geometries $M_{\theta}$.
\end{thm}

\smallskip

The main idea of the construction is to deform the standard
spectral triple describing the Riemannian geometry along a two
torus embedded in the isometry group, to a family of spectral
triples describing non-commutative geometries.

\smallskip

More precisely, under the assumption on the rank of the group of
isometries of the compact spin manifold $X$, there exists a
two-torus
$$ T^2 \subset {\rm Isom}(X), $$
where we identify $T^2 = \R^2 / (2\pi \Z)^2$. Let $U(s)$ be the
unitary operators in this subgroup of isometries, for
$s=(s_1,s_2)\in T^2$, acting on the Hilbert space $\cH = L^2(X,S)$
of the spectral triple $$(C^\infty(X),L^2(X,S),D,J).$$
Equivalently, we write $U(s)= \exp (i(s_1 P_1 + s_2 P_2))$, where
$P_i$ are the corresponding Lie algebra generators, with ${\rm
Spec}(P_i)\subset \Z$, satisfying $[D,P_i]=0$ and $P_i J = - J
P_i$, so that $[U(s),D]=[U(s),J]=0$.

\smallskip

The action $\alpha_s(T)=U(s) T U(s)^{-1}$ has the following
property. Any operator $T$ such that the map $s \mapsto
\alpha_s(T)$ is smooth can be uniquely written as a norm
convergent series
\begin{equation}\label{T-bidegree}
T = \sum_{n_1, n_2 \in \Z} \hat T_{n_1,n_2}
\end{equation}
where each term $\hat T_{n_1,n_2}$ is an operator of bi-degree
$(n_1,n_2)$, that is,
$$ \alpha_s(\hat T_{n_1,n_2})= \exp(i(s_1 n_1 + s_2 n_2)) \hat
T_{n_1,n_2}, $$ for each $s=(s_1,s_2)\in T^2$, and the sequence of
norms $\| \hat T_{n_1,n_2} \|$ is of rapid decay.

\smallskip

This property makes it possible to define left and right twists
for such operators $T$, defined as
\begin{equation}\label{T-left-twist}
\ell (T):= \sum_{n_1,n_2} \hat T_{n_1,n_2} \exp\left( 2 \pi i
\theta n_2 P_1 \right)
\end{equation}
and
\begin{equation}\label{T-right-twist}
r (T):=\sum_{n_1,n_2} \hat T_{n_1,n_2} \exp\left( 2 \pi i \theta
n_1 P_2 \right).
\end{equation}
Both series still converge in norm, since the $P_i$ are
self-adjoint operators.

\smallskip

It is then possible to introduce a (left) deformed product
\begin{equation}\label{left-twist-prod}
x * y = \exp ( 2 \pi i \theta n_1 ' n_2 ) xy,
\end{equation}
for $x$ a homogeneous operator of bi-degree $(n_1, n_2)$ and $y$ a
homogeneous operator of bi-degree $(n_1', n_2')$. A (right)
deformed product is similarly defined by setting $x *_r y = \exp (
2 \pi i \theta n_1  n_2' ) xy$. These deformed products satisfy
$\ell(x)\ell(y)= x* y$ and $r(x)r(y)=x *_r y$.

\smallskip

The deformed spectral triples are then obtained by maintaining the
same Hilbert space $\cH = L^2(X,S)$ and Dirac operator $D$, while
modifying the algebra $C^\infty(X)$ to the non-commutative algebra
$\cA_\theta :=\ell(C^\infty(X))$ and the involution $J$ that
defines the real structure to $J_\theta := \exp(2\pi i \theta
P_1P_2) J$.

\bigskip

\section{Algebraic deformations}\label{deform}

There is a very general context in which one constructs noncommutative
spaces via deformations of commutative algebras. Unlike the
isospectral deformations discussed in Section \ref{isodef}, here one
proceeds mostly at a formal algebraic level, without involving the
operator algebra structure and without invoking the presence of a
Riemannian structure.

The idea of deformation quantization originates in the idea that
classical mechanics has as setting a smooth manifold (phase space)
with a symplectic structure, which defines a Poisson bracket
$\{,\}$. The system is quantized by deforming the pointwise
product in the algebra $\cA=\cC^\infty(M)$ (or in a suitable
subalgebra) to a family $*_\hbar$ of products satisfying $f
*_\hbar g \to fg$ as $\hbar\to 0$, which are associative but no
longer necessarily commutative. These are also required to satisfy
$$ \frac{f *_\hbar g - g*_\hbar f}{i\hbar} \to \{ f, g \}, $$
as $\hbar\to 0$, namely, the ordinary product is deformed in the
direction of the Poisson bracket. On the algebra $\cC^\infty(M)$ a
Poisson bracket is specified by assigning a section $\Lambda$ of
$\Lambda^2(TM)$ with the property that
$$ \{ f, g \} = \langle \Lambda, df\wedge dg \rangle $$
satisfies the Jacobi identity. Typically, this produces a {\it
formal deformation}: a formal power series in $\hbar$. Namely, the
deformed product can be written in terms of a sequence of
bi--differential operators $B_k$ satisfying
\begin{equation}
f*g = fg + \hbar B_1(f,g) + \hbar^2 B_2(f,g) + \cdots
\label{Poisson-deform}
\end{equation}

\smallskip

Under this perspective, there is a good understanding of formal
deformations. For instance, Kontsevich \cite{KontIHES97} proved
that formal deformations always exist, by providing an explicit
combinatorial formula that generates all the $\{ B_2, B_3, \ldots
\}$ in the expansion from the $B_1$, hence in terms of the Poisson
structure $\Lambda$. The formal solution (\ref{Poisson-deform})
can then be written as
$$ \sum_{n=0}^\infty \hbar^n \sum_{\Gamma\in G[n]} \omega_\Gamma
B_{\Gamma,\Lambda}(f,g), $$ where $G[n]$ is a set of $(n(n+1))^n$
labeled graphs with $n+2$ vertices and $n$ edges, $\omega_\Gamma$
is a coeffcient obtained by integrating a differential form
(depending on the graph $\Gamma$) on the configuration space of
$n$ distinct points in the upper half plane, and
$B_{\Gamma,\Lambda}$ is a bi--differential operator whose
coefficients are derivatives of $\Lambda$ of orders specified by
the combinatorial information of the graph $\Gamma$.

\smallskip

A setting of deformation quantization which is compatible with
$C^*$--algebras was developed by Rieffel in \cite{Rieffel89}. We
recall briefly Rieffel's setting. For simplicity, we restrict to
the simpler case of a compact manifold.

\begin{defn}\label{Rieffel-def}
A strict (Rieffel) deformation quantization of $\cA=\cC^\infty(M)$
is obtained by assigning an associative product $*_\hbar$, an
involution (depending on $\hbar$) and a $C^*$--norm $\| \cdot
\|_{\hbar}$ on $\cA$, for $\hbar \in I$ (some interval containing
zero), such that:

(i) For $\hbar=0$ these give the $C^*$--algebra $C(M)$,

(ii) For all $f,g\in \cA$, as $\hbar \to 0$,
$$ \left\| \frac{f *_\hbar g - g*_\hbar f}{i\hbar} - \{ f, g \}
\right\|_{\hbar} \to 0. $$

One denotes by $\cA_\hbar$ the $C^*$--algebra obtained by
completing $\cA$ in the norm $\| \cdot \|_{\hbar}$.
\end{defn}

The functions of $\hbar$ are all supposed to be analytic, so that
formal power series expansions make sense.

\begin{rem}{\em The notion of a strict deformation quantization should be
regarded as a notion of integrability for a formal solutions. }
\end{rem}

\smallskip

Rieffel also provides a setting for compatible actions by a Lie
group of symmetries, and proves that non--commutative tori (also
of higher rank) are strict deformation quantizations of ordinary
tori, that are compatible with the action of the ordinary torus as
group of symmetry. Typically, for a given Poisson structure,
strict deformation quantizations are not unique. This happens
already in the case of tori.

\smallskip

In the same paper \cite{Rieffel89}, Rieffel uses a basic result of
Wassermann \cite{Wassermann88} to produce an example where formal
solutions are not integrable. The example is provided by the
two--sphere $S^2$. There is on $S^2$ a symplectic structure, and a
corresponding Poisson structure $\Lambda$ which is invariant under
$SO(3)$. Rieffel proves the following striking result (Theorem 7.1
of \cite{Rieffel89}):

\begin{thm} There are no $SO(3)$--invariant strict deformations of the
ordinary product on $\cC^\infty(S^2)$ in the direction of the
$SO(3)$--invariant Poisson structure.
\end{thm}

In fact, the proof of this result shows more, namely that no
$SO(3)$--invariant deformation of the ordinary product in $C(S^2)$
can produce a non--commutative $C^*$--algebra. This rigidity
result reflects a strong rigidity result for $SU(2)$ proved by
Wassermann \cite{Wassermann88}, namely the only ergodic actions of
$SU(2)$ are on von Neumann algebras of type I. The interest of
this result lies in the fact that there are formal deformations of
the Poisson structure that are $SO(3)$--invariant (see \eg
\cite{GuVai}, \cite{BayFron}), but these only exist as a formal
power series in the sense of (\ref{Poisson-deform}) and, by the
results of Wassermann and Rieffel are not integrable.

\smallskip

Summarizing, we have the following type of phenomenon: on the one
hand we have formal solutions, formal deformation quantizations
about which a lot is known, but for which, in general, there may
not be an integrability result. More precisely, when we try to
pass from formal to actual solutions, there are cases where
existence fails (the sphere), and others (tori) where uniqueness
fails. The picture that emerges is remarkably similar to the case
of formal and actual solutions of ordinary differential equations.

\smallskip

It is very instructive to build an analogy between the {\it problem
of ambiguity} for formal solutions of ODE's and the present
situation of formal non--commutative spaces and actual
non--commutative spaces. The main conclusion to be drawn from this
analogy is that there ought to be a {\it theory of ambiguity} which
formulates precisely the relation between the formal
non--commutative geometry and its integrated ($C^*$--algebraic)
version.

\smallskip

To illustrate this concept, we take a closer look at the analogous
story in the theory of ODE's. A good reference for a modern
viewpoint is \cite{Ramis}. A formal solution of a differential
equation is a power series expansion: for instance
$\sum_{n=0}^\infty (-1)^n n! x^{n+1}$ is a formal solution of the
Euler equation $x^2 y' +y=x$. Convergent series give rise to actual
solutions, and more involved summation processes such as Borel
summation can be used to transform  a given formal solution of an
analytic ODE into an actual solution on a sufficiently narrow sector
in $\C$ of sufficiently small radius, but such solution is in
general {\it not unique}. It is known from several classical methods
that some divergent series can  be ``summed'' modulo a function with
exponential decrease of a certain order. This property (Gevrey
summability) is also satisfied by formal solutions of analytic
ODE's, and, stated in a more geometric fashion, it is essentially a
cohomological condition. It also shows that, whereas on {\it small
sectors} one has existence of actual solutions but not uniqueness,
on {\it large sectors} one gains uniqueness, at the cost of possibly
loosing existence. A complete answer to summability of formal
solutions can then be given in terms of a more refined
multi--summability (combining Gevrey series and functions of
different order) and the Newton polygon of the equation.

\smallskip

The general flavor of this theory is surprisingly similar to the
problem of formal solutions in non--commutative geometry. It is to
be expected that an ambiguity theorem exists, which accounts for
the cases of lack of uniqueness, or lack of existence, of actual
solutions illustrated by the results of Rieffel.

\medskip

Already in dealing with our first truly non--trivial example of
noncommutative spaces, the noncommutative tori, we encountered
subtleties related to the difference between the quotient and the
deformation approach to the construction of non--commutative
spaces.

\smallskip

In fact, the non--commutative tori we described in Section
\ref{nctori} admit a description as algebras
obtained as deformations of the ordinary product of functions, by
setting
\begin{equation}\label{NCTdeform}
 (f*g)(x,y):= \left( e^{2\pi i\theta \frac{\partial}{\partial x}
\frac{\partial}{\partial y'} } f(x,y) g(x',y') \right)_{x=x',y=y'} =
\sum \frac{(i2\pi \theta)^n}{n!} D_1^n f D_2^n g.
\end{equation}

\smallskip

Notice however that while $U\,\frac{\partial}{\partial U}$ and
$V\,\frac{\partial}{\partial V}$ are derivations
 for the algebra of the non--commutative torus, this is {\it not}
 the case for
$\frac{\partial}{\partial U}$ and $\frac{\partial}{\partial V}$. The
same holds for the quantum plane (\cf \cite{[M]}) whose algebra of
coordinates admits two generators $u,v$ with relation
$$
u\,v=\,q\,v\,u\,.
$$
These generators can be rotated ($u\mapsto \lambda
\,u$, $v\mapsto \mu \,v$) without affecting the presentation but
translations of the generators are not automorphisms of the algebra.
 In other words, one can view the
non--commutative torus as a deformation of an ordinary torus, which
in turn is a quotient of the classical plane $\R^2$ by a lattice of
translations, but the action of translations does not extend to the
quantum plane. This is an instance of the fact that the general
operations of quotient and deformation, in constructing
non--commutative spaces, do not satisfy any simple compatibility
rules and need to be manipulated with care.

Moreover, phenomena like the Morita equivalence between, for
instance $\theta$ and $1/{\theta}$, are not  detectable in a purely
deformation theoretic perturbative expansion like the one given by
the Moyal product \eqref{NCTdeform}. They are non-perturbative and
cannot be seen at the perturbative level of the star product.

\medskip

In this respect, a very interesting recent result is that of Gayral,
Gracia-Bondia, Iochum, Sch\"ucker, and Varilly, \cite{GGISV}, where
they consider a version of the structure of spectral triple for
non-compact spaces. In that case, for instance, one no longer can
expect the Dirac operator to have compact resolvent and one can only
expect a local version to hold, \eg $a (D-i)^{-1}$ is compact for
$a\in \cA$. Other properties of Definitions \ref{def-triple} and
\ref{realstr} are easily adapted to a ``local version'' but become
more difficult to check than in the compact case. They show that the
Moyal product deformation of $\R^{2n}$ fits in the framework of
spectral triples and provides an example of such non-unital spectral
triples. Thus, it appears that the structure of noncommutative
Riemannian geometry provided by spectral triples should adapt nicely
to some classes of algebraic deformations.

\smallskip

It appears at first that spectral triples  may not be the right type
of structure to deal with noncommutative spaces associated to
algebraic deformations, because it corresponds to a form of
Riemannian geometry, while many such spaces originate from K\"ahler
geometry. However, the K\"ahler structure can often be also encoded
in the setting of spectral triple, for example by considering also a
second Dirac operator, as in \cite{BKLR} or through the presence of
a Lefschetz operator as in \cite{CM1}.

\medskip

Noncommutative spaces obtained as deformations of commutative
algebras fit in the context of a well developed algebraic theory of
noncommutative spaces (\cf \eg \cite{Kontsevich93}
\cite{Kontsevich00} \cite{Manin01RM} \cite{Manin91NC} \cite{[M]}
\cite{Rosenberg98} \cite{Rosenberg95}, \cite{Soibelman01}). This
theory touches on a variety of subjects like quantum groups and the
deformation approach to non--commutative spaces and is interestingly
connected to the theory of mirror symmetry. However, it is often not
clear how to integrate this approach with the functional analytic
theory of non--commutative geometry briefly summarized in section
\ref{road}. Only recently, several results confirmed the existence
of a rich interplay between the algebraic and functional analytic
aspects of noncommutative geometry, especially through the work of
Connes and Dubois-Violette (\cf \cite{CDV1}, \cite{CDV2}, \cite{CDV3}) 
and of Polishchuk (\cf \cite{Poli}). Also, the work of Chakraborty and Pal
\cite{ChakraPal} and Connes \cite{Co-Qgr} and more recently of van
Suijlekom, Dabrowski, Landi, Sitarz, and Varilly \cite{SDLSV},
\cite{SDLSV2} showed that quantum groups fit very nicely within the
framework of noncommutative geometry described by spectral triples,
contrary to what was previously belived. Ultimately, successfully
importing tools from the theory of operator algebras into the realm
of algebraic geometry might well land within the framework of what
Manin refers to as a ``second quantization of algebraic geometry''.

\bigskip

\section{Quantum groups}\label{qgroups}

For a long time it was widely believed that quantum groups could not
fit into the setting of noncommutative manifolds defined in terms of
spectral geometry. On the contrary, recent work of Chakraborti and Pal
showed in \cite{ChakraPal} that the quantum group
$SU_q(2)$, for $0\leq q < 1$, admits a spectral triple with Dirac
operator that is equivariant with respect to its own (co)action.

\smallskip

The algebra $\cA$ of functions on the quantum group $SU_q(2)$ is generated
by two elements $\alpha$ and $\beta$ with the relations
\begin{equation}\label{Qgr-alg}
\begin{array}{c}
\alpha^*\alpha +\beta^*\beta =1, \ \ \ \ \alpha\alpha^*+ q^2  \beta\beta^*
=1, \\[2mm] \alpha\beta = q\beta\alpha, \ \ \  \alpha\beta^* =q
\beta^*\alpha, \ \ \  \beta^*\beta =\beta\beta^*. \end{array}
\end{equation}

\smallskip

By the representation theory of the quantum group $SU_q(2)$ (\cf
\cite{KlimSchmu}) there exists a Hilbert space $\cH$ with orthonormal
basis $e^{(n)}_{ij}$, $n\in \frac{1}{2} \N$, $i,j\in \{ -n,\ldots, n \}$, and
a unitary representation
\begin{equation}\label{Qgr-rep}
\begin{array}{ll}
\alpha\, e^{(n)}_{ij} = & a_+(n,i,j)\, e^{(n+1/2)}_{i-1/2, j-1/2}   +
a_-(n,i,j)\, e^{(n-1/2)}_{i-1/2, j-1/2}  \\[4mm]
\beta\, e^{(n)}_{ij} = & b_+(n,i,j)\, e^{(n+1/2)}_{i+1/2, j-1/2}   +
b_-(n,i,j)\, e^{(n-1/2)}_{i+1/2, j-1/2},
\end{array}
\end{equation}
with coefficients
$$ \begin{array}{ll}
a_+(n,i,j) = & q^{2n+i+j+1}\, Q(2n-2j+2,2n-2i+2,4n+2,4n+4) \\[2mm]
 a_-(n,i,j) = & Q(2n+2j,2n+2i,4n,4n+2) \\[2mm]
 b_+(n,i,j) = & -q^{n+j} \, Q(2n-2j+2,2n+2i+2,4n+2,4n+4) \\[2mm]
 b_-(n,i,j) = & q^{n+i} \, Q(2n+2j,2n-2i,4n,4n+2),
\end{array} $$
where we use the notation
$$ Q(n,m,k,r)=\frac{(1-q^{ n })^{1/2}(1-q^{ m })^{1/2}}{(1-q^{ k
})^{1/2}(1-q^{ r })^{1/2}}. $$

\smallskip

Consider then, as in \cite{Co-Qgr}, the operator
\begin{equation}\label{Qgr-Dirac}
D\, e^{(n)}_{ij} = \left\{ \begin{array}{rr} -2n & n\neq i \\
2n & n=i. \end{array} \right.
\end{equation}
More generally, one can consider operators of the form
$D\, e^{(n)}_{ij} = d(n,i)\, e^{(n)}_{ij}$, as in \cite{ChakraPal},
with $d(n,i)$ satisfying the conditions
$d(n+1/2,i+1/2)-d(n,i) =O(1)$ and $d(n+1/2,i-1/2)-d(n,i) = O(n+i+1)$.
Then one has the following result (Chakraborti--Pal
\cite{ChakraPal}):

\begin{thm}\label{Qgr-Sp3}
The data $(\cA,\cH,D)$ ad above define an $SU_q(2)$ equivariant odd
3--summable spectral triple.
\end{thm}

The equivariance condition means that there is an action on $\cH$ of
the enveloping algebra $\cU=U_q(\SL(2))$, which commutes with the
Dirac operator $D$. This is generated by operators
$$ \begin{array}{rl} k\, e^{(n)}_{ij} =&  q^j \, e^{(n)}_{ij} \\[3mm]
e\, e^{(n)}_{ij} =& q^{-n+1/2} (1-q^{2(n+j+1)})^{1/2}
(1-q^{2(n-j)})^{1/2} (1-q^2)^{-1} \, e^{(n)}_{ij+1}, \end{array} $$
satisfying the relations
$$ ke=qek, \,\,\,\,  kf=q^{-1} fk, \,\,\,\,
[e,f]=\frac{k^2-k^{-2}}{q-q^{-1}},
$$
with $f=e^*$, and with coproduct
$$ \Delta(k)=k\otimes k, \,\,\,\, \Delta(e)=k^{-1}\otimes e + e\otimes
k, \,\,\,\, \Delta(f) =k^{-1}\otimes f + f\otimes k. $$

\smallskip

It is interesting that, while the classical $SU(2)$ is of
(topological and metric) dimension three, the topological dimension of
the algebra $\cA$ of $SU_q(2)$ drops to one (\cf \cite{ChakraPal}),
but the metric dimension of the spectral triple remains equal to
three as in the classical case.

\smallskip

Chakraborti and Pal showed in \cite{ChakraPal} that the Chern
character of the spectral triple is nontrivial. Moreover, Connes in
\cite{Co-Qgr} gave an explicit formula for its local index cocycle,
where a delicate calculation provides the cochain whose coboundary is
the difference between the Chern character and the local version in
terms of remainders in the rational approximation to the logarithmic
derivative of the Dedekind eta function.

\smallskip

The local index formula is obtained by constructing a symbol map
$$ \rho : \cB \to C^\infty(S^*_q), $$
where the algebra $C^\infty(S^*_q)$ gives a noncommutative version of
the cosphere bundle, with a restriction map
$r:C^\infty(S^*_q) \to C^\infty(D_{q+}^2 \times D_{q-}^2)$ to the
algebra of two noncommutative disks. Here
$\cB$ is the algebra generated by the elements $\delta^k(a)$, $a\in
\cA$, with $\delta(a) =[|D|,a]$. On the cosphere bundle there is a
geodesic flow, induced by the group of automorphisms $a\mapsto e^{it
|D|} \, a \, e^{-it |D|}$. Then $\rho(b)^0$ denotes the component of
degree zero with respect to the grading induced by this flow.

\smallskip

The algebra $C^\infty(D^2_q)$ is an extension
$$ 0 \to \cS \to C^\infty(D^2_q) \stackrel{\sigma}{\to} C^\infty(S^1)
\to 0, $$
where the ideal $\cS$ is the algebra of rapidly decaying
matrices. There are linear
functionals $\tau_0$ and $\tau_1$ on $C^\infty(D^2_q)$,
$$ \tau_1(a) = \frac{1}{2\pi} \int_0^{2\pi} \sigma(a) d\theta, $$
$$ \tau_0(a) = \lim_{N\to \infty} \sum_{k=0}^N \langle a\, \epsilon_k,
\epsilon_k \rangle - \tau_1(a) N, $$
where $\tau_0$ is defined in terms of the representation of
$C^\infty(D^2_q)$ on the Hilbert space $\ell^2(\N)$ with o.n.~basis $\{
\epsilon_k \}$.

\smallskip

Recall that (\cf \cite{Co94}) a cycle $(\Omega,d,\int)$ is a triple with
where $(\Omega , d )$ is a  graded differential algebra,
and $\int : \Omega^n \to \C$ is a closed graded trace on $\Omega$.
A cycle over an algebra $\cA$ is given by a cycle $(\Omega , d , \int
)$ together with a homomorphism $\rho :
\cA \to \Omega^0$.

In the case of the algebra $\cA$ of $SU_q(2)$, a cycle $(\Omega,d,\int)$
is obtained in \cite{Co-Qgr} by considering $\Omega^1
= \cA \oplus \Omega^{(2)}(S^1)$, with $\Omega^{(2)}(S^1)$ the space of
weight two differential forms $f(\theta) d\theta^2$, with the
$\cA$--bimodule structure
$$ a\, (\xi,f)=(a\xi,\sigma(a) f) \ \ \ \  (\xi,f)\,a = (\xi a,
-i\sigma(\xi) \sigma(a)' + f \sigma(a)), $$
with differential
$$ da = \partial a+ \frac{1}{2} \sigma(a)'' d\theta^2, $$
with $\partial$ the derivation $\partial = \partial_\beta -
\partial_\alpha$, and
$$ \int \, (\xi,f) = \tau(\xi) + \frac{1}{2\pi i} \int f\, d\theta, $$
where $\tau(a) = \tau_0(r_-(a^{(0)}))$, with $a^{(0)}$ the component
of degree zero for $\partial$ and $r_-$ the restriction to
$C^\infty(D^2_{q-})$. This definition of the cycle corrects for the
fact that $\tau$ itself (as well as $\tau_0$) fails to be a trace.

\smallskip

The following result then holds (Connes \cite{Co-Qgr}):

\begin{thm}\label{loc-Qgr}
\begin{enumerate}
\item The spectral triple $(\cA,\cH,D)$ of Theorem \ref{Qgr-Sp3} has
dimension spectrum $\Sigma=\{ 1,2,3 \}$.
\item The residue formula for pseudodifferential operators $a\in
\cB$ in terms of their symbol is given by
$$ \begin{array}{ll} \cutint\,\, a \, |D|^{-3} = & (\tau_1\otimes \tau_1)
(r\rho(a)^0) \\[2mm]
 \cutint\,\, a \, |D|^{-2} = & (\tau_1\otimes\tau_0+\tau_0\otimes
\tau_1)(r\rho(a)^0) \\[2mm]
\cutint\,\, a \, |D|^{-1} = & (\tau_0\otimes\tau_0)(r\rho(a)^0)
\end{array} $$
\item The character $\chi(a_0,a_1)= \int a_0 \, da_1$ of the cycle
$(\Omega,d,\int)$ is equal to the cocycle
$$ \psi_1(a_0,a_1)= 2 \cutint a_0 \delta (a_1) P |D|^{-1} - \cutint a_0
\delta^2 (a_1) P |D|^{-1}, $$
with $P=(1+F)/2$. The local index formula is given by
$$ \varphi_{odd} = \psi_1 + (b+B) \varphi_{even}, $$
where $\varphi$ is the local index cocycle.
\item The character $\Tr(a_0 [F,a_1])$ differs from the local form
$\psi_1$ by the coboundary $b\psi_0$, with $\psi_0(a) = 2 \Tr (a P
|D|^{-s})_{s=0}$. This cochain is determined by the values
$\psi_0((\beta^*\beta)^n)$, which are of the form
$$ \psi_0((\beta^*\beta)^n) = q^{-2n}(q^2 R_n(q^2) - G(q^2)), $$
where $G$ is the logarithmic derivative of the Dedekind eta function
\begin{equation}\label{DedEta}
 \eta(q^2)=q^{1/12} \prod_{k=1}^\infty (1-q^{2k}),
\end{equation}
and the $R_n$ are rational functions with poles only at roots of
unity.
\end{enumerate}
\end{thm}

More recently, another important breakthrough in the relation between
quantum groups and the formalism of spectral triples was obtained by
Walter van Suijlekom, Ludwik Dabrowski, Giovanni
Landi, Andrzej Sitarz, Joseph C. Varilly, in \cite{SDLSV} and
\cite{SDLSV2}.

They construct a $3^+$ summable spectral triple $(\cA,\cH,D)$,
where $\cA$ is, as before, the algebra of coordinates of the
quantum group $SU_q(2)$. The geometry in this case is an isospectral
deformation of the classical case, in the sense that the Dirac
operator is the same as the usual Dirac operator for the round metric
on the ordinary 3-sphere $S^3$. Moreover, the spectral triple
$(\cA,\cH,D)$ is especially nice, in as it is equivariant with
respect to both left and right action of the Hopf algebra
$\cU_q(su_q(2))$.

The classical Dirac operator for the round metric on $S^3$ has
spectrum $\Sigma=\Sigma_+ \cup \Sigma_-$ with $\Sigma_+=\{ (2j+3/2):
j=0,1/2,1,3/2,\ldots \}$ with multiplicities $(2j+1)(2j+2)$ and
$\Sigma_-=\{ -(2j+1/2): j=1/2,1,3/2,\ldots\}$ with multiplicities
$2j(2j+1)$. The Hilbert space is obtained by taking $V\otimes \C^2$,
where $V$ is the left regular representation of $\cA$. It is very
important here to take $V\otimes \C^2$ instead of $\C^2\otimes V$.
Not only the latter violates the equivariance condition, but it was
shown by Ghoswami that it produces unbounded commutators $[D,a]$,
hence one does not obtain a spectral triple in that way.

The spectral triple contructed in \cite{SDLSV} and \cite{SDLSV2} has
a real structure $J$  and the Dirac operator satisfies a weak form
of the ``order one condition'' (\cf Section \ref{singconic} above).
The local index formula of \cite{Co-Qgr} (\cf Theorem \ref{loc-Qgr}
above) extends to the spectral triple of \cite{SDLSV}, as proved in
\cite{SDLSV2} and the structures of the cotangent space and the
geodesic flow are essentially the same.

\bigskip
\bigskip

\section{Spherical manifolds}\label{ncspheres}

The noncommutative spheres $S^3_\varphi\subset{\mathbb R}^4_\varphi$
are obtained as solutions of a very simple problem, namely the vanishing
of the first component of the Chern character of a unitary $U\in
M_2(\cA)$ where $\cA$ is the algebra of functions on the sphere and
the Chern character is taken in the cyclic homology (b,B) bicomplex.
The origin of this problem is to quantize the volume form of a three
manifold (\cf \cite{CDV1}). The solutions are parameterized by three
angles $\varphi_k$, $k\in \{1,2,3\}$ and the corresponding algebras
are obtained by imposing the ``unit sphere relation"
\begin{equation}
\sum \,x_\mu^2=\,1
\end{equation}
to the four generators $x_0, x_1, x_2, x_3$ of the quadratic algebra
 $C_{\mathrm{alg}}({\mathbb R}^4_\varphi)$
with the six relations

\smallskip
\begin{equation}
\label{pres1} \sin (\varphi_k) \, [x_0 , x_k]_+ = i\; \cos
(\varphi_{\ell} - \varphi_m) \, [x_{\ell} , x_m]
\end{equation}
\begin{equation}
\label{pres2} \cos (\varphi_k) \, [x_0 , x_k] = i \;\sin
(\varphi_{\ell} - \varphi_m) \, [x_{\ell} , x_m]_+ \, ,
\end{equation}

\medskip

where  $ [a , b]_+=\,a\,b+\,b\,a$ is the anticommutator and by
convention the indices $k,l,m\in \{1,2,3\}$ always appear in cyclic
order.

The analysis of these algebras is a special case of the general
theory of central quadratic forms for quadratic algebras developed
in \cite{CDV2}, \cite{CDV3} and which we briefly recall below.

Let $\cala=A(V,R)=T(V)/(R)$  be a quadratic algebra where $V$ is the
linear span of the generators and $(R)\subset T(V)$ the ideal
generated by the relations. The geometric data $\{E\,,\,
\sigma\,,\,\call\}$ is given by an algebraic variety $E$, a
correspondence $\sigma$ on $E$ and a line bundle $\call$ over $E$.
These data are defined so as to yield an homomorphism $h$ from
$\cala$ to a crossed product algebra constructed from sections of
powers of the line bundle $\call$ on the graphs of the iterations of
the correspondence $\sigma$. This crossed product only involves the
positive powers of the correspondence $\sigma$ and thus remains
``triangular" and far remote from the ``semi-simple" set-up of
$C^*$-algebras.

This morphism  $h$ can be considerably refined using the notion of
positive central quadratic form.

\begin{defn} \label{cent}
Let $Q \in S^2(V)$ be a symmetric bilinear form on $V^\ast$ and $C$
a component of $E \times E$. We   say that $Q$ is
\underline{central} on $C$ iff for all ($Z,\,Z'$) in $C$ and
$\omega\in R$ one has,
\begin{equation} \label{defcentral}
\omega(Z,Z')\, Q(\sigma(Z'),\sigma^{-1}(Z))+Q(Z,Z')\,
\omega(\sigma(Z'),\sigma^{-1}(Z)) =0
\end{equation}
\end{defn}

This makes it possible to construct purely algebraically a crossed product
algebra and an homomorphism from $\cala=A(V,R)$ to this crossed
product \cite{CDV2}, \cite{CDV3}. 
The relation with $C^*$-algebras arises from
positive central quadratic forms which make sense on involutive
quadratic algebras.

Let $\cala=A(V,R)$ be an {\sl involutive} quadratic algebra \ie an
algebra over $\C$ which is a
 $\ast$-algebra with  involution $x\mapsto x^\ast$ preserving the subspace
  $V$ of the generators.
The real structure of $V$ is given by the antilinear involution
$v\mapsto j( v)$ restriction of $x\mapsto x^\ast$. As
$(xy)^\ast=y^\ast x^\ast$ for $x,y\in \cala$, the space  $R$ of
relations fulfills
\begin{equation}
(j \otimes j)( R)=t(R) \label{eq5.3}
\end{equation}
 in $V\otimes V$ where $t:V\otimes V\rightarrow V\otimes V$ is the
transposition $v\otimes w \mapsto t(v\otimes w)=w\otimes v$. This
implies that the characteristic variety is stable under the
involution $j$ and one has
$$
\sigma( j(Z)) =\, j( \sigma^{-1}(Z))
$$

\smallskip

Let then  $C$ be an invariant component of $E \times E$, we say that
$C$ is $j$-{\em real} when it is globally invariant under the
involution
\begin{equation}\label{tildej}
\tilde j(Z,\,Z'):=( j( Z'),\, j( Z))
\end{equation}
Let then $Q$ be a central qudratic form on $C$, we say that $Q$ is
positive on $C$ iff

$$
Q(Z,j(Z))> 0 \qqq Z\in K\,.
$$

\smallskip One can then endow the line bundle  $\call$ dual of the tautological bundle
on $P(V^*)$ with the hermitian metric given by
\begin{equation}\label{herm}
\langle f\,L,\,g\, L'\rangle_Q(Z) =\,f(Z)\,\overline{g(Z)} \,\frac{
L(Z)\,\overline{ L'(Z)}}{Q(Z,\,j(Z))} \qquad L, L' \in V,\quad Z \in
K\, \,.
\end{equation}
($\forall f,g \in C(K)$)

\medskip One then defines a generalized crossed product $C^*$-algebra
 $C(K) \times_{\sigma,\,\call}
\mathbb {Z} $ following M. Pimsner \cite{pims:1997}. Given a compact
space $K$, an homeomorphism $\sigma$ of $K$ and a hermitian line
bundle $\call$ on $K$ we define the  $C^\ast$-algebra $C(K)
\times_{\sigma,\,\call} \mathbb {Z} $ as the twisted cross-product
of $C(K)$ by the Hilbert $C^*$-bimodule associated to $\call$ and
$\sigma$ (\cite{aba-eil-exel:1998}, \cite{pims:1997}).

We let for each $n \geq 0$, $\call^{\sigma^n}$ be the hermitian line
bundle pullback of $\call$ by $\sigma^n$ and (cf.
\cite{art-tat-vdb:1990}, \cite{smi-sta:1992})
\begin{equation}
\call_n := \call \otimes \call^{\sigma} \otimes \cdots \otimes
\call^{\sigma^{n-1}} \label{gene2}
\end{equation}
We first define a $\ast$-algebra as the linear span of the monomials
\begin{equation}
\xi \, W^n\, , \quad W^{\ast n} \, \eta^\ast \,,\quad \xi\,,\eta \in
C(K,\call_n) \label{gene}
\end{equation}
with product given as in (\cite{art-tat-vdb:1990},
\cite{smi-sta:1992}) for $(\xi_1 \, W^{n_1})\,(\xi_2 \, W^{n_2})$ so
that
\begin{equation}
(\xi_1 \, W^{n_1})\,(\xi_2 \, W^{n_2}):= (\xi_1 \otimes
(\xi_2\circ{\sigma^{n_1}}) )\, W^{n_1+n_2} \label{gene3}
\end{equation}
We use the hermitian structure of $\call_n $ to give meaning to the
products $\eta^\ast \,\xi$ and $\xi \;\eta^\ast$ for $\xi\,,\eta \in
C(K,\call_n)$. The product then extends uniquely to an associative
product of $\ast$-algebra fulfilling the following additional rules
\begin{equation}
(W^{\ast k} \, \eta^\ast)\,( \xi \, W^k):= \, (\eta^\ast\, \xi)\circ
\sigma^{-k}\,,\qquad ( \xi \, W^k)\,(W^{\ast k} \, \eta^\ast)\,:= \,
\xi \;\eta^\ast \label{gene1}
\end{equation}

 The $C^\ast$-norm of $C(K) \times_{\sigma,\,\call} \mathbb {Z} $
 is defined as for ordinary cross-products and due to the amenability of the group $\mathbb {Z} $
there is no distinction between the reduced and maximal norms. The
latter is obtained as the supremum of the norms in involutive
representations in Hilbert space. The natural positive conditional
expectation on the subalgebra $C(K)$ shows that the $C^\ast$-norm
restricts to the usual sup norm on $C(K)$.

\medskip

\begin{thm}\label{C*}
 Let $K \subset E$ be a compact
$\sigma$-invariant subset and $Q$ be central and strictly positive
on $\{(Z,\,\bar Z);\, Z\in K\}$. Let $\call$ be the restriction to
$K$ of the dual of the tautological line bundle on $P(V^\ast)$
endowed with the hermitian metric $\langle\;,\; \rangle_Q$.

 (i) The equality $\sqrt{2}\,\theta(Y):= Y\, W + W^\ast\,\bar Y^\ast$
yields a $\ast$-homomorphism $$\theta:\cala=A(V,R) \to C(K)
\times_{\sigma,\,\call} \mathbb {Z} $$

 (ii) For any $Y \in V$
the $C^\ast$-norm of $\theta(Y)$ fulfills
$$\sup_K \|Y\|\leq \sqrt{2}\| \,\theta(Y)\|
\leq 2\sup_K \|Y\| $$

 (iii) If $\sigma^4 \neq \bbbone$, then $\theta(Q)= 1$ where
$Q$ is viewed as an element of $T(V)/(R)$.
\end{thm}

\bigskip

In the above case of the sphere $S^3_\varphi$ one lets $Q$ be the
quadratic form
\begin{equation}
 Q(x,\,x'):=\sum x_\mu\,x'_\mu
\label{quad}
\end{equation}

In the generic case one has :

\begin{prop}  1) The characteristic variety is the union of 4 points with an elliptic curve $F_\varphi$.

2) The quadratic form  $Q$ is central and positive on $F_\varphi
\times F_\varphi$.
\end{prop}

\smallskip
In suitable coordinates the equations defining the elliptic curve
$F_\varphi$ are
\begin{equation}
\frac{Z_0^2-Z_1^2}{
s_1}=\frac{Z_0^2-Z_2^2}{s_2}=\frac{Z_0^2-Z_3^2}{s_3} \label{charZ}
\end{equation}
where $s_k := 1 + t_\ell\, t_m \, , \: t_k:= {\rm tan}\,\varphi_k$.

The positivity of $Q$ is automatic since in the coordinates  $x$ the
involution $j_\varphi$
 of the $\ast$-algebra $C_{\mathrm{alg}}({\mathbb R}^4_\varphi)$
is just $j_\varphi(Z)= \bar Z$, so that $Q(X,\,j_\varphi(X))>0$ for
$X \neq 0$.

\begin{cor} \label{II}
Let $K \subset F_\varphi $  be a compact $\sigma$-invariant subset.
The homomorphism $\theta$ of Theorem \ref{C*} is a unital
 $\ast$-homomorphism from $ C_{\mathrm{alg}}(S^3_\varphi)$
to the cross-product $ C^{\infty}(K) \times_{\sigma,\,\call} \mathbb
{Z} $.
\end{cor}

 It follows  that one obtains a
non-trivial $C^\ast$-algebra $C^\ast(S^3_\varphi)$ as the completion
of  $ C_{\mathrm{alg}}(S^3_\varphi)$ for the semi-norm,
\begin{equation}
 \| P \|:= \sup \| \,\pi(P) \|
\label{norm}
\end{equation}
where $\pi$ varies through all unitary representations of $
C_{\mathrm{alg}}(S^3_\varphi)$. It was clear from the start that
\eqref{norm} defines a finite $C^\ast$-semi-norm on $
C_{\mathrm{alg}}(S^3_\varphi)$ since the equation of the sphere
$\sum x_\mu^2=1$ together with the self-adjointness
$x_\mu=\,x_\mu^*$ show that in any unitary representation one has
$$
\|\, \pi(x_\mu) \|\leq 1\qqq \mu\,.
$$
What the above corollary gives is a lower bound for the
$C^\ast$-norm such as that given by statement (ii) of Theorem
\ref{C*} on the linear subspace $V$ of generators.

\smallskip The correspondence
$\sigma$ on $F_\varphi $,
 is for generic  $\varphi$ a translation of module $\eta$ of the elliptic curve
 $F_\varphi$ and one distinguishes two cases :  the {\em even} case when it preserves
 the two real components of the curve $F_\varphi \cap
P_3(\mathbb {R})$ and the odd case when it permutes them.

\begin{prop} \label{siginvar} Let $\varphi  $ be generic and even.

(i) The cross-product $ C(F_\varphi) \times_{\sigma,\,\call} \mathbb
{Z} $ is isomorphic to the mapping torus of the automorphism $\beta
$ of the noncommutative torus ${\mathbb T}_{\eta}^2 = C_\varphi
\times_\sigma \mathbb {Z} $ acting on the generators by the  matrix
$\left[
\begin{array}{cc}
1& 4\\
0& 1
\end{array}
\right] $.

(ii) The crossed product $F_\varphi \times_{\sigma,\,\call} \mathbb
{Z} $ is a noncommutative $3$-manifold with an elliptic
 action of the three dimensional
Heisenberg Lie algebra $\frach_3$ and an invariant trace $\tau$.

\end{prop}

\medskip

It follows that one is exactly in the framework developed in
\cite{C3}. We refer to \cite{Rieffel89} and \cite{aba-exel:1997}
where these noncommutative manifolds were analyzed in terms of
crossed products by Hilbert $C^*$-bimodules.

 Integration on the translation invariant volume form $dv$
of $F_\varphi$ gives the $\frach_3$-invariant trace $\tau$,
\begin{eqnarray}
 \label{trace}
\tau(f)& = & \int f dv\,,\quad \forall f \in C^{\infty}(F_\varphi)\nonumber\\
\tau(\xi \, W^k)& = &\tau(W^{\ast k} \, \eta^\ast)\,=\,0\,,\quad
\forall k\neq 0
 \end{eqnarray}
It follows in particular that the results of \cite{C3} apply to
obtain the calculus. In particular the following gives the
``fundamental class" as a $3$-cyclic cocycle,
\begin{equation}\label{3trace}
\tau_3(a_0,\,a_1,\,a_2 ,\,a_3)=\,\sum
\epsilon_{ijk}\,\tau(a_0\,\delta_i(a_1)\,\delta_j(a_2)
\,\delta_k(a_3))
\end{equation}
where the $\delta_j$ are the generators of the action of $\frach_3$.

\medskip The relation between the noncommutative spheres
$S^3_\varphi$ and the noncommutative nilmanifolds $F_\varphi
\times_{\sigma,\,\call} \mathbb {Z} $ is analyzed in 
\cite{CDV2}, \cite{CDV3} thanks to the computation of the 
Jacobian of the homomorphism $\theta$.

\bigskip

\section{$\Q$-lattices}\label{Qlatt}

A class of examples of noncommutative spaces of relevance to
number theory is given by the moduli spaces of $\Q$-lattices up to
commensurability. These fall within the general framework of
noncommutative spaces obtained as quotients of equivalence relations
discussed in Section \ref{quotients}.

A $\Q$-lattice in $\R^n$ consists of a pair $( \Lambda , \phi) $
of a lattice $\Lambda\subset \R^n$ (a cocompact free abelian
subgroup of $\R^n$ of rank $n$) together with a system of labels
of its torsion points given by a homomorphism of abelian groups
\begin{equation}\label{phimap}
 \phi :  \Q^n/\Z^n \longrightarrow \Q\Lambda / \Lambda.
\end{equation}

Two $\Q$-lattices are commensurable,
$$ (\Lambda_1, \phi_1) \sim (\Lambda_2, \phi_2), $$
iff $\Q\Lambda_1=\Q\Lambda_2$ and
$$ \phi_1 = \phi_2 \mod \Lambda_1 + \Lambda_2 $$

In general, the map $\phi$ of \eqref{phimap} is just a group
homomorphism. A $\Q$-lattice is said to be {\em invertible}
is $\phi$ is an isomorphism. Two invertible $\Q$-lattices are
commensurable if and only if they are equal.

\begin{center}
\begin{figure}
\includegraphics[scale=0.95]{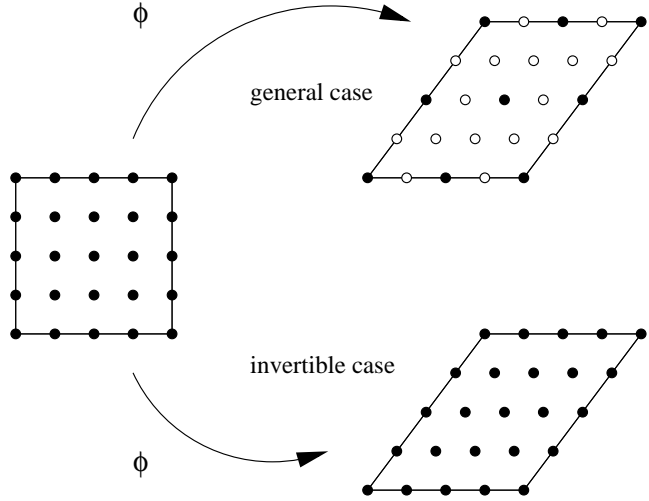}
\caption{$\Q$-lattices: generic and invertible case.
\label{FigQlat}}
\end{figure}
\end{center}

The space $\cL_n$ of commensurabilty classes of
$\Q$-lattices in $\R^n$ has the typical property of noncommutative
spaces: it has the cardinality of the continuum but one cannot
construct a countable collection of measurable functions that
separate points of $\cL_n$.
Thus, one can use noncommutative geometry to describe the quotient
space $\cL_n$ through a noncommutative $C^*$-algebra $C^*(\cL_n)$.

We consider especially the case of $n=1$ and $n=2$. One is also
interested in the $C^*$-algebras describing $\Q$-lattices up to
scaling, $\cA_1=C^*(\cL_1/\R_+^*)$ and $\cA_2=C^*(\cL_2/\C^*)$.

In the 1-dimensional case, a $\Q$-lattice can always be written in the
form
\begin{equation}\label{1dimQlat}
 ( \Lambda , \phi) \, = (\lambda\, \Z,\lambda\,\rho)
\end{equation}
for some $\lambda>0$ and some
\begin{equation}\label{rho}
\rho \in \Hom(\Q/\Z,\Q/\Z)=\varprojlim \Z/n\Z = \hat\Z.
\end{equation}
By considering lattices up to scaling, we eliminate the factor
$\lambda>0$ so that 1-dimensional $\Q$-lattices up to scale are
completely specified by the choice of the element $\rho\in
\hat\Z$. Thus, the algebra of coordinates of the space of
1-dimensional $\Q$-lattices up to scale is the commutative
$C^*$-algebra
\begin{equation}\label{ChatZ}
C(\hat\Z)\simeq C^*(\Q/\Z),
\end{equation}
where we use Pontrjagin duality to get the identification in
\eqref{ChatZ}.

The equivalence relation of commensurability is implemented by the
action of the semigroup $\N^\times$ on $\Q$-lattices. The
corresponding action on the algebra \eqref{ChatZ} is by
\begin{equation}\label{Nact}
 \alpha_n(f) (\rho)=\left\{ \begin{array}{lr} f(n^{-1} \rho) & \rho \in
 n\hat\Z \\ 0 & \text{otherwise.} \end{array}\right.
\end{equation}

Thus, the quotient of the space of 1-dimensional $\Q$-lattices up
to scale by the commensurability relation and its algebra of
coordinates is given by the semigroup crossed product
\begin{equation}\label{semicross}
 C^*(\Q/\Z)\rtimes \N^\times .
\end{equation}
This is the Bost--Connes $C^*$-algebra introduced in \cite{BC}.

It has a natural time evolution given by the covolume of a pair of
commensurable $\Q$-lattices. It has symmetries (compatible with the
time evolution) given by the group $\hat\Z^*=\GL_1(\A_f)/\Q^*$ and
the KMS (Kubo--Martin--Schwinger) equilibrium states of the system
have interesting arithmetic properties. Namely, the partition function
of the system is the Riemann zeta function. There is a unique KMS
state for sufficiently high temperature, while at low temperature the
system undergoes a phase transition with spontaneous symmetry
breaking. The pure phases (estremal KMS states) at low temperature
are parameterized by elements in $\hat\Z^*$. They have an explicit
expression in terms of polylogarithms at roots of unity. At zero
temperature the extremal KMS states, evaluated on the elements of a
rational subalgebra affect values that are algebraic numbers.
The action on these values of the Galois group $\Gal(\bar\Q/\Q)$
factors through its abelianization and
is obtained (via the class field theory isomorphism $\hat\Z^* \cong
\Gal(\Q^{ab}/\Q)$) as the action of symmetries on the algebra (\cf
\cite{BC}, \cite{CoMa}, \cite{CoMajgp1} for details).

\medskip

In the 2-dimensional case, a $\Q$-lattice can be written in the form
$$ (\Lambda,\phi)=(\lambda (\Z+\Z\tau),\lambda\rho), $$
for some $\lambda\in \C^*$, some $\tau\in \H$, and some $\rho\in
M_2(\hat\Z)=\Hom(\Q^2/\Z^2,\Q^2/\Z^2)$. Thus, the space of
2-dimensional $\Q$-lattices up to the scale factor $\lambda\in
\C^*$ and up to isomorphisms, is given by
\begin{equation}\label{2dQlat}
 M_2(\hat\Z)\times \H  \mod \Gamma=\SL_2(\Z).
\end{equation}
The commensurability relation giving the
space $\cL_2/\C^*$ is implemented by the partially
defined action of $\GL_2^+(\Q)$.

One considers in this case the quotient of the space
\begin{equation}\label{tildeU}
\tilde \cU:= \{ (g,\rho,\alpha)\in \GL_2^+(\Q)\times
M_2(\hat\Z)\times \GL_2^+(\R)\, :\,\, g\rho\in M_2(\hat\Z) \}
\end{equation}
by the action of $\Gamma \times \Gamma$ given by
\begin{equation}\label{Gamma2action}
(\gamma_1,\gamma_2) \, (g,\rho,\alpha)=(\gamma_1 g \gamma_2^{-1},
\gamma_2 \rho, \gamma_2 \alpha).
\end{equation}

The groupoid $\cR_2$ of the equivalence relation of
commensurability on 2-dimensional $\Q$-lattices (not considered up
to scaling for the moment) is a locally compact groupoid, which
can be parameterized by the quotient of \eqref{tildeU} by
$\Gamma\times \Gamma$ via the map $r: \tilde \cU \to \cR_2$,
\begin{equation}\label{mapgr}
r(g,\rho,\alpha)=\left( (\alpha^{-1}g^{-1} \Lambda_0, \alpha^{-1}
\rho), (\alpha^{-1}\Lambda_0, \alpha^{-1} \rho)\right).
\end{equation}

We then consider the quotient by scaling.
The quotient $\GL_2^+(\R)/\C^*$ can be identified with the
hyperbolic plane $\H$ in the usual way.
If $(\Lambda_k,\phi_k)$ $k=1,2$ are a pair of commensurable
2-dimensional $\Q$-lattices, then for any $\lambda\in \C^*$, the
$\Q$-lattices $(\lambda \Lambda_k, \lambda \phi_k)$ are also
commensurable, with
$$ r(g,\rho,\alpha\lambda^{-1})=\lambda r(g,\rho,\alpha). $$
However, the action of $\C^*$ on $\Q$-lattices is not free due to
the presence of lattices (such as $\Lambda_0$ above) with nontrivial
automorphisms. Thus, the quotient $Z=\cR_2/\C^*$ is no longer a
groupoid. Still, one can define a convolution algebra for $Z$ by
restricting the convolution product of $\cR_2$ to homogeneous
functions of weight zero, where a function $f$ has weight $k$ if
it satisfies
$$ f(g,\rho,\alpha\lambda)=\lambda^k f(g,\rho,\alpha), \ \ \
\forall \lambda \in \C^*. $$
The space $Z$ is the quotient of the space
\begin{equation}\label{Uspace}
 \cU:=\{ (g,\rho,z) \in \GL_2^+(\Q)\times M_2(\hat\Z)\times \H |
g\rho\in M_2(\hat\Z) \}
\end{equation}
by the action of $\Gamma\times \Gamma$. Here the space
$M_2(\hat\Z)\times \H$ has a partially defined action of
$\GL_2^+(\Q)$ given by
$$ g (\rho,z)=(g\rho, g(z)), $$
where $g(z)$ denotes action as fractional linear transformation.

Thus, the algebra of coordinates $\cA_2$ for the noncommutative space
of commensurability classes of 2-dimensional $\Q$-lattices up
to scaling is given by the following convolution algebra.

Consider the space $C_c(Z)$ of continuous compactly supported
functions on $Z$. These can be seen, equivalently, as functions on
$\cU$ as in \eqref{Uspace} invariant under the
$\Gamma\times\Gamma$ action $(g,\rho,
z)\mapsto(\gamma_1g\gamma_2^{-1},\gamma_2 z)$.
One endows $C_c(Z)$ with the convolution product
\begin{equation}\label{Heckeprod2}
(f_1*f_2)(g,\rho,z)= \displaystyle{\sum_{s\in \Gamma\backslash
\GL_2^+(\Q): s\rho\in M_2(\hat\Z)}} f_1(gs^{-1},s\rho,s(z))
f_2(s,\rho,z)
\end{equation}
and the involution
$f^*(g,\rho,z)=\overline{f(g^{-1},g\rho,g(z))}$.

Again there is a time evolution on this algebra, which is given by the
covolume,
\begin{equation}\label{evolution2}
\sigma_t(f) (g,\rho,z) = \det(g)^{it} \, f(g,\rho,z).
\end{equation}
The partition function for this $\GL_2$ system is given by
\begin{equation}\label{partition2}
Z(\beta)= \sum_{m\in \Gamma\backslash M_2^+(\Z)}\det(m)^{-\beta}
=\sum_{k=1}^\infty \sigma(k)\,
k^{-\beta}=\zeta(\beta)\zeta(\beta-1),
\end{equation}
where $\sigma(k)=\sum_{d|k} d$. The form of the partition
function suggests the possibility that two distinct phase
transitions might happen at $\beta=1$ and $\beta=2$.

The structure of KMS states for this system is analysed in
\cite{CoMa}. The main result is the following.

\begin{thm}\label{GL2KMS}
The KMS$_\beta$ states of the $\GL_2$-system have the following
properties:
\begin{enumerate}
\item In the range $\beta\leq 1$ there are no KMS states.
\item In the range $\beta>2$ the set of extremal KMS states is
given by the classical Shimura variety
\begin{equation}\label{Ekmsbeta2}
\cE_\beta \cong \GL_2(\Q)\backslash \GL_2(\A) /\C^*.
\end{equation}
\end{enumerate}
\end{thm}

The symmetries are more complicated than in the Bost--Connes case. In
fact, in addition to symmetries given by automorphisms that commute
with the time evolution, there are also symmetries by {\em
endomorphisms} that play an important role. The resulting symmetry
group is the quotient $\GL_2(\A_f)/\Q^*$. An important result of
Shimura \cite{Shi} shows that this group is in fact the Galois group
of the field $F$ of modular functions. The group $\GL_2(\A_f)$
decomposes as a product
\begin{equation}\label{SymmGL2}
\GL_2(\A_f)=\GL_2^+(\Q) \GL_2(\hat\Z),
\end{equation}
where $\GL_2(\hat\Z)$ acts by automorphisms related to the
deck transformations of the tower of the modular curves, while
$\GL_2^+(\Q)$ acts by endomorphisms that move across levels in
the modular tower.

The modular field $F$ is the field of modular functions over
$\Q^{ab}$, namely the union of the fields $F_N$ of modular
functions of level $N$ rational over the cyclotomic field
$\Q(\zeta_n)$, that is, such that the $q$-expansion in powers of
$q^{1/N}=\exp(2\pi i \tau/N)$ has all coefficients in $\Q(e^{2\pi
i/N})$.

The action of the Galois group $\hat\Z^* \simeq \Gal(\Q^{ab}/\Q)$
on the coefficients determines a homomorphism
\begin{equation}\label{cyclhom}
{\rm cycl}: \hat\Z^* \to \Aut(F).
\end{equation}

If $\tau\in \H$ is a generic point, then the evaluation map
$f\mapsto f(\tau)$ determines an embedding $F\hookrightarrow \C$. We
denote by $F_\tau$ the image in $\C$. This yields an
identification
\begin{equation}\label{GalFtau}
 \theta_\tau: \Gal(F_\tau/\Q)
\stackrel{\simeq}{\to} \Q^* \backslash \GL_2(\A_f).
\end{equation}

There is an arithmetic algebra $\cA_{2,\Q}$ (defined over $\Q$) of
unbounded multipliers of the $C^*$-algebra $\cA_2$, obtained by
considering continuous functions on $Z$ (\cf
\eqref{Uspace}), with finite support in the variable $g\in
\Gamma\backslash\GL_2^+(\Q)$ and with the following properties.
Let $p_N: M_2(\hat\Z)\to M_2(\Z/N\Z)$ be the canonical
projection. With the notation $f_{(g,\rho)}(z) = f(g,\rho,z)$,
we say that $f_{(g,\rho)}\in C(\H)$ is of level $N$ if
$$ f_{(g,\rho)} = f_{(g,p_N(\rho))} \ \ \ \ \forall (g,\rho). $$
We require that elements of $\cA_{2,\Q}$ have the $f_{(g,\rho)}$ of
finite level with $ f_{(g,m)} \in F$ for all $(g,m)$. We also
require that the action \eqref{cyclhom} on the coefficients
of the q-expansion of the $f_{(g,m)}$ satisfies
$$ f_{(g,\alpha(u)m)} = {\rm cycl}(u) \, f_{(g,m)}, $$
for all $g\in \GL_2^+(\Q)$ diagonal and all $u\in \hat\Z^*$, with
$$ \alpha(u)=\begin{pmatrix} u& 0 \\ 0 & 1 \end{pmatrix}, $$
to avoid some ``trivial'' elements that would spoil the Galois
action on values of states (\cf \cite{CoMa}, \cite{CoMajgp1}).
The action of symmetries extends to $\cA_{2,\Q}$.
We have then the following result (\cite{CoMa}):

\begin{thm}\label{GalGL2infty}
Consider a state $\varphi=\varphi_{\infty,L}\in \cE_\infty$, for a
generic invertible $\Q$-lattice $L=(\rho,\tau)$. Then the
values of the state on elements of the arithmetic
subalgebra generate the image in $\C$ of the modular field,
\begin{equation}\label{values}
\varphi(\cA_{2,\Q})\subset F_\tau,
\end{equation}
and the isomorphism
\begin{equation}\label{thetaphi1}
\theta_\varphi : \Gal(F_\tau/\Q)
\stackrel{\simeq}{\longrightarrow} \Q^* \backslash \GL_2(\A_f),
\end{equation}
given by
\begin{equation}\label{thetaphi2}
 \theta_\varphi(\gamma)=\rho^{-1}\, \theta_\tau(\gamma) \, \rho,
\end{equation}
for $\theta_\tau$ as in \eqref{GalFtau}, intertwines the Galois
action on the values of the state with the action of symmetries,
\begin{equation}\label{intertwineGL2}
\gamma\, \varphi(f) = \varphi( \theta_\varphi(\gamma) f), \ \ \ \
\forall f\in \cA_{2,\Q}, \ \ \forall\gamma\in \Gal(F_\tau/\Q).
\end{equation}
\end{thm}

\medskip

A notion analogous to that of $\Q$-lattices can be given for other
number fields $\K$. This notion was used in \cite{CMR} to construct a
quantum statistical mechanical system for $\K$ an imaginary quadratic
field. This system shares properties with both the Bost--Connes system
of \cite{BC} and the $\GL_2$ system (2-dimensional $\Q$-lattices) of
\cite{CoMa}.

We assume that $\K=\Q(\sqrt{-d})$, $d$ a positive integer. Let
$\tau\in \H$ be such that
$\K=\Q(\tau)$ and $\cO = \Z + \Z \tau$ is the ring of integers of $\K$.

A 1-dimensional $\K$-lattice
$(\Lambda,\phi)$ is a finitely generated $\cO$-submodule $\Lambda\subset \C$,
such that $\Lambda\otimes_\cO \K \cong \K$, together with a morphism of
$\cO$-modules
\begin{equation}\label{Kphi}
 \phi : \K/\cO \to \K\Lambda/\Lambda.
\end{equation}
A 1-dimensional $\K$-lattice is {\em invertible} if $\phi$ is an
isomorphism of $\cO$-modules.
A 1-dimensional $\K$-lattice is, in particular, a 2-dimensional
$\Q$-lattice.

We consider the
notion of commensurability as in the case of $\Q$-lattices.
Two 1-dimensional $\K$-lattices
$(\Lambda_1,\phi_1)$ and
$(\Lambda_2,\phi_2)$ are commensurable if $\K\Lambda_1=\K\Lambda_2$
and $\phi_1=\phi_2$ modulo $\Lambda_1+\Lambda_2$.
In particular, two 1-dimensional $\K$-lattices are commensurable iff
the underlying $\Q$-lattices are commensurable.

The algebra of the corresponding noncommutative space is a restriction
of the algebra of the $\GL_2$-system to the subgroupoid of the
equivalence of commensurability restricted to $\K$-lattices.
The time evolution is also a restriction from the $\GL_2$-system.

The resulting system has partition function the Dedekind zeta function
$\zeta_\K(\beta)$ of the number field $\K$. Above the critical temperature
$T=1$ there is a unique KMS state, while at lower temperatures the
extremal KMS states are parameterized by elements of $\A_\K^*/\K^*$,
where $\A_K=\A_{K,f}\times \C$ are the adeles of $\K$, with
$\A_{K,f}=\A_f\otimes \K$. The KMS states at zero temperature,
evaluated on the restriction to $\K$-lattices of the arithmetic
algebra of the $\GL_2$-system, have an action of the Galois group
$\Gal(\K^{ab}/\K)$ realized (via the class field theory isomorphism)
through the action of symmetries (automorphisms and endomorphisms) of
the system (\cf \cite{CMR}).

\bigskip

\section{Modular Hecke algebras}\label{hecke}

Connes and Moscovici \cite{Co-Mosc1} defined modular Hecke algebras
$\cA(\Gamma)$ of level $\Gamma$, a congruence subgroup of
$\PSL_2(\Z)$. These extend both the ring of classical Hecke operators
and the algebra of modular forms.

Modular Hecke algebras  encode two \textit{a priori} unrelated
structures on modular forms, namely the algebra structure given by
the pointwise product on one hand, and the action of the Hecke
operators on the other. To any congruence subgroup $\Gamma$ of $\,
\SL_2(\Z)$ corresponds a crossed product algebra
$\cA(\Gamma)$, the \textit{modular Hecke algebra} of level $\Gamma$,
which is a direct extension of both the ring of classical Hecke
operators and of the algebra $\cM(\Gamma)$ of $\Gamma$-modular
forms.

These algebras can be obtained by considering the action of
$\GL^+_2(\Q)$ on the algebra of modular forms on the full (adelic)
modular tower, which yields the ``holomorphic part" of the ``ring of
functions" of the noncommutative space of commensurability classes of
2-dimensional $\Q$-lattices, introduced in Section \ref{Qlatt}.

With $\cM$ denoting the algebra of modular forms of arbitrary
level, the elements of $\cA(\Gamma)$ are maps with finite support
$$
 F: \Gamma\backslash \GL^+ (2, \mathbb{Q}) \to \cM \, ,
\qquad \Gamma  \alpha \mapsto F_{\alpha} \in \cM  \, ,
$$
satisfying the covariance condition
\begin{equation}
F_{\alpha \gamma} \, = \, F_{\alpha} \vert \gamma \qqq\alpha \in
\GL^+ (2, \mathbb{Q}) \, , \gamma \in \Gamma \, \nonumber
\end{equation}
and their product is given by convolution.

More in detail,
let $G=\PGL^+_2(\Q)$ and $\Gamma \subset \PSL_2(\Z)$ a finite
index subgroup. The quotient map $\Gamma \backslash G \to \Gamma
\backslash G / \Gamma$ is finite to one, and $\Gamma$ acts on $\C
[\Gamma \backslash G]$. Let $\cH_k$ be the space of holomorphic
functions $f: \H \to \C$ with polynomial growth, and with the
action $|_k$, for $k\in 2\Z$, of $\PGL^+_2(\R)$ of the form
$$ \left(f|_k \left(\begin{array}{cc} a& b\\ c& d \end{array}\right)\right)
(z)= \frac{(ad-bc)^{k/2}}{(cz+d)^k} \, f\left( \frac{az+b}{cz+d}
\right).
$$
This determines induced actions of $G$ and $\Gamma$ on $\cH_k$.
The space of modular forms is obtained as $\cM_k(\Gamma)=
\cH_k^\Gamma$, the invariants of this action.

One can then define
\begin{equation}\label{AkGamma}
 \cA_k (\Gamma) := \left( \C [ \Gamma \backslash G ] \otimes_\C
\cH_k \right)^\Gamma,
\end{equation}
with respect to the right action of $\Gamma$,
$$ \gamma : \sum_i (\Gamma g_i) \otimes f_i \mapsto \sum_i \Gamma g_i
\gamma \otimes (f_i |_k \gamma).
$$
One considers the graded vector space $\cA_*(\Gamma) = \oplus_k
\cA_k(\Gamma)$. The elements of $\cA_k(\Gamma)$ can be thought of
as finitely supported $\Gamma$-equivariant maps
$$ \phi: \Gamma\backslash G \to \cH_k  \ \ \ \sum_i (\Gamma
g_i)\otimes f_i \mapsto f_i. $$
We can embed
$$ \cA_*(\Gamma)\subset \hat\cA_*(\Gamma):= \Hom_\Gamma
(\C[\Gamma\backslash G], \cH_k), $$ where we think of
$\cA_*=\cH_*[\Gamma\backslash G]$ as polynomials in
$\Gamma\backslash G$ with $\cH_*$ coefficients, and of $\hat\cA_*=
\cH_* [[ \Gamma\backslash G ]]$ as formal power series, that is,
$\Gamma$-equivariant maps $\phi: \Gamma \backslash G \to \cH_k$.
There is on $\cA_*(\Gamma)$ an associative multiplication (\cf
\cite{Co-Mosc1}), which makes
$\cA_*(\Gamma)$ into a noncommutative ring. This is given by a
convolution product. For any $\phi \in \cA_k(\Gamma)$, we have
$\phi_g = \phi_{\gamma g}$, with $\phi_g=0$ off a finite subset of
$\Gamma \backslash G$, and $\phi_g | \gamma = \phi_{g\gamma}$, so
these terms are left $\Gamma$-invariant and right
$\Gamma$-equivariant. For $\phi \in \cA_k(\Gamma)$ and $\psi \in
\cA_\ell (\Gamma)$ we then define the convolution product as
\begin{equation}\label{convAprod}
( \phi * \psi )_g := \sum_{(g_1,g_2)\in G\times_\Gamma G, g_1 g_2
=g} (\phi_{g_1} | g_2) \phi_{g_2},
\end{equation}

The algebra $\cA_*(\Gamma)$ constructed this way has two
remarkable subalgebras.
\begin{itemize}
\item $\cA_0(\Gamma)= \C [\Gamma \backslash G / \Gamma]$ is the
algebra $\bT$ of Hecke operators.
\item $\cM_k(\Gamma)\subset \cA_k(\Gamma)$ also gives a subalgebra
$\cM_*(\Gamma) \subset \cA_*(\Gamma)$.
\end{itemize}

In particular observe that all the coefficients $\phi_g$ are
modular forms. In fact, they satisfy $\phi_g | \gamma =
\phi_{g\gamma}$, hence, for $\gamma \in \Gamma$, this gives
$\phi_g | \gamma = \phi_g$.

Notice however that the convolution product on $\cA_*$ does not
agree with the Hecke action, namely the diagram
\begin{eqnarray*}
\diagram \cM_*(\Gamma) \otimes \bT \dto^{\iota}\rto^{H} &
\cM_*(\Gamma) \dto^{\iota} \\
\cA_*\otimes \cA_* \rto^{*}& \cA_*,
\enddiagram
\end{eqnarray*}
with $\iota$ the inclusion of subalgebras and $H$ the Hecke
action, is {\em not} commutative, nor is the symmetric one
\begin{eqnarray*}
\diagram \bT  \otimes  \cM_*(\Gamma) \dto^{\iota}\rto^{H} &
\cM_*(\Gamma) \dto^{\iota} \\
\cA_*\otimes \cA_* \rto^{*}& \cA_*.
\enddiagram
\end{eqnarray*}
To get the correct Hecke action on modular forms from the algebra
$\cA_*(\Gamma)$, one needs to introduce the augmentation map
$$ \epsilon : \C [\Gamma \backslash G] \to \C $$
extended to a map
$$ \epsilon\otimes 1 : \C[\Gamma \backslash G] \otimes \cH_k \to
\cH_k  \ \ \ \sum [g]\otimes \phi_g \mapsto \sum \phi_g. $$

One then obtains a commutative diagram
\begin{eqnarray*}
\diagram \bT  \otimes  \cM_*(\Gamma) \dto^{\iota}\rto^{H} &
\cM_*(\Gamma)  \\
\cA_*\otimes \cA_* \rto^{*}& \cA_* \uto^{\epsilon \otimes 1}.
\enddiagram
\end{eqnarray*}

In \cite{CoMoHopf} Connes and Moscovici introduced a Hopf algebra $\cH_1$
associated to the transverse geometry of codimension one foliations.
This is the universal enveloping algebra of a Lie algebra with basis
$\{ \cX,\cY,\delta_n \, n\geq 1 \}$ satisfying, for $n,k,\ell \geq 1$,
\begin{equation}\label{Liebra}
[\cY,\cX]=\cX, \,\,\, [\cY,\delta_n] = n\,\delta_n, \,\,\,
[\cX,\delta_n] =\delta_{n+1}, \,\,\, [\delta_k,\delta_\ell]=0,
\end{equation}
with coproduct an algebra homomorphism $\Delta: \cH_1 \to \cH_1
\otimes \cH_1$ satisfying
\begin{equation}\label{coprod-H1}
\begin{array}{ll}
\Delta \cY=& \cY \otimes 1 + 1 \otimes \cY, \\ \Delta \delta_1 =&
\delta_1 \otimes 1 + 1 \otimes \delta_1, \\ \Delta \cX =& \cX\otimes 1
+ 1 \otimes \cX + \delta_1 \otimes \cY,
\end{array}
\end{equation}
antipode the anti-isomorphism satisfying
\begin{equation}\label{antipode-H1}
S(\cY)=-\cY, \,\,\,\, S(\cX)=-\cX+\delta_1\cY, \,\,\,\,
S(\delta_1)=-\delta_1,
\end{equation}
and co-unit $\epsilon(h)$ the constant term of $h\in \cH_1$.

The Hopf algebra $\cH_1$ acts as symmetries of the modular Hecke
algebras. This is a manifestation of the general fact that,
while symmetries of ordinary commutative spaces are encoded by group
actions, symmetries of noncommutative spaces are given by Hopf
algebras.

By comparing the actions of the Hopf algebra $\cH_1$, it is possible to
derive an analogy (\cf \cite{Co-Mosc1}) between the modular Hecke
algebras and the crossed product algebra of the action of a discrete
subgroup of ${\rm Diff}(S^1)$ on polynomial functions on the frame
bundle of $S^1$.

In fact, for $\Gamma$ a discrete subgroup of ${\rm Diff}(S^1)$, and
$X$ a smooth compact 1-dimensional manifold, consider as in
\cite{Co-Mosc1} the algebra
\begin{equation}\label{Jet-alg}
\cA_\Gamma = C^\infty_c (J^1_+(X)) \rtimes \Gamma,
\end{equation}
where $J^1_+(X)$ is the oriented 1-jet bundle. This has an action of
the Hopf algebra $\cH_1$ by
\begin{equation}\label{H1act1}
\begin{array}{ll}
\cY(f U_\phi^*)=& y_1 \frac{\partial f}{\partial y_1} U_\phi^* \\[3mm]
\cX(f U_\phi^*)=& y_1 \frac{\partial f}{\partial y} U_\phi^* \\[2mm]
\delta_n (f U_\phi^*)=& y_1^n \frac{d^n}{d y^n}\left( \log
\frac{d\phi}{dy} \right) f U_\phi^*,
\end{array}
\end{equation}
with coordinates $(y,y_1)$ on $ J^1_+(X)\simeq X\times \R^+$.
The trace $\tau$ defined by the volume form
\begin{equation}\label{trace-tau}
\tau(f U_\phi^*)= \left\{ \begin{array}{ll} \int_{J^1_+(X)} f(y,y_1)
\, \frac{dy\wedge dy_1}{y_1^2} & \phi=1 \\[3mm]
0 & \phi\neq 1 \end{array}\right.
\end{equation}
satisfies
\begin{equation}\label{tau-nu}
\tau(h(a)) = \nu(h) \tau(a) \, \, \, \, \forall h \in\cH_1,
\end{equation}
with $\nu\in \cH_1^*$ satisfying
\begin{equation}\label{nu}
\nu(\cY)=1, \,\,\,\,\,\, \nu(\cX)=0, \,\,\,\,\,\, \nu(\delta_n)=0.
\end{equation}
The twisted antipode $\tilde S=\nu * S$ satisfies $\tilde S^2=1$ and
\begin{equation}\label{twistS}
\tilde S(\cY)=-\cY+1, \,\,\,\,\,\, \tilde S(\cX)=-\cX +\delta_1 \cY,
\,\,\,\,\,\, \tilde S(\delta_1)=-\delta_1.
\end{equation}

The Hopf cyclic cohomology of a Hopf algebra is another fundamental
tool in noncommutative geometry, which was developed by Connes and
Moscovici in \cite{CoMoHopf}. They applied it to the
computation of the local index formula for tranversely
hypoelliptic operators on foliations. An action of a Hopf algebra on
an algebra induces a characteristic map from the Hopf cyclic
cohomology of the Hopf algebra to the cyclic cohomology of the
algebra, hence the index computation can be done in terms of Hopf
cyclic cohomology. The periodic Hopf cyclic cohomology of the Hopf
algebra of transverse geometry is related to the Gelfand-Fuchs
cohomology of the Lie algebra of formal vector fields \cite{CoMoHopf2}.

In the case of the Hopf algebra $\cH_1$, there are three basic
cyclic cocycles, which in the original context of transverse geometry
correspond, respectively, to the Schwarzian derivative, the
Godbillon-Vey class, and the transverse fundamental class.

In particular, the Hopf cyclic
cocycle associated to  the Schwarzian derivative is of the form
\begin{equation}\label{delta2'}
\delta_2':= \delta_2 - \frac 12 \delta_1^2
\end{equation}
with
\begin{equation}\label{delta2'act}
 \delta_2' (f U_\phi^*)= y_1^2 \left\{ \phi(y); y\right\} \, f
U_\phi^*
\end{equation}
\begin{equation}\label{schw-der}
\left\{ F ; x \right\} : = \frac{d^2}{dx^2} \left( \log \frac{dF}{dx}
\right) - \frac 12 \left( \frac{d}{dx} \left( \log \frac{dF}{dx}
\right) \right)^2.
\end{equation}

The action of the Hopf algebra $\cH_1$ on the modular Hecke algebra
described in \cite{Co-Mosc1} involves the natural derivation on
the algebra of modular forms initially introduced by Ramanujan, which
corrects the ordinary
differentiation by a logarithmic derivative of the Dedekind $\eta$
function,
\begin{equation}\label{Heck-rep}
\cX:= \frac{1}{2\pi i} \frac{d}{dz} - \frac{1}{2\pi i}
\frac{d}{dz}(\log \eta^4) \,\, \cY, \ \ \
\cY (f) =\frac k2 \, f, \,\,\, \forall f \in \cM_k .
\end{equation}
The element $\cY$ is the grading operator
that multiplies by $k/2$ forms of weight $k$, viewed as sections of
the $(k/2)$th power of the line bundle of 1-forms. The element
$\delta_1$ acts as multiplication by a form-valued cocycle on
$\GL^+_2(\Q)$, which measures the
lack of invariance of the section $\eta^4 dz$.
More precisely, one has the following action of
$\cH_1$ (Connes--Moscovici \cite{Co-Mosc1}):

\begin{thm}\label{actH1hecke}
There is an action of the Hopf algebra $\cH_1$ on the modular Hecke
algebra $\cA(\Gamma)$ of level $\Gamma$, induced by an action on
$\cA_{G^+(\Q)}:= \cM \rtimes G^+(\Q)$, for
$\cM=\varinjlim_{N\to\infty} \cM(\Gamma(N))$, of the form
\begin{equation}\label{actH1GQ}
\begin{array}{ll}
\cY(fU^*_\gamma)= & \cY(f)\, U^*_\gamma \\[3mm]
\cX(fU^*_\gamma)= & \cX(f)\, U^*_\gamma \\[2mm]
\delta_n(fU^*_\gamma)= & \frac{d^n}{dZ^n} \left( \log \frac{d(Z |_0
\gamma)}{dZ} \right) (dZ)^n \, fU^*_\gamma,
\end{array}
\end{equation}
with $\cX(f)$ and $\cY(f)$ as in \eqref{Heck-rep}, and
\begin{equation}\label{coordZ}
 Z(z) = \int_{i\infty}^z \eta^4 dz.
\end{equation}
\end{thm}

The cocycle \eqref{delta2'} associated to the Schwarzian derivative is
represented by an inner derivation of $\cA_{G^+(\Q)}$,
\begin{equation}\label{delta2'heck}
\delta_2'(a) = [a, \omega_4],
\end{equation}
where $\omega_4$ is the weight four modular form
\begin{equation}\label{omega4}
\omega_4 = -\frac{E_4}{72}, \,\,\,\, \text{ with } \,\,\, E_4(q)= 1 +
240 \sum_{n=1}^\infty n^3 \frac{q^n}{1-q^n}, \,\,\,\,\, q=e^{2\pi i z},
\end{equation}
which is expressed as a Schwarzian derivative
\begin{equation}\label{Somega4}
 \omega_4 = (2\pi i)^{-2} \, \{ Z; z\}.
\end{equation}

This result is
used in \cite{Co-Mosc1} to investigate perturbations of the Hopf algebra
action. The freedom one has in modifying the action by
a 1-cocycle corresponds exactly to the data introduced by Zagier
in \cite{Zagier}, defining "canonical" Rankin-Cohen algebras,
with the derivation $\partial$ and the element $\Phi$ in Zagier's
notation corresponding, respectively, to the action of the generator
$\cX$ on modular forms and to $\omega_4=2\Phi$.

The cocycle associated to the Godbillon-Vey class is described in
terms of a 1-cocycle on $\GL^+_2(\Q)$ with values in
Eisenstein series of weight two, which measures the lack of
$\GL^+_2(\Q)$-invariance of the connection associated to the
generator $\cX$. The authors derive from this an arithmetic presentation
of the rational Euler class in $H^2(\SL_2(\Q),\Q)$ in
terms of generalized Dedekind sums.

The cocycle associated to the transverse fundamental class, on the
other hand, gives rise to a natural extension of the first
Rankin-Cohen bracket \cite{Zagier} from modular
forms to the modular Hecke algebras.

\medskip

Rankin--Cohen algebras can be treated in different perspectives:
Zagier introduced them and studied them with a
direct algebraic approach (\cf \cite{Zagier}). There appears to be an
interesting and deep connection to vertex operator algebras, which
manifests itself in a form of duality between these two types of
algebras.

The Rankin--Cohen brackets are a family of brackets
$[f,g\,]^{(k,\ell)}_n$ for $n\geq 0$, defined for $f,g \in R$, where
$R$ is a graded ring with a derivation $D$. For
$R=\oplus_{k\geq 0} R_k$, $D: R_k \to R_{k+2}$, $f\in R_k$ $g\in
R_\ell$, the brackets $[,]_n : R_k\otimes R_\ell \to R_{k+\ell+2}$
are given by
\begin{equation}\label{RCbraD}
[f,g\,]^{(k,\ell)}_n = \sum_{r+s=n} (-1)^r \left(\begin{array}{c}
n+k-1 \\ s \end{array}\right) \left(\begin{array}{c} n+\ell -1 \\
r \end{array}\right) D^r f \, D^s g.
\end{equation}

These Rankin--Cohen brackets induced by $(R_*, D) \Rightarrow
(R_*, [,]_*)$ give rise to a {\em standard Rankin--Cohen algebra}, in
Zagier's terminology (\cf \cite{Zagier}).
There is an isomorphism of categories between graded rings with a
derivation and standard Rankin--Cohen algebras.

In the case of Lie algebras, one can first define a {\em standard
Lie algebra} as the Lie algebra associated to an associative
algebra $(A,*)\Rightarrow (A,[,])$ by setting $[X,Y]= X*Y-Y*X$ and
then define an abstract Lie algebra as a structure $(A,[,])$ that
satisfies {\em all} the algebraic identities satisfied by a
standard Lie algebras, though it is not necessarily induced by an
associative algebra. It is then a theorem that the antisymmetry of
the bracket and the Jacobi identity are sufficient to determine
all the other algebraic identities, hence one can take these as a
definition of an abstract Lie algebra.

Just as in the case of Lie algebras, we can define a Rankin--Cohen
algebra $(R_*, [,]_*)$ as a graded ring $R_*$ with a family of
degree $2n$ brackets $[,]_n$ satisfying all the algebraic
identities of the standard Rankin--Cohen algebra. However, in this
case there is no simple set of axioms that implies all the
algebraic identities.

The motivation for this structure lies in the fact that there is a
very important example of a Rankin--Cohen algebra which is in fact
non-standard. The example is provided by modular forms (\cf
\cite{Zagier}).

If $f\in \cM_k$ is a modular form satisfying
$$ f\left( \frac{az+b}{cz+d} \right) = (cz+d)^k f(z), $$
then it derivative is no longer a modular form, due to the
presence of the second term in
$$ f' \left(\frac{az+b}{cz+d} \right) =(cz+d)^{k+2} f(z) + kc
(cz+d)^{k+1} f(z).
$$
On the other hand, if we have $f\in \cM_k$ and $g\in \cM_\ell$,
the bracket
$$ [f,g](z) := \ell f'(z) g(z) - k f(z) g'(z) $$
is a modular form in $\cM_{k+\ell +2}$. Similarly, we can define
an $n$-th bracket $[,]_n: \cM_k \otimes \cM_\ell \to
\cM_{k+\ell+2}$. Here's the first few brackets:
$$ \begin{array}{ll} [f,g]_0 =& f g \\[2mm]
[f,g]_1 =& k f g' - \ell f' g \\[2mm]
[f,g]_2 = & \left(\begin{array}{c} k+1\\ 2\end{array}\right) f g''
- (k+1) (\ell +1) f'g' + \left(\begin{array}{c} \ell +1\\
2\end{array}\right)f'' g. \end{array} $$

Notice that for the graded ring of modular forms we have $\cM_*
(\Gamma)\subset \cH$, where $\cH$ is the vector space $\cH = {\rm
Hol}(\H)_{polyn}$ of holomorphic functions on the upper half
plane $\H$ with polynomial growth. This is closed under
differentiation and $(\cH, D)$ induces a standard Rankin--Cohen
algebra $(\cH, [,]_*)$. The inclusion $(\cM_*,[,]_*) \subset (\cH,
[,]_*)$ is not closed under differentiation but it is closed under
the brackets.

A way of constructing non-standard Rankin--Cohen algebras is
provided by Zagier's {\em canonical construction} (\cf
\cite{Zagier}). One considers
here the data $(R_*,D,\Phi)$, where $R_*$ is a graded ring with a
derivation $D$ and with a choice of an element $\Phi \in R_4$, the
{\em curvature}. One then defines the brackets by the formula
\begin{equation}\label{RCbraDPhi}
[f,g]^{(k,\ell)}_n = \sum_{r+s=n} (-1)^r \left(\begin{array}{c}
n+k-1 \\ s \end{array}\right) \left(\begin{array}{c} n+\ell -1 \\
r \end{array}\right)\, f_r \, g_s,
\end{equation}
where $f_0 = f$ and
\begin{equation}\label{DPhi}
f_{r+1} = D f_r + r(r+1) \Phi f_{r-1}.
\end{equation}
The structure $(R_*,[,]_*)$ obtained this way is a Rankin--Cohen
algebra (see \cite{Zagier}).

There is a gauge action on the curvature $\Phi$, namely, for any
$\varphi \in R_2$ the transformation $D \mapsto D'$ and $\Phi
\mapsto \Phi'$ with
$$ D'(f) = D(f) + k \varphi f $$
for $f \in \cM_k$ and
\begin{equation}\label{gaugePhi}
 \Phi' = \Phi + \varphi^2 - D(\varphi)
\end{equation}
give rise to {\em the same} Rankin--Cohen algebra. Thus, all the
cases where the curvature $\Phi$ can be gauged away to zero
correspond to the standard case.

\medskip

The modular form $\omega_4$ of \eqref{Somega4} provides
the curvature element $\omega_4 = 2\Phi$, and the gauge equivalence
condition \eqref{gaugePhi} can be rephrased in terms of Hopf
algebras as the freedom to change the $\cH_1$ action by a
cocycle. In particular (\cf \cite{Co-Mosc1}), for the specified
action, the resulting Rankin--Cohen structure is canonical but not
standard, in Zagier's terminology.

The 1-form $dZ=\eta^4 dz$ is, up to scalars, the only holomorphic
differential on the elliptic curve $E=X_{\Gamma(6)} \cong
X_{\Gamma_0(36)}$ of equation $y^2 = x^3 +1$, so that
$dZ=\frac{dx}{y}$ in Weierstrass coordinates.

\medskip

The Rankin--Cohen brackets on modular forms can be extended to
brackets $RC_n$ on the modular Hecke algebra, defined in terms of the
action of the Hopf algebra $\cH_1$ of transverse geometry. In fact,
more generally, it is shown in \cite{Co-Mosc2} that it is possible to
define such Rankin--Cohen brackets on any associative algebra $\cA$
endowed with an action of the Hopf algebra $\cH_1$ for which there
exists an element $\Omega\in \cA$ such that
\begin{equation}\label{innerdelta2'}
\delta_2'(a) = [ \Omega, a ], \,\,\,\, \forall a\in \cA,
\end{equation}
and with $\delta_2'$ as in \eqref{delta2'}, and
\begin{equation}\label{deltan0}
\delta_n(\Omega)=0, \,\,\, \forall n\geq 1.
\end{equation}

\smallskip

Under these hypotheses, the following result holds (Connes--Moscovici
\cite{Co-Mosc2}):

\begin{thm}\label{RCdef}
Suppose given an associative algebra $\cA$ with an action of the Hopf
algebra $\cH_1$ satisfying the conditions \eqref{innerdelta2'} and
\eqref{deltan0}.
\begin{enumerate}
\item There exists Rankin--Cohen brackets $RC_n$ of the form
\begin{equation}\label{RCnH1}
RC_n(a,b)= \sum_{k=0}^n \frac{A_k}{k!} (2\cY+k)_{n-k}(a)\,
\frac{B_{n-k}}{(n-k)!} (2\cY+n-k)_k (b),
\end{equation}
with $(\alpha)_r=\alpha(\alpha+1)\cdots(\alpha +r-1)$ and the
coefficients $A_{-1}=0$, $A_0=1$, $B_0=1$, $B_1=\cX$,
$$ A_{n+1}= S(\cX) A_n -n \Omega^0 \left(\cY - \frac{n-1}{2}\right)
A_{n-1}, $$
$$ B_{n+1}= \cX B_n - n \Omega \left(\cY-\frac{n-1}{2}\right) B_{n-1},
$$
and $\Omega^0$ the right multiplication by $\Omega$.
\item When applied to the modular Hecke algebra $\cA(\Gamma)$, with
$\Omega= \omega_4 = 2\Phi$, the
above construction yields brackets \eqref{RCnH1} that are completely
determined by their restriction to modular forms where they agree with
the Rankin--Cohen brackets \eqref{RCbraDPhi}.
\item The brackets \eqref{RCnH1} determine associative deformations
\begin{equation}\label{assoc-def}
a*b = \sum_n \hbar^n \, RC_n(a,b).
\end{equation}
\end{enumerate}
\end{thm}

\medskip

The first of the steps described in Section \ref{road}, namely
resolving the diagonal in $\cA(\Gamma)$, is not yet
done for the modular Hecke algebras and should shed light on the
important number theoretic problem of the interrelation of the Hecke
operators with the algebra structure given by the pointwise product.

The algebra $\cA(\Gamma)$ is deeply related to the algebra of the
space of two dimensional $\Q$-lattices of Section \ref{Qlatt}.

\section{Noncommutative moduli spaces, Shimura
varieties}\label{shimura}

It appears from the study of some significant cases that an important
source of interesting noncommutative spaces is provided by the
``boundary'' of classical (algebro-geometric) moduli spaces, when one
takes into account the possible presence of degenerations of classical
algebraic varieties that give rise to objects no longer defined within
the context of algebraic varieties, but which still make sense as
noncommutative spaces.

\smallskip

An example of algebro-geometric moduli spaces which is sufficiently
simple to describe but which at the same time exhibits a very rich
structure is that of the modular curves. The geometry of modular
curves has already appeared behind our discussion of the 2-dimensional
$\Q$-lattices and of the modular Hecke algebras, through an
associated class of functions: the modular functions that appeared in
our discussion of the arithmetic algebra for the quantum statistical
mechanical system of 2-dimensional $\Q$-lattices and the
modular forms in the modular Hecke algebras.

The modular curves, quotients of the hyperbolic plane $\H$ by the
action of a subgroup $\Gamma$ of finite index of $\SL_2(\Z)$, are complex
algebraic curves, which admit an arithmetic structure, as they are
defined over cyclotomic number fields $\Q(\zeta_N)$. They are also
naturally moduli spaces. The object they parameterize are elliptic
curves (with some level structure). The modular curves have an
algebro-geometric compactification obtained by adding finitely many
cusp points, given by the points in $\P^1(\Q)/\Gamma$. These
correspond to the algebro-geometric degeneration of the elliptic
curve to $\C^*$. However, in addition to these degenerations, one can
consider degenerations to noncommutative tori, obtained by a limit
$q\to \exp(2\pi i \theta)$ in the modulus $q=\exp(2\pi i \tau)$ of the
elliptic curve, where now $\theta$ is allowed to be also irrational.
The resulting boundary $\P^1(\R)/\Gamma$ is a noncommutative space (in
the sense of Section \ref{quotients}). It appeared in the string
theory compactifications considered in \cite{CDS}. The arithmetic
properties of the noncommutative spaces $\P^1(\R)/\Gamma$ were studied
in \cite{ManMar}, \cite{Mar2} \cite{Mar3}.

The modular curves, for varying finite index $\Gamma\subset
\SL_2(\Z)$, form a tower of branched coverings. The projective limit
of this tower sits as a connected component in the more refined {\em
adelic} version of the modular tower, given by the quotient
\begin{equation}\label{Sh2}
\GL_2(\Q)\backslash \GL_2(\A) /\C^*,
\end{equation}
where, as usual, $\A=\A_f\times \R$ denotes the adeles of $\Q$, with
$\A_f=\hat\Z\otimes \Q$ the finite adeles.

The space \eqref{Sh2} is also a moduli space. In fact, it belongs to
an important class of algebro-geometric moduli spaces of great
arithmetic significance, the Shimura varieties $Sh(G,X)$, where the
data $(G,X)$ are given by a reductive algebraic group $G$ and a
hermitian symmetric domain $X$. The pro-variety \eqref{Sh2}
is the Shimura variety $Sh(\GL_2,\H^\pm)$, where $\H^\pm=\GL_2(\R)
/\C^*$ is the union of upper and lower half plane in $\P^1(\C)$.

\smallskip

We have mentioned above that the spaces $\P^1(\R)/\Gamma$ describe
degenerations of elliptic curves to noncommutative tori. This type of
degeneration corresponds to degenerating a lattice $\Lambda=\Z+\Z\tau$
to a pseudolattice $L=\Z+\Z\theta$ (see \cite{Manin01RM} for a detailed
discussion of this viewpoint and its implications in noncommutative
geometry and in arithmetic). In terms of the space \eqref{Sh2}, it
corresponds to degenerating the archimedean component, namely
replacing $\GL_2(\R)$ by $M_2(\R)^\cdot=M_2(\R)\smallsetminus \{ 0
\}$, nonzero $2\times 2$-matrices. However, when one is working with
the adelic description as in \eqref{Sh2}, one can equally consider the
possibility of degenerating a lattice at the non-archimedean
components. This brings back directly the notion of $\Q$-lattices of
Section \ref{Qlatt}.

In fact, it was shown in \cite{CMR2} that the notions of 2-dimensional
$\Q$-lattices and commensurability can be reformulated in terms of
Tate modules of elliptic curves and isogeny. In these terms, the space
of $\Q$-lattices corresponds to non-archimedean degenerations of the
Tate module, which correspond to the ``bad quotient''
\begin{equation}\label{Sh2nc}
\GL_2(\Q)\backslash M_2(\A_f)\times \GL_2(\R) /\C^*.
\end{equation}
The combination of these two types of degenerations yields a
``noncommutative compactification'' of the Shimura variety
$Sh(\GL_2,\H^\pm)$ which is the algebra of the ``bad quotient''
\begin{equation}\label{Sh2ncC}
\GL_2(\Q)\backslash M_2(\A)^\cdot /\C^*,
\end{equation}
where $M_2(\A)^\cdot$ are the elements of $M_2(\A)$ with nonzero
archimedean component. One recovers the Shimura variety
$Sh(\GL_2,\H^\pm)$ as the set of classical points (extremal KMS states
at zero temperature) of the quantum statistical mechanical system
associated to the noncommutative space \eqref{Sh2nc} (\cf \cite{CoMa},
\cite{CMR2}).

\medskip

More generally, Shimura varieties are moduli spaces for certain
types of motives or as moduli spaces of Hodge structures (\cf \eg
\cite{Mil2}). A Hodge structure $(W,h)$ is a pair of a finite
dimensional $\Q$-vector space $W$ and a homomorphism $h:\bS\to
\GL(W_\R)$, of the real algebraic group
$\bS=\Res_{\C/\R}\bG_m$, with
$W_\R=W\otimes \R$. This determines a decomposition
$W_\R\otimes \C=\oplus_{p,q} W^{p,q}$ with
$\overline{W^{p,q}}=W^{q,p}$ and $h(z)$ acting on $W^{p,q}$ by
$z^{-p}\bar z^{-q}$. This gives a Hodge filtration and a weight
filtration $W_\R=\oplus_k W_k$ where $W_k=\oplus_{p+q=m} W^{p,q}$.
The Hodge structure $(W,h)$ has weight $m$ if $W_\R=W_m$. It is
rational if the weight filtration is defined over $\Q$. A Hodge
structure of weight $m$ is polarized if there is a morphism of
Hodge structures $\psi:W\otimes W\to \Q(-m)$, such that
$(2\pi i)^m\psi(\cdot, h(i)\cdot)$ is symmetric and positive
definite. Here $\Q(m)$ is the rational Hodge structure of weight
$-2m$, with $W=(2\pi i)^m \Q$ with the action $h(z)=(z\bar z)^m$.
For a rational $(W,h)$, the subspace of $W\otimes \Q(m)$ fixed by
the $h(z)$, for all $z\in \C^*$, is the space of ``Hodge cycles''.

One can then view Shimura varieties $Sh(G,X)$ as moduli spaces of
Hodge structures in the following way. Let $(G,X)$ be a Shimura datum
and $\rho:G \to \GL(V)$ a faithful representation. Since $G$ is
reductive, there is a finite family of tensors $\tau_i$ such that
\begin{equation}\label{tensors}
G=\{ g\in \GL(V): g\tau_i=\tau_i \}.
\end{equation}
A point $x\in X$ is by construction a $G(\R)$ conjugacy class of
morphisms $h_x: \bS\to G$, with suitable properties.

Consider data of the form $((W,h),\{ s_i \},\phi)$, where
$(W,h)$ is a rational Hodge structure, $\{ s_i \}$ a finite
family of Hodge cycles, and $\phi$ a $K$-level structure, for some
$K\subset G(\A_f)$, namely a $K$-orbit of $\A_f$-modules {\em
isomorphisms} $\phi: V(\A_f) \to W(\A_f)$, which maps $\tau_i$ to
$s_i$. Isomorphisms of such data are isomorphisms $f:W\to W'$ of
rational Hodge structures, sending $s_i\mapsto s_i'$,
and such that $f\circ \phi=\phi' k$,  for some $k\in K$.

We assume that there exists an isomorphism of $\Q$-vector spaces
$\beta:W\to V$ mapping $s_i \mapsto \tau_i$ and $h$ to $h_x$, for some
$x\in X$.

One denotes by ${\rm Hodge}(G,X,K)$ the set data $((W,h),\{ s_i \},\phi)$.
The Shimura variety $$Sh_K(G,X)=G(\Q)\backslash X\times
G(\A_f)/K$$ is the moduli space of the isomorphism classes of data
$((W,h),\{ s_i \},\phi)$, namely there is a map of ${\rm Hodge}(G,X,K)$
to $Sh_K(G,X)$ (seen over $\C$), that descends to a bijection on
isomorphism classes ${\rm Hodge}(G,X,K)/\sim$.

In such cases also one can consider degenerations of these data, both
at the archimedean and at the nonarchimedean components. One then
considers data $((W,h),\{ s_i \}, \phi,\tilde\beta)$, with
a non-trivial {\em homomorphism} $\tilde\beta: W\to V$, which is a
morphism of Hodge structures, and such that $\tilde\beta(\ell_{s_i})\subset
\ell_{\tau_i}$. This yields noncommutative spaces inside which the
classical Shimura variety sits as the set of classical points.

Quantum statistical machanical systems associated to Shimura varieties
have been recently studied by Ha and Paugam in \cite{HaPa}. Given a
faithful representation $\rho:G\to \GL(V)$ as above, there is an
``enveloping semigroup'' $M$, that is, a normal irreducible semigroup
$M\subset \End(V)$ such that $M^\times =G$. Such semigroup can be used
to encode the degenerations of the Hodge data described above.
The data $(G,X,V,M)$ then determine a noncommutative space which
describes the ``bad quotient'' $Sh^{nc}_K(G,X)=G(\Q)\backslash X\times
M(\A_f)$ and is a moduli space for the
possibly non-invertible data $((W,h),\{ s_i \},\phi)$.
Its ``set of classical points'' is the Shimura variety
$Sh(G,X)$.  The actual construction of the algebras involves
some delicate steps, especially to handle the presence of stacky
singularities (\cf \cite{HaPa}).

\section{The adeles class space and the spectral realization}\label{adeles}

In this section we   describe a noncommutative space, the adele
class space $X_\K$, associated to any global field $\K$, which leads
to a spectral realization of the zeros of the Riemann zeta function
for $\K=\Q$ and more generally of $L$-functions associated to Hecke
characters. It also gives a geometric interpretation of the
Riemann-Weil explicit formulas of number theory as a trace formula.
This space is closely related for $\K=\Q$ with the space of
commensurability classes of $\Q$-lattices described above. Rather
than starting directly with its description we   first put the
problem of finding the geometry of the set of prime numbers in the
proper perspective.

\medskip
{\bf The set of primes}
\medskip

One of the main problems of arithmetic is to understand the
distribution of the set of prime numbers
$$
\{2, 3, 5, 7, 11, 13, 17, 19, 23, 29, 31, 37, 41, 43, 47, 53, 59,
61, 67, 71, \ldots \}
$$
as a subset of the integers. To that effect one introduces the
counting function
$$
\pi(n)=\,{\rm number}\;  {\rm of}\;  {\rm primes}\; p\leq n
$$
The problem is to understand the behavior of $\pi(x)$ when $x\to
\infty$. One often hears that the problem comes from the lack of a
``simple" formula for $\pi(x)$. This is not really true and for
instance in 1898 H. Laurent \cite{laurent} gave the following
formula whose validity is an easy exercise in arithmetic,
\begin{small}
$$
\pi(n)=\, 2+\,\sum_{k=5}^n\,\frac{e^{2\pi i \Gamma(k)/k}-1}{e^{-2\pi i
/k}-1} \,,
$$
\end{small}
where $ \Gamma(k)=(k-1)!$ is the Euler Gamma function.

\begin{figure}
\begin{center}
\includegraphics[scale=0.80]{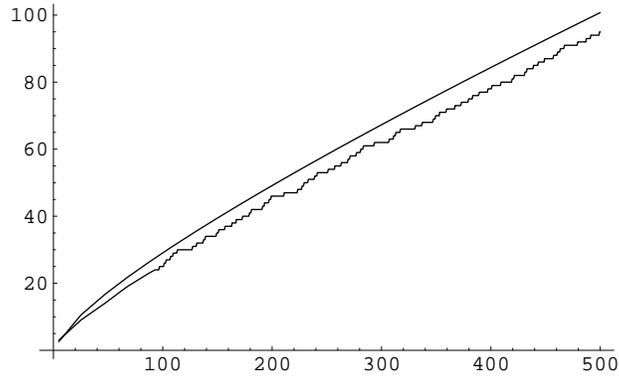}
\end{center}
\caption{Graphs of $\pi(x)$ and  Li$(x)$ \label{pili1}}
\end{figure}

The problem with this formula is that it has no bearing on the
asymptotic expansion of $\pi(x)$ when $x\to \infty$. Such an
expansion was guessed by Gauss in the following form,
$$
\pi(x)=\, \int_0^x\,\frac{du}{\log(u)}\,+\;R(x)
$$
where the integral logarithm admits the asymptotic expansion,
$$Li(x) =\int_0^x\,\frac{du}{\log(u)}
\sim \sum \,(k-1)!\,\frac{x}{\log(x)^k}$$

\begin{figure}
\begin{center}
\includegraphics[scale=0.80]{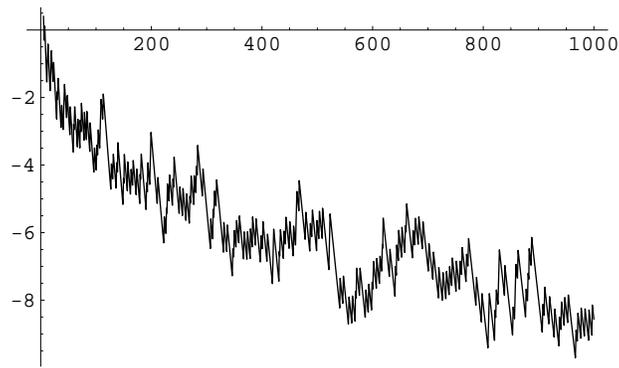}
\end{center}
\caption{Graph of $\pi(x)-$ Li$(x)$ \label{pili2}}
\end{figure}

The key issue then is the size of the remainder $R(x)$.

\medskip
 {\bf The Riemann hypothesis}
\medskip

It asserts that this size is governed by
\begin{equation}\label{rh}
R(x)=\,O(\sqrt{x}\,\log(x)) .
\end{equation}
The Riemann Hypothesis is in fact a conjecture on the zeros of the
zeta function
\begin{equation}\label{zeta}
\zeta(s)=\,\sum_1^\infty\,n^{-s} ,
\end{equation}
whose definition goes back to Euler,
who showed the fundamental factorization
\begin{equation}\label{eulerprod}
 \zeta(s)=\,\prod_{\mathcal P}\,(1-\,p^{-s})^{-1}\,.
\end{equation}

It extends to a meromorphic function in the whole complex plane $\C$
and fulfills the functional equation
\begin{equation}\label{functequ}
\pi^{-s/2}\,\Gamma(s/2)\,\zeta(s)=\,\pi^{-(1-s)/2}\,\Gamma((1-s)/2)\,\zeta(1-s),
\end{equation}
so that the function
\begin{equation}\label{zetaq}
\zeta_\Q(s)=\pi^{-s/2}\,\Gamma(s/2)\,\zeta(s)
\end{equation}
admits the symmetry $s\mapsto 1-s$.

The Riemann conjecture asserts that all zeros of  $\zeta_\Q$ are on
the critical line $\frac{1}{2} +\,i\,\R$. The reason why the
location of the zeros of $\zeta_\Q$ controls the size of the
remainder in \eqref{rh} is the explicit formulas that relate primes
with these zeros. Riemann proved the following first instance of an
``explicit formula"
\begin{equation}\label{explicit1}
\pi'(x)=\,Li(x)-\,\sum_\rho\,Li(x^\rho)+\,\int_x^\infty\,\frac{1}{u^2-1}\;\frac{du}{u
\log u}-\,\log 2
\end{equation}
where the sum is over the non-trivial (\ie complex) zeros of the
zeta function, and where the function $\pi'(x)$ is given by
$$
\pi'(x)=\,\pi(x)\,+\,\,\frac{1}{2}\,\pi(x^{\frac{1}{2}})+
\,\frac{1}{3}\,\pi(x^{\frac{1}{3}})\,+\cdots
$$
which gives by the Moebius inversion formula
$$
\pi(x)=\,\sum \,\mu(m)\,\frac{1}{m}\,\pi'(x^{\frac{1}{m}}) .
$$

\bigskip
\begin{figure}
\begin{center}
\includegraphics[scale=0.80]{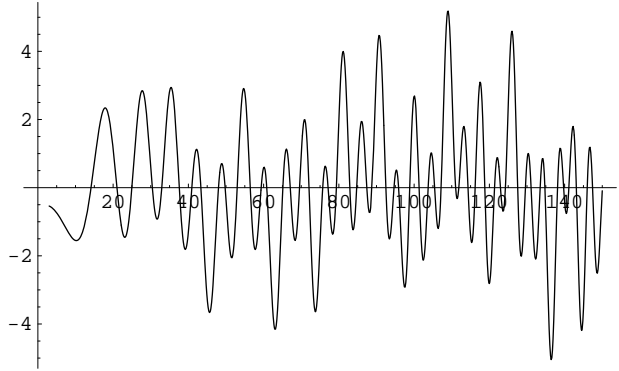}
\end{center}
\caption{Zeros of zeta \label{zeros}}
\end{figure}

\medskip
{\bf The generalized Riemann hypothesis}
\medskip

The explicit formulas of Riemann were put in more modern form by A.
Weil, as
\begin{equation}\label{explicit}
\,\wh h(0) + \wh h(1) -\sum_\rho\,\wh h(\rho)=\,\sum_{v}
\int'_{\K^*_v} \, \frac{h(u^{-1})}{
 |1-u|}\, d^* u \, ,
\end{equation}
where $\K$ is now an arbitrary global field, $v\in \Sigma_\K$ varies
among the places of $\K$ and the integral is taking place over the
locally compact field $\K_v$ obtained by completion of $\K$ at the
place $v$. Also $\int'$ is  the pairing with the distribution on
$\K_v$ which agrees with ${du \over \vert 1-u \vert}$ for $u \not= 1$
and whose Fourier transform (relative to a selfdual choice of
additive characters $\alpha_v$) vanishes at $1$. By definition a
{\em global field} is a (countable) discrete cocompact subfield in a
locally compact ring. This ring depends functorially on $\K$ and is
called the ring $\A_\K$ of adeles of $\K$. The quotient

\begin{equation}\label{ideleclass}
C_\K=\,{\rm GL}_1(\A_\K)/{\rm GL}_1(\K)
\end{equation}
is the locally compact group of idele classes of $\K$ which plays a
central role in class field theory. In Weil's explicit formula the
test function $h$ is in the Bruhat-Schwartz space $\cS(C_\K)$. The
multiplicative groups ${\rm GL}_1(\K_v)=\K_v^*$ are embedded
canonically as cocompact subgroups of $C_\K$. The sum on the left
hand side takes place over the zeros of $L$-functions associated to
Hecke characters. The function $\wh h$ is the Fourier transform of
$h$. The generalized Riemann conjecture asserts that all the zeros
of these $L$-functions are on the critical line $\frac{1}{2}
+\,i\,\R$. This was proved by Weil when the global field $\K$ has
non-zero characteristic, but remains open in the case when $\K$ is of
characteristic zero, in which case it is a number field \ie a finite
algebraic extension of the field $\Q$ of rational numbers.

\medskip
 {\bf
Quantum Chaos $\to$ Riemann Flow ? }
\medskip

For $E>0$ let $N(E)$ be the number of zeros of the Riemann zeta
function $\zeta_\Q$ whose imaginary parts are in the open interval
$]0,E[$. Riemann proved that the step function $N(E)$ can be written
as the sum
$$N(E)=\,  \langle N(E) \rangle  + N_{\rm osc} (E)$$
of a smooth approximation $ \langle N(E) \rangle$  and a purely
oscillatory function $N_{\rm osc} (E)$ and gave the following
explicit form

\begin{equation}\label{riemcount}  \langle N(E) \rangle = \,\frac{E}{2\,\pi}\,(\log
\frac{E}{2\,\pi} -1) + \frac{7}{8} +\,o(1)
\end{equation}

for the smooth approximation.

There is a striking analogy between the behavior of the step
function $N(E)$ and that of the function counting the number of
eigenvalues of the Hamiltonian $H$ of the quantum system obtained
after quantization of a chaotic dynamical system, which is at
centerstage in the theory of {\em quantum chaos}. A  comparison of
the asymptotic expansions of the oscillatory terms in both cases,
namely

$$N_{\rm osc} (E) \sim \frac{1}{ \pi} \sum_{\gamma_p} \sum_{m=1}^{\infty}
\frac{1}{ m} \, \frac{1}{ 2{\rm sh} \left( \frac{m\lambda_p }{
2}\right)} \, \sin (m \, E \, T_{\gamma}^{\#})
$$

for the quantization of a chaotic dynamical system, and

$$
N_{\rm osc} (E) \sim \frac{-1 }{ \pi} \sum_p \sum_{m=1}^{\infty}
\frac{1}{ m} \, \frac{1 }{ p^{m/2}} \, \sin \, (m \, E \, \log \, p)
\,
$$

for the Riemann zeta function, gives precious indications on the
hypothetical {\em Riemann flow} that would make it possible to identify the
zeros of zeta as the spectrum of an Hamiltonian. For instance the
periodic orbits of the flow should be labeled by the prime numbers
and the corresponding periods $T_p$ should be given by the $\log p$.
However a closer look reveals an overall minus sign that forbids any
direct comparison.

\bigskip
\begin{figure}
\begin{center}
\includegraphics[scale=0.90]{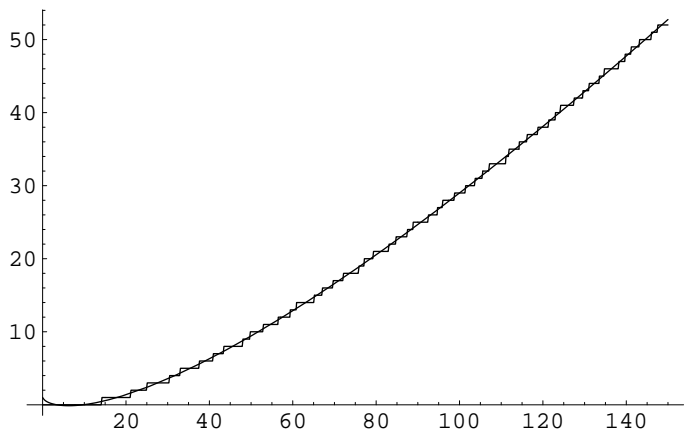}
\end{center}
\caption{Counting zeros of zeta \label{zeros1}}
\end{figure}

\medskip
 {\bf
Spectral realization as an absorption spectrum}
\medskip

The above major sign obstruction was bypassed in \cite{CoRH} using
the following basic distinction between observed spectra in physics.
When the light coming from a hot chemical element is decomposed
through a prism, it gives rise to bright emission lines on a dark
background, and the corresponding frequencies are a signature of its
chemical composition. When the light coming from a distant star is
decomposed through a prism, it gives rise to dark lines, called
absorption lines,  on a white background. The spectrum of the light
emitted by the sun was the first observed example of an absorption
spectrum. In this case the absorption lines were discovered by
Fraunhofer. The chemicals in the outer atmosphere of the star absorb
the corresponding frequencies in the white light coming from the
core of the star.

\bigskip
\begin{figure}
\begin{center}
\includegraphics[scale=0.8]{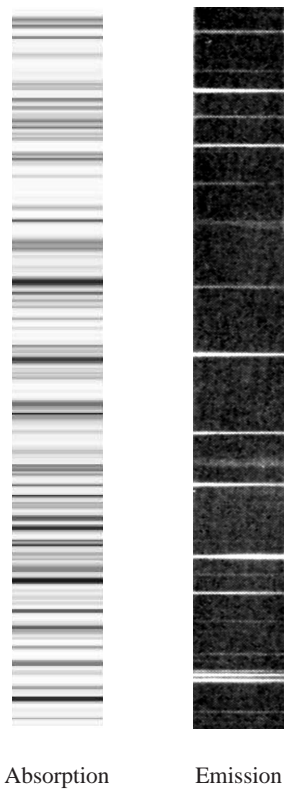}
\end{center}
\caption{The two kinds of Spectra \label{spectra}}
\end{figure}

The simple idea then is that, because of the minus sign above, one
should look for the spectral realization of the zeros of zeta not as
a usual emission spectrum but as an absorption spectrum. Of course
by itself this idea does not suffice to get anywhere since one needs
the basic dynamical system anyway. The  adele class space, namely the
quotient
$$
X_\K=\A_\K/ \K^*
$$
introduced in \cite{CoRH}, does the job as shown there. The action of
the idele class group $C_\K$ on the adele class space is simply given
by multiplication. In particular the idele class group $C_\K$ acts on
the suitably defined Hilbert space $L^2(X_\K)$ and the  zeros of
$L$-functions give the absorption spectrum, with non-critical zeros
appearing as resonances.

\medskip

Exactly as adeles, the adele class space $X_\K$ involves all the
places of $\K$. If in order to simplify, one restricts to a finite
set of places one still finds a noncommutative space but one can
analyze the action of the analogue of $C_\K$ and compute its trace
after performing a suitable cutoff (necessary in all cases to see
the missing lines of an absorption spectrum). One gets a trace
formula
$$
{\rm Trace} \, (R_{\Lambda} \, U(h)) = 2h (1) \log' \Lambda +
\sum_{v \in S} \int'_{\K^*_v} \frac{h(u^{-1})}{  \vert 1-u \vert} \,
d^* u + o(1),
$$
 where the terms in the right hand side are exactly the same as  in Weil's explicit
 formula \eqref{explicit}. This is very encouraging since at least it
 gives  geometric meaning to the complicated terms of \eqref{explicit}
 as the contributions of the periodic orbits to the computation of
 the trace.

In particular it gives a perfect interpretation  of the smooth
function $  \langle N(E) \rangle $ approximating the counting $N(E)$
of the zeros of zeta, from counting the number of states of the
 one dimensional quantum system with Hamiltonian
 $$
h (q, p) =  \, 2\pi q\,p
$$
which is just the generator of the scaling group. Indeed the
function $\frac{E}{2\,\pi}\,(\log \frac{E}{2\,\pi} -1)$ of Riemann
formula \eqref{riemcount} appears as the number of missing degrees
of freedom in the number of quantum states for the above system, as
one obtains from the simple computation of the area of the region
$$
B_+=\,\{(p,q)\in[0,\Lambda]^2\,;\,h(p,q)\leq E\}
$$

$${\rm Area}(B_+)=\frac{E}{2\pi} \times 2 \log \Lambda\,
 - \frac{E}{2\pi} \left( \log \frac{E}{2\pi} - 1 \right)
$$

while the term $\frac{E}{2\pi} \times 2 \log \Lambda$ corresponds to
the number of degrees of freedom of white light. A careful
computation gives not only the correction term of $\frac{7}{8}$ in
\eqref{riemcount} but all the remaining $o(1)$ terms.

\bigskip
\begin{figure}
\begin{center}
\includegraphics[scale=0.8]{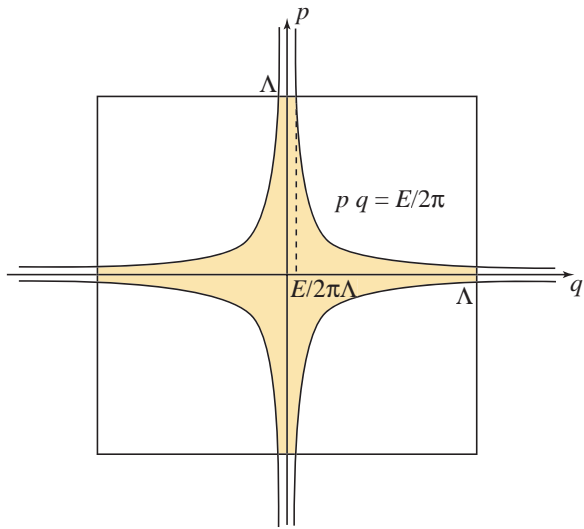}
\end{center}
\caption{Counting quantum states}\label{countingstates}
\end{figure}

Finally it was shown in \cite{CoRH} that the generalized Riemann
hypothesis is equivalent to the validity of a {\em global} trace
formula, but this is only one of many equivalent reformulations of the
Riemann hypothesis.\footnote{There is a running joke, inspired
by the European myth of Faust, about a mathematician trying to
bargain with the devil for a proof of the Riemann hypothesis ...}

\medskip

\section{Thermodynamics of endomotives and the Tehran program}\label{program}

In many ways the great virtue of a problem like RH comes from the
developments that it generates. At first sight it does not appear as
having any relation with geometry, and its geometric nature
gradually emerged in the twentieth century mainly because of the
solution of Weil in the case of global fields of positive
characteristic.

 We    outline a program, current joint
work of Katia Consani and the two authors, to adapt Weil's proof for
the case of global fields of positive characteristic to the case of
number fields.\footnote{This program was first announced in a
lecture at IPM Tehran in September 2005, hence we refer to it as ``the Tehran
program''.}

\medskip
 {\bf
Function Fields}
\medskip

Given a global field $\K$ of positive characteristic, there exists a
finite field $\F_q$ and a  smooth projective curve $C$ defined over
$\F_q$ such that $\K$ is the field of $\F_q$-valued rational
functions on $C$. The analogue (Artin, Hasse, Schmidt) of the zeta
function is
$$
\zeta_\K(s)=\,\prod_{v\in\Sigma_K}\,(1-q^{-f(v)s})^{-1}
$$
where  $\Sigma_\K$ is the set of places of $\K$ and $f(v)$ is the
degree (see below) of the place $v\in \Sigma_\K$.

The functional equation takes the form
$$
q^{(g-1)(1-s)}\,\zeta_\K(1-s)=\,q^{(g-1)s} \,\zeta_\K(s)
$$
where $g$ is the genus of  $C$.

The analogue of the Riemann conjecture for such global fields was
proved by Weil (1942) who developed algebraic geometry in that
context. Weil's proof rests on two steps.

\begin{itemize}

\item (A) Explicit Formula
\medskip

\item (B) Positivity

\end{itemize}

Both are based on the geometry of the action of the Frobenius on the
set $C(\bar \F_q)$ of points of $C$ over an algebraic closure $\bar
\F_q$ of  $\F_q$. This set $C(\bar \F_q)$ maps canonically to the
set $\Sigma_\K$ of places of $\K$ and the degree of a place $v\in
\Sigma_\K$ is the number of points in the orbit of the
 Frobenius acting on the fiber of the projection
$$
C(\bar \F_q) \to \Sigma_\K\,.
$$
The analogue
$$\#\{C(\F_{q^j})\}=\,\sum\,(-1)^k\,
{\rm Tr}({\rm Fr}^{* j}|H^k_{{\rm et}}(\bar C,\Q_\ell))$$
of the Lefschetz fixed point formula makes it possible to compute
the number $\#\{C(\F_{q^j})\}$ of points with coordinates in the
finite extension $\F_{q^j}$ from the action of ${\rm Fr}^*$ in the
etale cohomology group $H^1_{{\rm et}}(\bar C,\Q_\ell)$, which does
not depend upon the choice of the $\ell$-adic coefficients
$\Q_\ell$.

This shows that the zeta function is a rational fraction
$$
\zeta_\K(s)=\,\frac{P(q^{-s})}{(1-q^{-s})(1-q^{1-s})}
$$

where the polynomial $P$ is the characteristic polynomial of  the
action of ${\rm Fr}^*$ in $H^1$.

The analogue of the Riemann conjecture for global fields of
characteristic  $p$ means that its eigenvalues  \ie the complex
numbers $\lambda_j$ of the factorization
$$
P(T)=\,\prod \,(1-\lambda_j\, T)
$$
are of modulus $|\lambda_j|=q^{1/2}$.

The main ingredient in the proof of Weil is the notion of
correspondence, given by divisors in $C\times C$. They can be viewed
as multivalued maps,
$$Z\;:\;C\to C,\ \ \ \  P\mapsto Z(P).$$

Two correspondences are equivalent if
they differ by a principal divisor,
$$
U\sim V \Leftrightarrow U-V=\,(f)
$$
The composition of correspondences is
$$
Z=\,Z_1 \star Z_2 \,, \ \ \  Z_1 \star Z_2(P)=\,Z_1 ( Z_2(P))
$$
and the adjoint is given using the transposition
$\sigma(x,y)=\,(y,x)$ by
$$
Z'=\, \sigma(Z)\,.
$$
The degree $d(Z)$ of a correspondence is defined, independently of a
generic point $P\in C$ by,
$$
d(Z)=\,Z \,\bullet \,(P\times C)\,,
$$
where $\bullet$ is the intersection number. One has a similar
definition of the codegree
$$ \quad d^{\,'}(Z)=\,Z
\,\bullet \,(C\times P) $$

 Weil defines the {\em Trace} of a
correspondence as follows
$$
{\rm Trace}(Z)=\,d(Z)+\,d^{\,'}(Z)-\,Z\bullet \Delta
$$
where $\Delta$ is the identity correspondence. The main step in
Weil's proof is

\smallskip

\begin{thm} (Weil)
The following positivity holds : ${\rm Trace}(Z \star Z')>0$ unless
$Z$ is a trivial class.
\end{thm}

\smallskip

 Clearly
if one wants to have any chance at imitating the steps of Weil's
proof of RH for the case of number fields one needs to have an
analogue of the points of $C(\bar \F_q)$ and the action of the
Frobenius, of the etale cohomology and of the unramified extensions
$\K\otimes_{\F_q}\,\F_{q^n}$ of $\K$.

\medskip
 {\bf
 Endomotives and Galois action}
\medskip

The adele class space $X_\K$ of a global field admits a natural
action of the idele class group $C_\K$ and as such is on the adelic
side of the class field theory isomorphism. In order to obtain a
description of this space which is closer to geometry one needs to
pass to the Galois side of class field theory. In the case
 $\K=\Q$, it is possible to present the adele class space in a fairly simple
 manner not involving adeles thanks to its intimate relation with
 the  space of $1$-dimensional $\Q$-lattices of section \ref{Qlatt}.
The direct interpretation of the action of the Galois group of $\bar
\Q/\Q$ on the values of fabulous states for the BC-system (section
\ref{Qlatt}) then suggests that one should be able to construct
directly the space $X_\Q$ with a canonical action of the Galois
group of $\bar \Q/\Q$.

This was done in \cite{CCM} thanks to an extension of the notion of
Artin motives, called {\em endomotives}. Following Grothendieck, one
can reformulate Galois theory over a field $\K$ as the equivalence of
the category of reduced commutative finite dimensional algebras over
$\K$ with the category of continuous actions of the Galois group $G$
of $\bar \K/\K$ on finite sets. By construction the algebra of the
BC-system is a crossed product of a commutative algebra $A$ by a
semi-group. When working over $\K=\Q$ which is essential in the
definition of fabulous states, the algebra $A$ is simply the group
ring $\Q[\Q/\Z]$ of the torsion group $\Q/\Z$. Thus
$$
A=\, \varinjlim \,A_n \,,\quad A_n=\,\Q[\Z/n\Z]
$$
and we are dealing with a projective limit of Artin motives. The key
point then is to keep track of the corresponding action of the
Galois group $G$ of $\bar \K/\K$, $\K=\Q$. The Galois-Grothendieck
correspondence associates to a reduced commutative finite
dimensional algebra $B$ over $\K$ the set of characters of $B$ with
values in $\bar \K$ together with the natural action of $G$. This
action is non-trivial for the algebras $A_n=\,\Q[\Z/n\Z]$ where it
corresponds to the cyclotomic theory.

One can then recover the Bost--Connes system with its natural Galois symmetry
in a conceptual manner which extends to the general context of
semigroup actions on projective systems of Artin motives. These
typically arise from self-maps of algebraic varieties.
Given a pointed algebraic
variety $(Y,y_0)$ over a field $\K$ and a countable unital abelian
semigroup $S$ of finite endomorphisms of $(Y,y_0)$, unramified over
$y_0\in Y$ one constructs a projective system of Artin motives $X_s$
over $\K$ from these data as follows. For $s\in S$, one sets
\begin{equation}\label{Xs}
X_s=\{ y\in Y:\, s(y)=y_0 \}.
\end{equation}
For a pair $s,s'\in S$, with $s'=sr$, the map $\xi_{s',s}: X_{sr}\to
X_s$ is given by
\begin{equation}\label{Xsr}
X_{sr} \ni y \mapsto r(y)\in X_s.
\end{equation}
This defines a projective system indexed by the semigroup $S$ itself
with partial order given by divisibility. We let $X=\varprojlim_s
X_s$.

Since $s(y_0)=y_0$, the base point $y_0$ defines a component $Z_s$
of $X_s$ for all $s\in S$. Let $\xi_{s',s}^{-1}(Z_s)$ be the inverse
image of $Z_s$ in $X_{s'}$. It is a union of components of $X_{s'}$.
This defines a projection $e_s$ onto an open and closed subset
$X^{e_s}$ of the projective limit $X$. One then shows (\cite{CCM})
that the semigroup $S$ acts on the projective limit $X$ by partial
isomorphisms $\rho_s: X \to X^{e_s}$ defined by the property that
\begin{equation}\label{endoS}
\xi_{su}(\rho_s(x))=\xi_u(x), \, \forall u\in S, \forall x\in X.
\end{equation}

The BC-system is obtained from the pointed algebraic variety
$(\bG_m(\Q),1)$ where  the affine group scheme $\bG_m$ is the
multiplicative group. The semigroup $S$ is the semigroup of non-zero
endomorphisms of $\bG_m$. These correspond to maps of the form
$u\mapsto u^n$ for some non-zero $n\in \Z$, and one restricts to
$n\in \N^*$.

In this class of examples one has an ``equidistribution'' property,
by which the uniform normalized counting measures $\mu_s$ on $X_s$
are compatible with the projective system and define a probability
measure on the limit $X$. Namely, one has
\begin{equation}\label{measlim}
\xi_{s',s} \mu_s = \mu_{s'},\ \ \ \  \forall s,s'\in S.
\end{equation}
This follows from the fact that the number of preimages of a point
under $s\in S$ is equal to $\deg s$. This provides exactly the data
which makes it possible to perform the thermodynamical analysis of such
endomotives. This gives a rather unexplored new territory since even
the simplest examples beyond the BC-system  remain to be
investigated. For instance let $Y$ be an elliptic curve defined over
$\K$. Let $S$ be the semigroup of non-zero endomorphisms of $Y$.
This gives rise to an example in the general class described above.
When the elliptic curve has complex multiplication, this gives rise
to a system which, in the case of a maximal order, agrees with the
one constructed in \cite{CMR}. In the case without complex
multiplication, this provides an example of a system where the
Galois action does not factor through an abelian quotient.

\medskip
 {\bf
 Frobenius as dual of the time evolution}
\medskip

The Frobenius is such a universal symmetry in characteristic $p$,
owing to the linearity of the map $x\mapsto x^p$ that it is very hard to
find an analogue of such a far reaching concept in characteristic
zero. As we   now explain, the  classification of type III
factors provides the basic ingredient which when combined with
cyclic cohomology makes it possible to analyze  the thermodynamics of a
noncommutative space and get an analogue of the action of the
Frobenius on etale cohomology.

The key ingredient is that noncommutativity generates a time
evolution at the ``measure theory" level.  While it had been long
known by operator
 algebraists that the theory of
 von-Neumann algebras represents
 a far reaching extension of measure theory, the
  main surprise which occurred at the
beginning of the seventies in (Connes \cite{Co_2}) following
Tomita's theory is that such an algebra $M$ inherits
 from its noncommutativity a god-given time evolution:
\begin{equation}
\delta\,: \quad \R \quad\longrightarrow \quad{\rm Out}\,M=\,{\rm
Aut}\,M/{\rm Inn}\,M
\end{equation}

where ${\rm Out}\,M ={\rm Aut}\,M/{\rm Inn}\,M$ is the quotient of
the group of automorphisms of $M$ by the normal subgroup of inner
automorphisms. This led in \cite{Co_2} to the reduction from type
III to type II and their automorphisms and eventually to the
classification of injective factors.

\noindent They are classified by their module,
\begin{equation}
{\rm Mod} (M) \mathop{\subset}_{\sim} \ \R_+^* \, , \label{eq:(20)}
\end{equation}
which is a virtual closed subgroup of $\R_+^*$ in the sense of
G.~Mackey, i.e. an ergodic action of $\R_+^*$, called the flow of
weights \cite{CT}. This invariant was first defined and used in
\cite{Co_2}, to show in particular the existence of hyperfinite
factors which are not isomorphic to Araki-Woods factors.

The ``measure theory" level \ie the set-up of von-Neumann algebras
does not suffice to obtain the relevant cohomology theory and one
needs to be given a weakly dense subalgebra $\cA\subset M$ playing
the role of smooth functions on the noncommutative space. This
algebra will play a key role when cyclic cohomology is used at a
later stage. At first one only uses its norm closure $A=\bar \cA$ in
$M$ and assumes that it is globally invariant under the modular
automorphism group $\sigma_t^\varphi$ of a faithful normal state
$\varphi$ on $M$. One can then proceed with the thermodynamics of
the $C^*$ dynamical system $(A,\sigma_t)$. By a very simple
procedure assuming that KMS states at low temperature are of type I,
one obtains a ``cooling morphism" $\pi$ which is a morphism of
algebras from the crossed product $\hat \cA=\,\cA
\rtimes_{\sigma}\R$ to a type I algebra of compact operator valued
functions on a canonical $\R^*_+$-principal bundle
$\tilde\Omega_\beta$ over the space $\Omega_\beta$ of type I
extremal KMS$_\beta$ states fulfilling a suitable regularity
condition (\cf \cite{CCM}). Any $\varepsilon\in \Omega_\beta$ gives
an irreducible representation $\pi_\varepsilon$ of $\cA$ and the
choice of its essentially unique extension to $\hat \cA$ determines
the fiber of the $\R^*_+$-principal bundle $\tilde\Omega_\beta$. The
cooling morphism is then given by,
\begin{equation}
\pi_{\varepsilon,H}(\int \,x(t)\,U_t\,dt)=\,\int
\,\pi_\varepsilon(x(t))\,e^{itH}\,dt .
\end{equation}
This morphism is equivariant for the dual action $\theta_\lambda \in
\Aut(\hat \cA)$ of $\R^*_+$,
\begin{equation}\label{dualaction}
\theta_\lambda(\int \,x(t)\,U_t\,dt)=\,\int
\,\lambda^{it}\,x(t)\,U_t\,dt.
\end{equation}

The key point is that the range of the morphism $\pi$ is contained
in an algebra of functions on $\tilde\Omega_\beta$ with values in
trace class operators. In other words modulo a Morita equivalence
one lands in the commutative world provided one lowers the
temperature.

The interesting space is obtained by ``distillation" and is simply
given by the cokernel of the cooling morphism $\pi$ but this does
not make sense in the category of algebras and algebra homomorphisms
since the latter is not even an additive category. This is where
cyclic cohomology enters the scene : the category of cyclic modules
is an abelian category with a natural functor from the category of
algebras and algebra homomorphisms.

Cyclic modules are modules of the cyclic category  $\Lambda$ which
is a small category, obtained by enriching with {\it cyclic
morphisms} the familiar {\it simplicial category} $\Delta$ of
totally ordered finite sets and increasing maps. Alternatively,
$\Lambda$ can be defined by means of its ``cyclic covering'', the
category $E \Lambda$. The latter has one object $(\Zb , n)$ for each
$n \geq 0$ and the morphisms $f : (\Zb , n) \to (\Zb , m)$ are given
by non decreasing maps $f : \Zb \to \Zb \ $, such that $ f(x+n) =
f(x)+m \qqq  x \in \Zb$. One has $\Lambda = E \Lambda / \Zb$, with
respect to the obvious action of $\Zb$ by translation. To any
algebra $\cA$ one associates a module $\cA^{\natural}$ over the category
$\Lambda$ by assigning to each $n$ the $(n+1)$ tensor power
$\cA\otimes \cA\cdots \otimes \cA$. The cyclic morphisms correspond to the
cyclic permutations of the tensors while the face and degeneracy
maps correspond to the algebra product of consecutive tensors and
the insertion of the unit. The corresponding functor $\cA\to
\cA^\natural$ gives a linearization of the category of associative
algebras and cyclic cohomology appears as a derived functor.

One can thus define the {\em distilled} module $D(\cA,\varphi)$ as
the cokernel of the cooling morphism  and consider the action of
$\R^*_+$ (obtained from the above equivariance) in the cyclic
homology group $HC_0(D(\cA,\varphi))$. As shown in \cite{CCM} this
in the simplest case of the BC-system gives a cohomological
interpretation of the above spectral realization of the zeros of the
Riemann zeta function (and of Hecke L-functions).

One striking feature is that the KMS strip (\cf Figure \ref{FigKMS})
becomes canonically identified in the process with the critical
strip of the zeta function (recall that $\beta >1$) by
multiplication by $i=\sqrt{ -1}$.

This cohomological interpretation combines with the above theory of
endomotives to give a natural action of the Galois group $G$ of
$\bar\Q/\Q$ on the above cohomology. This action factorizes  to the
abelianization $G^{ab}$ and the corresponding decomposition
according to characters of $G^{ab}$ corresponds to the spectral
realization of L-functions.

The role of the invariant $S(M)$ in the classification of factors or
of the more refined flow of weights mentioned above, is very similar
to the role of the module of local or global fields and the Brauer
theory of central simple algebras. In fact there is a striking
parallel (see \cite{CCM}) between the lattice of unramified
extensions $\K\to \,\K\otimes_{\F_q}\,\F_{q^n}$ of a global field of
characteristic $p$ and the lattice of extensions of a factor $M$ by
the crossed product algebras  $M\to \,M\rtimes_{ \sigma_T}\Z$. Using
the algebraic closure of $\F_q$ \ie the operation $\K\to
\K\otimes_{\F_q}\,\bar \F_{q}$ corresponds to passing to the dual
algebra $M\to \,M\rtimes_\sigma \R$ and the dual action corresponds
to the Frobenius automorphism when as above the appropriate
cohomological operations (distillation and $HC_0)$ are performed.

\bigskip
\bigskip
\begin{center}
\begin{tabular}{|c|c|}
\hline &  \\
 \bf{Global field $\K$} &\ \bf{Factor $M$} \\
&\\ \hline &  \\
${\rm Mod}\,\K\subset \R_+^*$  &\ ${\rm Mod}\,M\subset \R_+^*$  \\
&\\
\hline &  \\
 $\K\to \,\K\otimes_{\F_q}\,\F_{q^n}
$ &\ $M\to \,M\rtimes_{ \sigma_T}\Z$ \\
&\\ \hline &  \\ $K\to
\K\otimes_{\F_q}\,\bar \F_{q}$ &\ $M\to \,M\rtimes_\sigma \R$ \\
&\\ \hline &  \\ \bf{Points}
$C(\bar \F_{q})$ &\ $\Gamma\subset$ \bf{$X_\Q$} \\
&\\ \hline
\end{tabular}
\end{center}

\bigskip
\bigskip
\bigskip

\medskip
A notable difference with the original Hilbert space theoretic
spectral realization of \cite{CoRH} is that while in the latter case
only the critical zeros were appearing directly (the possible
non-critical ones appearing as resonances), in the cyclic homology
set-up it is more natural to use everywhere the ``rapid decay"
framework (advocated in \cite{Meyer}) so that all zeros appear on
the same footing. This eliminates the difficulty coming from the
potential non-critical zeros, so that the trace formula is much
easier to prove and reduces to the Riemann-Weil explicit formula.
However, it was not obvious how to obtain a direct geometric proof
of this formula from the $S$-local trace formula of \cite{CoRH}.
This was done in \cite{Meyer}, showing that the noncommutative
geometry framework makes it possible to give a geometric
interpretation of the Riemann-Weil explicit formula. While the
spectral side of the trace formula is given by the action on the
cyclic homology of the distilled space, the geometric side is given
as follows \cite{CCM}.

\medskip

\begin{thm} \label{geom} Let $h
\in S (C_\K)$. Then the following holds:
\begin{equation}\label{trace3}
\Trace \; \left(\vartheta(h)|_{ \cH^1} \right)=\,\wh h(0) + \wh h(1)
-\,\Delta\bullet \Delta\;h(1) -\,\sum_{v} \int_{(\K^*_v,e_{\K_v})}'
\, \frac{h(u^{-1})}{
 |1-u|}\, d^* u \, .
\end{equation}
\end{thm}

\medskip

We refer to \cite{CCM} for the detailed notations which are
essentially those of  \cite{CoRH}. The origin  of the terms in the
geometric side of the trace formula  comes from the Lefschetz
formula by Atiyah-Bott \cite{atbo} and its adaptation by
Guillemin-Sternberg (\cf \cite{gui}) to the distribution theoretic
trace for flows on manifolds, which is a variation on the theme of
\cite{atbo}. For the action of $C_\K$ on the adele class space
$X_\K$ the relevant periodic orbits on which the computation
concentrates turn out to form also the classical points of the
noncommutative space $X_\K \backslash C_\K$ distilled in the above
sense from $X_\K$. This ``classical" subspace of $X_\K\backslash
C_\K$ is given by
\begin{equation}\label{defcurve}
\Gamma_\K=\,\cup \; C_\K\,[v]\,,\quad v\in \Sigma_\K .
\end{equation}
where for each place $v\in\Sigma_\K$, one lets $[v]$ be the adele
\begin{equation}\label{base}
[v]_w=1\qqq w\neq v\,,\quad [v]_v=0\,.
\end{equation}

In the function field case, one has a {\em non-canonical}
isomorphism of the following form.

\begin{prop} \label{Zspace}
Let $\K$ be the function  field of an algebraic curve $C$ over
$\F_q$. Then the  action of the Frobenius on $Y=C(\bar \F_q)$ is
isomorphic to the action of $q^\Z$ on the quotient
$$
\Gamma_\K/C_{\K,1}.
$$
\end{prop}

\begin{center}
\begin{figure}
\includegraphics[scale=1]{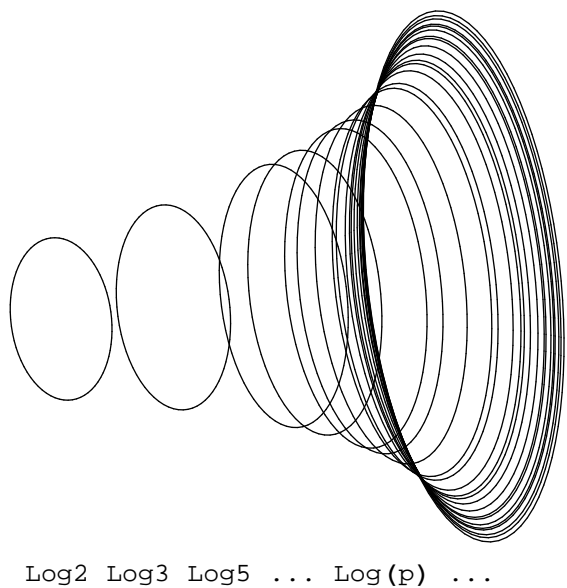}
\caption{The classical points of the adeles class space
\label{Figpoints}}
\end{figure}
\end{center}

In the case $\K=\Q$ the space $\Gamma_\Q/C_{\Q,1}$ appears as the
union of periodic orbits of period $\log p$ under the action of
$C_\Q/C_{\Q,1}\sim \R$ (\cf Figure \ref{Figpoints}). This gives a
first approximation to the sought for space $Y=C(\bar \F_q)$ in
characteristic zero. One important refinement is obtained from the
subtle nuance between the adelic description of $X_\Q$ and the finer
description in terms of  the endomotive obtained from the pointed
algebraic variety $(\bG_m(\Q),1)$. The second description keeps
track of the Galois symmetry and as in proposition \ref{Zspace} the
isomorphism of the two descriptions is non canonical.

\smallskip

At this point we have, in characteristic zero,  several of the
geometric notions which are the analogues of the ingredients of
Weil's proof and it is natural to try and imitate the steps of his
proof. The step (A), \ie the explicit formula is taken care of by
Theorem \ref{geom}. That what remains is to prove a positivity is a
well known result of A. Weil (\cf \cite{EB}) which states  that RH
is equivalent to the positivity of the distribution entering in the
explicit formulae. Thanks to the above $H^1$ obtained as the cyclic
homology of the distilled module, Weil's reformulation can be stated
as follows.

\begin{thm}  \label{pos} The following two conditions are
equivalent.
\begin{itemize}
\item All $L$-functions with Gr\"ossencharakter on $\K$ satisfy the Riemann
Hypothesis.
\item $\Trace\,\vartheta(f\,\star\,
f^\sharp)|_{ \cH^1} \, \geq \,0$, for all $f\in\cS (C_\K)$.
\end{itemize}
\end{thm}

Here we used the notation
\begin{equation}\label{conv}
f=\,f_1 \star f_2 \,,\ \ \ \text{ with } \ (f_1 \star f_2)(g)=\,\int
\,f_1(k)\, f_2(k^{-1}\,g) \,d^*g
\end{equation}
for the convolution of functions, using the multiplicative Haar
measure $d^*g$, and for the adjoint
\begin{equation}\label{adj1}
f\to f^{\sharp}\,,\quad f^{\sharp}(g)=\,\vert\, g\vert^{-1}\,\bar
f(g^{-1}) .
\end{equation}

The role of the specific correspondences used in Weil's proof of RH
in positive characteristic is played by the test functions $ f\in
\cS (C_\K) $. More precisely the scaling map which replaces $f(x)$
by $f(g^{-1} x)$ has a graph, namely the set of pairs $(x, g^{-1} x)
\in X_\K \times X_\K$, which we view as a correspondence $Z_g$.
Then, given a test function $f$ on the ideles classes, one assigns
to $f$ the linear combination
\begin{equation}\label{corresp}
Z(f)=\,\int f(g) Z_g d^*g
\end{equation}
of the above graphs, viewed as a ``divisor" on $X_\K \times X_\K$.

The analogs of the degrees $d(Z)$ and codegrees $d^{\,'}(Z)=\,d(Z')$
of correspondences in the context of Weil's proof are given, for the
degree, by
\begin{equation}\label{degree}
d(Z(h))=\,\wh h(1)=\,\int\,h(u)\, \vert u \vert  \, d^* \, u ,
\end{equation}
so that the degree $d(Z_g)$ of the correspondence $Z_g$ is equal to
$|g|$. Similarly, for the codegree one has
\begin{equation}\label{codegree}
d^{\,'}(Z(h))=\,d(Z(\bar h^\sharp))=\,\int\,h(u)\,  d^* \, u =\,\wh
h(0) ,
\end{equation}
so that the codegree $d^{\,'}(Z_g)$ of the correspondence $Z_g$ is
equal to $1$.

One of the major difficulties is to find the replacement for the
principal divisors which in Weil's proof play a key role as an ideal
in the algebra of correspondences on which the trace vanishes. At
least already one can see that there is an interesting subspace $V$
of the linear space of correspondences described above on which the
trace also vanishes. It is given by the subspace
\begin{equation}\label{V}
V\subset \cS(C_\K) \,,\quad V=\{\,g(x)=\, \sum\,\xi(k\,x)\,|\,\xi\in
\cS(\A_\K)_0\} ,
\end{equation}
where the subspace $\cS(\A_\K)_0\subset \cS(\A_\K)$ is defined by
the two boundary conditions
$$
\xi(0)=0 \,,\quad \int \,\xi(x)\,dx\,=\,0 .
$$

\medskip

\begin{lem} \label{vanish}
For any $f\in V\subset \cS (C_\K)$, one has
$$
\vartheta(f)|_{ \cH^1}=\,0 .
$$
\end{lem}

\medskip

This shows that the Weil pairing of Theorem \ref{pos} admits a huge
radical given by all functions which extend to adeles and gives
another justification for working with the above cohomology $H^1$.
  In particular one can
modify arbitrarily the degree and codegree of the correspondence
$Z(h)$ by adding to $h$ an element of the radical $V$ using a subtle
failure of Fubini's theorem. We will show in a forthcoming paper
\cite{CoCM} that several of the steps of Weil's proof can be
transposed in the framework described above.

This constitutes a clear motivation to develop noncommutative
geometry much further. One can write a very tentative form of a
dictionary from the language of algebraic geometry (in the case of
curves) and that of noncommutative geometry. The dictionary is
summarized in the following table. It should be stressed that the
main problem is to find the correct translation in the right column
(non-commutative geometry) of the well established notion of
principal divisor in the (algebraic geometry) left column. The table
below is too rough in that respect since one does not expect to be
able to work in the usual ``primary" theory which involves periodic
cyclic homology and index theorems. Instead one expects that both
the unstable cyclic homology and the finer invariants of spectral
triples arising from transgression will play an important role. Thus
the table below should be taken as a very rough first approximation,
and a motivation for developing the missing finer notions in the
right column.

\bigskip
\bigskip\bigskip\bigskip\bigskip\bigskip

\begin{small}
\begin{center}
\begin{tabular}{|c|c|}
 \hline & \\
Virtual correspondences  & bivariant $K$-theory class $\Gamma$ \\ & \\
\hline & \\
Modulo torsion & $KK(A,B\otimes {\rm II}_1)$ \\ & \\
\hline & \\
Effective correspondences & Epimorphism of $C^*$-modules \\
& \\ \hline & \\
Principal correspondences & Compact morphisms  \\
& \\ \hline &  \\
Composition & cup product in $KK$-theory \\
& \\ \hline & \\
Degree of  correspondence & Pointwise index $d(\Gamma)$ \\
& \\ \hline & \\
$\deg D(P)\geq g \Rightarrow \sim$ effective & $d(\Gamma)>0
\Rightarrow \exists K, \Gamma+K$ onto \\
& \\ \hline & \\
Adjusting the degree  & Fubini step  \\
by trivial correspondences & on the test functions \\
& \\ \hline & \\
Frobenius correspondence &  Correspondence $Z_g$ \\
& \\ \hline & \\
Lefschetz formula & bivariant Chern  of $Z(h)$ \\
&
(localization on  graph $Z(h)$) \\  & \\ \hline &  \\
Weil trace unchanged & bivariant Chern unchanged
\\  by principal divisors & by compact
perturbations
\\
 & \\ \hline
\end{tabular}
\end{center}
\end{small}

\end{document}